\documentclass[12pt,leqno]{article}
\usepackage{amsmath, amscd, amsthm, amssymb, graphics, xypic, mathrsfs, setspace, fancyhdr, bm, pdfsync, enumitem, times, endnotes}
\usepackage[french, german, english]{babel}
\usepackage[usenames, dvipsnames, svgnames, table]{xcolor}
\usepackage[letterpaper,top=1.1in, bottom=1.1in, left=1.1in, right=1.1in]{geometry}
\usepackage[colorlinks=true,pagebackref=true]{hyperref}
\usepackage{comment}
\usepackage[textsize=scriptsize]{todonotes}
\usepackage{tikz-cd}
\hypersetup{backref,hidelinks}
\usepackage{etoolbox}
\usepackage{epigraph}
\usepackage[utf8]{inputenc}

\newlength{\bibitemsep}
\newlength{\bibparskip}
\let\oldthebibliography\thebibliography
\renewcommand\thebibliography[1]{%
	\oldthebibliography{#1}%
	\setlength{\parskip}{\bibitemsep}%
	\setlength{\itemsep}{\bibparskip}%
}

\newtoggle{final}
\toggletrue{final}

\iftoggle{final} {
\newcommand{\aravind}[1]{} 
\newcommand{\tom}[1]{} 
\newcommand{\mike}[1]{} 
\newcommand{\NB}[1]{}
\newcommand{\TODO}[1]{}
\renewcommand{\todo}[1]{}
}{ 
\newcommand{\aravind}[1]{\textcolor{red}{#1}} 
\newcommand{\tom}[1]{\todo[color=blue!40]{#1}} 
\newcommand{\mike}[1]{\textcolor{green}{#1}} 
\newcommand{\NB}[1]{\todo[color=gray!40]{#1}}
\newcommand{\TODO}[1]{\todo[color=red]{#1}}
}

\newcommand{\tensor}{\otimes}

\renewcommand{\hom}{\operatorname{Hom}}

\newcommand{\cplx}{{\mathbb C}}
\newcommand{\Q}{{\mathbb Q}}
\newcommand{\Z}{{\mathbb Z}}

\newcommand{\Addresses}{{
		\bigskip
		\footnotesize
		
		A.~Asok, Department of Mathematics, University of Southern California, 3620 S.~Vermont Ave., Los Angeles, CA 90089-2532, United States; E-mail address: asok@usc.edu
}}

\newcounter{intro}
\theoremstyle{plain}
\newtheorem{theorem}{Theorem}[subsection]

\newtheorem*{claim*}{Claim} 

\newtheorem*{thm*}{Theorem}
\newtheorem*{problem*}{Problem}

\theoremstyle{definition}

\theoremstyle{remark}

\newtheorem{entry}[theorem]{}

\numberwithin{equation}{subsection}

\begin{document}
\pagestyle{fancy}
\renewcommand{\sectionmark}[1]{\markright{\thesection\ #1}}
\fancyhead{}
\fancyhead[LO,R]{\bfseries\footnotesize\thepage}
\fancyhead[LE]{\bfseries\footnotesize\rightmark}
\fancyhead[RO]{\bfseries\footnotesize\rightmark}
\chead[]{}
\cfoot[]{}
\setlength{\headheight}{1cm}

\author{Aravind Asok}

\title{{\bf Constructing projective modules} \\ \begin{large}elements of a social history\end{large}}

\date{}
\maketitle

\epigraph{Guidance by abstract ideas is a dangerous business when not controlled by strong personal relations.}{Paul Feyerabend, {\em The Tyranny of Science}}

\begin{abstract}
We discuss elements of a social history of the theory of projective modules over commutative rings.  We attempt to study the question: how did the theory of projective modules become one of ``mainstream'' focus in mathematics?  To do this, we begin in what one might call the pre-history of projective modules, describing the mathematical culture into which the notion of projective module was released.  These recollections involve four pieces: (a) analyzing aspects of the theory of fiber bundles, as it impinges on algebraic geometry, (b) understanding the rise of homological techniques in algebraic topology, (c) describing the influence of category-theoretic ideas in topology and algebra and (d) revisiting the story of the percolation of sheaf-theoretic ideas through algebraic geometry.  

We will then argue that it was this unique confluence of mathematical events that allowed projective modules to emerge as objects of central mathematical importance.  More precisely, we will first argue that, in the context of social currents of the time, projective modules initially were isolated as objects of purely technical convenience reflecting the aesthetic sensibilities of the creators of the fledgling theory of homological algebra.  Only later did they transcend this limited role to become objects of ``mainstream importance'' due to influence from the theory of algebraic fiber bundles and the theory of sheaves.  Along the way, we aim to show how strong personal ties emanating from the Bourbaki movement and its connections in mathematical centers including Paris, Princeton and Chicago were {\em essential} to the entrance, propagation and mainstream mathematical acceptance of the theory.  
\end{abstract}

\newpage
\begin{footnotesize}
\tableofcontents
\end{footnotesize}

\section{Prolegomena}
\label{s:prolegomena}
\addtoendnotes{\noindent \begin{large}{\bf Prolegomena}\end{large} \vskip .2em}
This text grew from my attempt first to justify, but then to simply analyze, the phrase:
\begin{quote}
``projective modules are important objects of study''.
\end{quote}
Initially, this justification was intended to be a small detour before a detailed presentation of the true content: applications of the Morel--Voevodsky motivic homotopy theory \cite{MV} to the problem of constructing projective modules over smooth affine algebras.  However, the more I analyzed the phrase above, the more it left me with a profound sense of unease: what did I mean by ``important''?\footnote{I use quotation marks throughout the text in several different ways, which are potentially confusing.  Of course, I use them for direct quotes.  However, I also use quotation marks in the style in which they are sometimes referred to as scare quotes; this practice serves several independent purposes: I use it to, alternatively, draw attention to the words, to indicate that I want to question standard usage, but also to indicate that the words are potentially controversial or possibly ambiguous.}  

My first attempts attacked this question like I might in any mathematics article I might write.  I aimed to justify the word ``important'' by writing a ``historical'' introduction, informed by the treatments of history that I had seen.  Typically, this involves highlighting ``important'' papers and ``ideas''--stripping away one story from the meta-narrative of mathematics--but that approach seemed circular and left me with the feeling that I was just dancing around the real question.  

While it is perhaps unsurprising to a mathematician, most of the ``historical'' treatments of mathematics I'd read were written by practicing mathematicians.  These treatments frequently seemed to just push the question further back, suggesting a long-forgotten primordial source of ``importance''.  Simultaneously, it started to seem like these treatments purposefully smoothed over many cracks, both minor and major.  In brief, these ``historical'' treatments were problematic for many reasons, some of which I will discuss later.  However, to start the discussion, let me counterpoise two statements; the tension should be evident.

\subsection{Knowledge production in mathematics; from certainty to confusion}
\label{ss:knowledgeproduction}
\addtoendnotes{\vskip .2em\noindent {\bf Knowledge production..} \vskip .2em}
As a working mathematician my view of projective modules, and the context in which they are analyzed, reflects something of what David Mumford writes in his {\em Foreword for non-mathematicians} to Carol Parikh's biography of Oscar Zariski \cite{Parikh}:
\begin{quote}
	To be a mathematician is to be an out-and-out Platonist. The more you study mathematical constructions, the more you come to believe in their objective and prior existence. Mathematicians view themselves as explorers of a unique sort, explorers who seek to discover not just one accidental world into which they happen to be born, but the universal and unalterable truths of all worlds.
\end{quote}
I write ``something of'' to highlight the {\em sentiment} of Mumford's description: the mathematical objects that we study and the ways in which we reason with them seem to have some kind of special nature.  The more we interact with them, the more they seem to take on a universal or transcendental character.  This seems especially true if one considers only ``elementary'' results like the Pythagorean theorem or the quadratic formula.\endnote{There is a lot to unpack here, depending on what your view of what mathematics ``is'' is.  Is mathematics a science?  Is it an art?  For example, if mathematics is a science, then the idea that mathematical knowledge has a special character sounds a lot like what is called ``scientism'', about which there are many strong opinions.  However, whether or not mathematics is a science, there are difficult problems. 

Unsurprisingly, many of the problems about the nature of mathematical knowledge that have been confronted appear (already?) in elementary mathematics. J.S. Mill asked \cite{Mill}, in effect: in what sense do we understand large numbers?  Much later, Frege \cite{Frege} had opinions strongly opposed to the views put forward by Mill.  I really want to separate myself from these kinds of questions about the nature of elementary mathematical objects because mathematicians frequently don't find such questions ``exciting''.  
 
Taken at face-value, Mumford's claim seems to espouse a view of what mathematicians do that I don't think any mathematician would agree with.  For example, mathematicians do not just aim to write down ``true'' statements.  For one thing, what does one mean by ``true statements"?  After G\"odel and the failure of Hilbert's axiomatic program for mathematics, one could give this meaning as ``true in some axiomatic system.''  However, I doubt there is any mathematician who would characterize mathematics as the act of writing down formally true statements in some axiomatic system, which, frankly, sounds like rather dull activity.  On this point, A.N. Whitehead famously wrote \cite{Whitehead}
	\begin{quote}
		But in the real world it is more important that a proposition be interesting than that it be true; the importance of truth is that it adds to interest.
	\end{quote}

Furthermore, so as not to inadvertently single out Mumford's view, we observe that this kind of belief is far from unique among mathematicians: one can find any number of assertions of Platonism in the mathematical literature.  Another famous one appears in Hardy's famous apology \cite[p. 123]{Hardy}: 
\begin{quote}
	I believe that mathematical reality lies outside us, that our function is to discover or observe it...the theorems which we prove, and which we describe grandiloquently as our `creations', are simply the notes of our observations.
\end{quote}
However, mathematicians are not a monolith and there are a large collection of theories of the nature of mathematical knowledge; among the most popular beyond Platonism are empiricism, formalism, and intuitionism.  Formalism in particular will play a role in our discussion because of its association with Hilbert and subsequently Bourbaki and intuitionism is frequently counterposed as its opposite in the Hilbert--Brouwer debate.  Nevertheless, as Michael Crowe discusses in \cite[\S 10]{CroweMisconceptions}, restricting one's scope to these categories is, well, limiting.  He highlights Reuben Hersch's \cite{Hersh} approach, wherein Hersch asks: 
\begin{quote}
	Do we really have to choose between a formalism that is falsified by our everyday experience, and a Platonism that postulates a mythical fairyland
	where the uncountable and the inaccessible lie waiting to be observed ... ?
\end{quote}

I mention all these things only to say that I have less interest in trying to find a precise philosophical name for the amalgam of views espoused by practicing mathematicians than in posing the question: in what (if any!) sense is mathematical knowledge universal?}  On the other hand, historian of mathematics David Rowe in his review \cite{RoweBullReview} of Jeremy Gray's book {\em Plato's Ghost} \cite{Gray} writes:
\begin{quote}
	For while mathematicians generally presume the universal validity of their work, some even ascribing to it an eternal, transcendent quality, historians are bound to the opposite prejudice: for us, mathematics is a highly stable, yet contingent body of knowledge that remains tied to its producers and practitioners. Studying the history of modern mathematics therefore requires scrutiny of the social and cultural conditions under which mathematics was made.
\end{quote}
Here, I draw attention to Rowe's qualifying use of the word {\em modern}, which, to my {\em mathematical} ears, evokes considerable distance from the elementary results described above.\endnote{Following up on the previous point, Rowe also attributes this kind of Platonism to Andr\'e Weil, who will be a central figure in our discussion, writing \cite[p. 11]{Rowe} ``the Platonist viewpoint Weil espouses has been the accepted orthodoxy for many years, not just among mathematicians...'' and it is this latter characterization that will be important for the mathematics that I discuss.  I think this is something of an oversimplification and one can find evidence for many different views of the nature of mathematical knowledge.

Yuri Manin reviewed Plato's Ghost for the Notices of the American Mathematical Society \cite{ManinPlato}, a discussion that will be reprised shortly.  As regards Mathematical Platonism, he has this to say:
\begin{quote}
	So far as I can judge, “Platonism” of working mathematicians is based on a feeling that important mathematical facts are discoveries rather than inventions.
\end{quote}
I highlight this here because the distinction between invention and discovery will appear repeatedly in what follows, but Manin continues avowing some kind of Platonism himself:
\begin{quote}
	What about Galois groups? If you feel that they were discovered by Évariste Galois, rather than invented by him, you are in a sense a Platonist. 
	
	I will call such an attitude emotional Platonism in order to stress that (in my view) it is intellectually indefensible, but not to the least degree invalidated by this fact, since our emotions happily resist rational arguments. 
	
	Being such an emotional Platonist myself, I do not want to say that all mathematics is a discovery of a Platonic world, whatever that could mean.
\end{quote}
But Manin already raises the importance of the irrational in mathematics.} What does Rowe mean by modern?  I suspect that an answer is provided by the amorphous characterization of modernism given by Jeremy Gray himself \cite[p. 1]{Gray}: 
\begin{quote}
	Here, modernism is defined as an autonomous body of ideas, having little or no outward reference, placing considerable emphasis on formal aspects of the work and maintaining a complicated—indeed, anxious—rather than a na\"ive relationship with the day-to-day world...
\end{quote}
Mathematics, especially of the pure sort, was undoubtedly in the process of moving away from the experiences and phenomena of the everyday world.  The mathematical ideas that appeared were interwoven in a meta-narrative of universality that provided meaning and gave birth to internal-to-mathematics notions of ``importance.''

In any event, placing Mumford and Rowe's statements side-by-side demonstrates a stark contrast in points of view.\endnote{As background to what I discuss, but not something about which I will say much, let me mention the ``two cultures'' drama that perhaps begins with C.P. Snow's famous ``Two cultures'' essay suggesting a disjunction between the sciences and the humanities \cite{Snow}.  The tension increased significantly in the 1970s and 80s around the bugaboo of ``relativism'' and ``post-modernism'', used intentionally vaguely here, tracing through the ``Sokal Affair'' era; books around this period predominantly written by ``science-aligned'' folks included  \cite{BigherSuperstition} and \cite{BricmontSokal}.  
	
This ``two cultures'' dichotomy is, bluntly, over-simplistic and reductive.  In the direction of texts pertaining to the ``Sokal affair'', Michael Harris wrote several things that I have found compelling including \cite{HarrisIntelligencer} and I refer to his work for some discussion of the Strong Programme as well.  Shortly after Bloor's work appeared, so did numerous criticisms; I refer the interested reader to \cite{Laudan1981-LAUITP} for criticisms from a Philosopher of science, and the second edition of {\em Knowledge and Social Imagery} begins with a response to critics. This is less interesting to me than the whole genre of books aimed at rejecting the conception of scientific or mathematical knowledge as ``contingent'' in any sense; see for example \cite{BrownSocial} for one such book written by a mathematician.

Nevertheless, methodologically, sociologists and anthropologists have ``come back from the brink'' somewhat, since the objections to relativism above were raised, especially when analyzing science.  For example, B. Latour's ``We have never been modern'' \cite{LatourWHNBM} or L Boltanski and L. Thevenet's ``On justification'' \cite{boltanskithevenet} deal with issues arising in the sciences.  While we will consult these texts, I don't believe any part of the argument hinges on such points of view.}  I will admit that when I started doing mathematics, I would have strongly endorsed something like Mumford's mathematical Platonism and I probably did explicitly do so, even if I could not have formulated it precisely.  But many years and many interactions later, I have a much harder time agreeing with Mumford's point of view, or even the sentimental version of it that I suggested.  Undoubtedly, I am not alone in this belief, and it strikes me as methodologically naive to assume that all or even most mathematicians believe in this kind of mathematical Platonism.  Thus, I mention Mumford's stated view as one of many, but one that, I aim to demonstrate, plays a considerable role in the development of projective modules.\endnote{George Lakoff and Rafael Nu{\~n}ez \cite[p. xv]{LakoffNunez} discuss a related notion of ``romance'' of mathematics, a mythology including at least the following ingredients: ``Mathematics is abstract and disembodied—yet it is real; Mathematics has an objective existence/ providing structure to this universe and any possible universe/ independent of and transcending the, existence of human beings or any beings at all; Human mathematics is just a part of abstract, transcendent mathematics; Hence, mathematical proof allows us to discover transcendent truths of the universe.''}

Describing mathematical knowledge as ``highly stable, yet contingent'' places mathematical knowledge squarely within the domain of what is typically called ``situated knowledge'' in philosophy or sociology of science.\endnote{
	It is also probably necessary at this point to make some distinction between mathematics and other sciences.  I do not want to engage with facile discussions about scientific constructionism akin to Dawkins' famous relativist at 30000 feet.  Science has clearly allowed us to exert some control over our natural environment, but I think it is unarguable that the methods that have allowed this control need not bear any resemblance to idealized versions that scientists might describe, e.g., the ``scientific method'' as it is taught to my children in school.  For the sake of context, let me just say that I have found the discussions of experimental physics in \cite{Mulkay} or \cite{KnorrCetina} elucidating.  
	
	In a mathematical direction, the discussion of Mill and Frege in \cite{Bloor} and the discussion of Hamilton's quaternions in \cite{PickeringMangle} were also very informative. Nevertheless, I am not going to assert my fealty to any particular doctrinal social constructivism; I hope the subsequent argument speaks for itself.  One related avenue of dissecting mathematical practice that I have not explored stems from questions of mathematics education.  I suspect that there is some value in formulating a ``standard model of mathematical practice'' by collecting stories we tell to those learning mathematics about how we do mathematics.}  However, practicing mathematicians routinely claim not to care about the nature of mathematical knowledge\endnote{How many times has it been claimed that the job of mathematics is just to ``prove theorems''.  More prosaically, in the introduction to {\em Foundations of Algebraic Geometry} \cite[p. vii]{WeilFoundations}, A. Weil writes:
	\begin{quote}
		At the same time, it should always be remembered that it is the duty, as it is the business, of the mathematician to prove theorems, and that this duty can never be disregarded for long without fatal effects. 
\end{quote}
Weil makes this remark in the context of the use of ``intuition'' in Italian algebraic geometry, which he defends against proponents of the ``axiomatic creed'' mentioning the importance of insight provided by careful study of examples.  As such, I am taking it out of context in using it with respect to philosophical investigation.  However, there is another plausible interpretation of this dictum to ``just prove theorems'', which I also take from Weil.  In January 1954, Weil writes to Cartan \cite{CartanWeil}:\footnote{En tout cas l'exp\'erience me semble prouver abondamment que tous ceux qui continuent \`a jouer un grand r\^ole dans le monde math\'ematiques tout en ayant cess\'e de faire des math. ou tout au moins de se tenir s\'erieusement au courant jouent t\^ot ou tard un r\^ole n\'efaste en tombant dans la basse intrique...}\begin{quote}
	In any case, experience seems to me to prove abundantly that all those who continue to play a major role in the mathematical world while having ceased to do mathematics or at least to keep themselves seriously up to date, sooner or later play a harmful role by falling into petty squabbles...\end{quote}
} while simultaneously speaking of true statements or beautiful concepts.  Reuben Hersch summarized this as: ``the typical ‘‘working mathematician’’ is a Platonist on weekdays and a formalist on Sundays'' \cite[p. 32]{Hersh}.  In the background, we should remember that debates about the nature of mathematical knowledge have been going on for centuries, so why should I as a modern mathematician secretly bristle at the idea that mathematical knowledge is situated?

What of the assertion that mathematical knowledge ``remains tied to its producers and practitioners''?  As a student, I learned mathematics from books written by people and from lectures given by people; I ended up in algebraic geometry/algebraic topology because the styles of argument and the types of objects considered in these subjects aligned with my {\em tastes}.  My graduate student cohort and I were fascinated by figures such as Évariste Galois and Alexander Grothendieck; to us, they were romantic/tragic heroes and we regaled each other with our distilled understandings of their lives.  Who could not read Grothendieck's {\em R\'ecoltes et Semailles} \cite{GrothendieckRS} as an entirely human text?  We told each other stories about the crystalline purity of the writings of Bourbaki, while simultaneously discussing the ``right'' way to think about proofs, the ``canonicality'' of constructions, or ``correctness'' of definitions.  It seems self-evident now that speaking and writing like this are a choice, though one that I was perhaps unconsciously making.\endnote{Robert Langlands, reflecting on beauty in mathematics, offers this comment on the human nature of mathematics \cite[pp. 5-6]{LanglandsBeauty}:
\begin{quote}
	There are excellent mathematicians who are
	persuaded that mathematics too is divine, in the sense that its beauties are the work of God, although they can only be discovered by men. That also seems to me too easy, but I do not have an alternative view to offer. Certainly it is the work of men, so that it is has many flaws and many deficiencies.
	
	Humans are of course only animals, with an animal’s failings, but more dangerous to
	themselves and the world. Nevertheless they have also created — and destroyed — a great deal of beauty, some small, some large, some immediate, some of enormous complexity and fully accessible to no-one. It, even in the form of pure mathematics, partakes a little of our very essence, namely its existence is, like ours, like that of the universe, in the end inexplicable.
	\end{quote}} 

Likewise, as professional mathematicians, we routinely talk about who or what {\em influenced us} to think about particular problems or to move in a particular research direction.  What we actually write in a paper might make no mention of numerous false starts or partial answers to questions that were later deemed unrelated to our framing of a narrative.  The entire content of the conversations influencing our approach and presentation in a paper might be summarized in brief and cryptic acknowledgments.  All these parts of the practice of doing mathematics are relegated to the cutting room floor.  How did this become the way we write?  To what end?  Does it make the presentation seem more ``objective'' or ``rational''?  If we do not forget all of these things, it seems impossible that personalities do not influence the transmission and perception of mathematics. 

How are we to square this with the unspoken truth that social factors should be divorced from the discussion of ``actual'' mathematics?  Sal Restivo attributes to Jean Dieudonn\'e \cite[p. 140]{RestivoSM} the statement:
\begin{quote}
	Celui qui m'expliquera pourquoi le milieu social des petites cours allemandes du XVIIIe siecle ou vivait Gauss devait inevitablement le conduire a s'occuper de la construction du polygone regulier a 17 cotes, eh bien, je lui donnerai une medaille ou chocolat.
\end{quote}
This statement reads to me as a veritable sneer in the direction of the idea that social factors could affect mathematics; but I highlight the use of the word {\em inevitablement}, which I want to analyze in more detail.\endnote{Richard Brown writes, immediately after quoting Hardy that ``mathematics is timeless'' \cite[p. 149]{BrownSocial}.  He continues by further quoting Hardy's opinion on the relative merits of Archimedes and Aeschylus (the latter will be forgotten) and then concluding that mathematics ``cannot be explained by social or historical factors.'' Brown, a mathematician with historical interests and training seems to exemplify something identified by W. Aspray and P. Kitcher \cite{AsprayKitcher} who observe that 
\begin{quote}Social histories of modern mathematics are relatively uncommon, probably because in comparison with other sciences mathematics is regarded as least affected by factors beyond its intellectual content...

Although mathematicians and historians have come to understand the value of studying professional societies, journals, prizes, institutions, funding agencies, and curricula, they have considerably less appreciation for the study of the social roots of the form and content of mathematics. This is evidence of the firmly seated belief that mathematicians but not their ideas may be affected by external factors.\end{quote}}  In an oft-quoted eassy, William Thurston was willing to admit social factors played a role in the evolution of mathematics, writing \cite{Thurston}: 
\begin{quote}
	but no matter how the process of verification plays out, it helps illustrate how mathematics evolves by rather organic psychological and social processes.
\end{quote}
But the phrase ``no matter how the process of verification plays out'' seems to downplay the importance of these features, which are rarely present in final products.\footnote{Exposition of inspiration is rare enough that it is usually memorable.  Famously, in one of the landmark papers in algebraic K-theory, R. Thomason wrote \cite{ThomasonTrobaugh}:\begin{quote}The first author must state that his coauthor and close friend, Tom Trobaugh, quite intelligent, singularly original, and inordinately generous, killed himself consequent to endogenous depression. Ninety-four days later, in my dream, Tom's simulacrum remarked, "The direct limit characterization of perfect complexes shows that they extend, just as one extends a coherent sheaf." Awaking with a start, I knew this idea had to be wrong, since some perfect complexes have a non-vanishing $K_0$ obstruction to extension. I had worked on this problem for 3 years, and saw this approach to be hopeless. But Tom's simulacrum had been so insistent, I knew he wouldn't let me sleep undisturbed until I had worked out the argument and could point to the gap. This work quickly led to the key results of this paper. To Tom, I could have explained why he must be listed as a coauthor.\end{quote}}  Why is the suggestion that social processes might {\em direct} ``pure mathematical thought'' anathema to so many mathematicians?\endnote{This is really a statement about intended audience: while I hope that sociologists of science/mathematics might be enticed to care about the modern ``pure mathematics'' I discuss, I get the impression that for reasons of ``methodological richness'' this is not the case. For example, I. Grattan-Guinness writes in his review of Restivo's book for MathSciNet:
	\begin{quote}
		The important fact is rightly emphasised that the preference for pure mathematics correlates positively with the growth of professionalisation of mathematics from the mid-19th century on, but surprisingly little is said about the sociology of modern applied mathematics. The history of mechanics and mathematical physics, for example, offers more rich issues than many of the pure cases discussed here. It would have been even more interesting to tackle areas of mathematics where sociology and/or society itself were motivating agents: the mathematics of engineering, say, where issues of economy and efficiency play a role; mathematical economics (both in the influence of algebra on the growth of Mediterranean commerce in the late Middle Ages, and of mechanics-inspired mathematical economics in the late 19th century); or the very late rise of probability and statistics as a discipline within the last 100 years, when it finally came to form a third stream apart from the disparate communities of pure and applied mathematicians.
	\end{quote}
To the pure mathematicians that remain, let me remark that there is some interest among historians/sociologists on questions of pure mathematics, and I refer them to Alma Steingart's interesting book {\em Axiomatics} \cite{Steingart}, especially Chapters 2 and 4.} 

To this discussion of ``inevitability'' and ``universality'' let me interpose the ``anxiety'' suggested by Gray.  Consider Arthur Jaffe's article \cite{JaffeProof} entitled {\em Proof and the evolution of mathematics}, which arose in the wake of fears that standards of rigor in mathematics were changing because of influence from (among other things) theoretical physics.\footnote{This episode begins with the article by Jaffe and Frank Quinn \cite{JaffeQuinn}, to which numerous mathematicians, physicists and historians responded \cite{RespJaffeQuinn}; the text \cite{Thurston} was also partial response to the Jaffe--Quinn argument.} Jaffe starts by asserting that the ``nature of what one accepts as a complete account of an argument changes over time'' and the factors affecting this change are ``human'': ``mathematics is a subject done by real people, with our own interests, our personal perspectives, our own ways of working, our pressures of life, our strengths and our shortcomings.'' Later, he writes: ``I would propose that these subtleties suggest the answer that mathematical proof is the manifestation of the highest form of human intelligence.  Moreover, mathematics is different from science because it lasts an eternity.''  Thus, social factors influence the process of mathematics, standards of proof can vary with time, but nevertheless the finished product is ``eternal''.  From the extreme phrasing of ``highest form of human intelligence'', or ``A mathematical proof has the highest degree of certainty possible for man'' one gets the sense that ``universality'' of mathematics is more tied to our self-worth as mathematicians than to mathematics itself.


In conjunction with the difficulties that many of us regularly encounter in explaining to students how to construct ``rigorous proofs'' it starts to feel like our standards of ``rigor'' and ``correctness'' are baked into our training as mathematicians.\footnote{I have overheard or myself participated in the lionizing of ``clear thinking'' and ``elegant proofs''.  We talk about the ``damage to professional reputation'' incurred by writing proofs with gaps or mistakes.  As a provocation: is this anything but reinforcing the social norms around proofs?}   I then found myself shaking my head in affirmation while reading:
\begin{quote}
	mathematical concepts and rules become self-evident only if and when they turn into “institutions”—when they are inserted in a network of concepts and practices supported by the collective interests of a social group;
\end{quote}
a summary by Massimo Mazzotti in his book {\em Reactionary Mathematics} \cite[p. 7]{Mazzotti} of the point of view of David Bloor \cite{Bloor}, one of the prime movers of the ``strong programme'' of sociology of scientific knowledge.  Imagine my surprise.

OK, so suppose you're now willing consider that social factors might have an effect on the content and direction of mathematics.  How might we tease out the influence of social factors from a literature that has been effectively scrubbed of its traces?  To answer that question, I started reading everything I could find that seemed, even vaguely, related to the theory of projective modules, beginning by rereading parts of Cartan--Eilenberg's classic {\em Homological algebra} \cite{CE} where the notion of projective module appears to be put in print for the first time.  At the same time I started writing to all of original historical figures I could, and reading the treasure trove of letters, autobiographies, and reviews that have now been collected that were unavailable years ago, e.g., \cite{CartanWeil}, \cite{GrothendieckSerre}.  Placed side-by-side, one confronts an immediate disjunction between the dry mathematical prose of the Cartan--Eilenberg text and the frankly human prose, alternatively exasperated, joking, hyperbolic,...of the correspondence/reviews/autobiographies/interviews.  

Simultaneously, I followed the references from Gray, Rowe, Mazzotti, etc.  As I kept reading, I found myself agreeing with so many things I would have found surprising twenty years ago. What part of my early acceptance of naive mathematical Platonism was due to my socialization as a mathematician, i.e., learning to take part in the professional duties of mathematics, and what part of my understanding of those professional duties depended on the time and place where I learned those ideals?  What part of this acceptance stems from a purely psychological anxiety growing from the necessity to view my direction of study and objects of fascination as worthwhile in the first place?  In the words of the critic John Guillory: ``professional training produces a certain bias of perspective, a way of seeing the world from within an occupational enclosure'' \cite[p. 3-4]{Guillory}, forming a part of what he calls a {\em d\'eformation professionelle}, a turn of phrase that very much appeals to my mathematical tastes.  

In the end, I gave up on writing a ``history'' of the theory of projective modules in the initial ``mathematical'' sense.  My original conception seemed problematic.  Furthermore, it seemed dishonest to write anything without considering the social, cultural and professional factors that surrounded the development of the ideas.

\subsection{What happens here or an outline of the argument}
\label{ss:whathappenshere}
\addtoendnotes{\vskip .2em\noindent {\bf What happens here..} \vskip .2em}
I imagine that any mathematician reading my title would envision that I was about to give the reader tools for building projective modules over commutative rings.  If this is what you thought, and you haven't already steered away, then it must appear that I have veered wildly off course.  I can now confess that my title pays homage to both Andrew Pickering's fantastic book {\em Constructing Quarks} \cite{Pickering} subtitled ``a sociological history of particle physics'' and a chapter from John Milnor's treatise {\em Introduction to Algebraic K-theory} \cite[\S 2]{MilnorKtheory} which bears the same name.  

Pickering's book is nominally about the developments related to the standard model in particle physics, focusing on ``practice'' in high energy physics.  Pickering later studied ``practice'' in a more purely mathematical context with his analysis of Hamilton's invention of quaternions in \cite[Chapter 4]{PickeringMangle}.  To a mathematician, this notion of practice likely requires comment.  Pickering differentiates two notions of practice; I will summarize them here as I will use them in mathematics: a ``generic'' one, in this context, the work of the production of mathematical knowledge so to speak and a ``specific'' one, the specific (sequences of) activities that mathematicians rely on in their work.  Regarding the latter, Pickering underscores that it ``falls into the sphere of culture'' or, said differently, practice is grounded in culture.  In {\em Constructing Quarks}, Pickering writes about high energy physics that ``the emphasis is on practice, and the practice is irredeemably esoteric''.  These words strike me as applying equally well to the practice of mathematics, as I hope some of the comments above have suggested.   

By contrast, Milnor's name is routinely trotted out when discussion turns to examples of great mathematical exposition.  What supports this opinion?  I, for example, routinely recommend these books to students as primary texts for directed reading courses.  Just as routinely, they find them difficult and not enlightening in the way that I have internalized they should expect.  Are they doing something wrong, perhaps simply ``reading'' the text instead of ``actively reconstructing'' it?\footnote{In the off chance you are not a mathematician, let me add a word about what I mean.  One of the first things I tell people taking advanced mathematics courses is that mathematics textbooks are not meant to be read linearly, even if they are constructed to reflect a particular progression of ideas.  Instead, one has to develop a feeling for the notions involved by a process akin to experimentation.  Outside of doing exercises laid out in a text, here are some possible strategies. Read a definition, then shut the book and try to deduce logical consequences.  Pick a ``key'' example, and understand how features of a definition are abstracted from the example or, understand how features of a theorem are reflected in the example.  Pick a theorem and then try to ``break'' the proof: understand how the hypotheses of the theorem are used in the proof and what happens if they are dropped.  This process is very active and involves engaging with a great many ideas lying around the text.}  Writing style is a matter of making choices.  I can now confess that even if the mathematics is presented in a way that I appreciate now, these style of these texts left me with the ineffable sense that ``this is not mathematics {\em I} can or will make''.  I am thus led to ask: why do Milnor's stylistic choices reflect the values of contemporary mathematicians?  

The year 1950 will become something of a focal-point in our discussion for a few reasons: it was marked by two events that represent themes I'd like to explore further.  The International Congress of Mathematicians (ICM) was held in Cambridge, MA and, simultaneously, the National Science Foundation was established.  These two events can be viewed as lenses refracting the changing nature of mathematical practice.  The first aligns with what I feel is a turning point in the professionalization of mathematics, a topic that I will explore during the text in many ways.  The second, which will appear only tangentially throughout the text, but looms large in the background, is the role of research funding in charting the course of mathematics.

At the 1950 ICM, Raymond Wilder, a topologist at the University of Michigan of some repute, trained in the E.H. Moore school, delivered a plenary lecture on {\em The cultural basis of mathematics}, which eventually became \cite{WilderICM1950}, a short note which he describes as growing out of his interactions with anthropologists. For Wilder, ``culture'' was ``the collection of customs, rituals, beliefs, tools, mores, etc., which we may call cultural elements, possessed by a group of people''.  Granting this, he wrote
\begin{quote}
	As mathematicians, we share a certain portion of our cultures which is called ``mathematical.'' We are influenced by it, and in turn we influence it. As individuals we assimilate parts of it, our contacts with it being through teachers, journals, books, meetings such as this, and our colleagues. We contribute to its growth the results of our individual syntheses of the portions that we have assimilated.
\end{quote}
Even though Wilder is writing about mathematics in general, it is to me telling that he speaks about a plurality of culture{\em s}.  

In light of my assertion of the collective avoidance of discussion of socio-cultural factors in contemporary mathematics, a {\em mainstream} mathematician talking about mathematical culture in a prominent international setting seems extremely noteworthy to me.  Moreover, in this era, Wilder was not even alone among mainstream mathematicians in considering the importance of ideas of social science on the development of mathematics: Dirk Struik had written a call for a sociology of mathematics by 1942 \cite{Struik}.  I have already parenthetically noted some of the enmities that have arisen between sociologists of mathematics and mathematicians, but I think it's worthwhile to ask: what has changed since then?

More recently, Mazzotti describes ``mathematical cultures'' in the following terms \cite[p. 7]{Mazzotti}:
\begin{quote}
	There is a long tradition of identifying styles in mathematics, a notion that has been variously understood as having to do with individual temperament, national character, the Zeitgeist, or the direct expression of a culture in a Spenglerian sense.
\end{quote}
Mazzotti refers to work of Rowe \cite{Rowe} for a view of ``mathematical culture'' that echoes Wilder: ``the set of resources available to a specific group of mathematicians''.\footnote{My hope is that ``culture'' used in this latter sense is uncontroversial.  For example, I think most professional mathematicians of some level of seniority would agree that any given field of mathematics has a culture, at least in the following limited professional sense: practitioners are familiar with major results, commonly used tools, standard methods of argument, etc., even if it is difficult to demarcate the boundaries of a field.  Whether or not two or more mathematicians, even in the same field, would agree precisely on what constitutes the ``culture of a field'' naturally leads to a pluralistic notion of culture.}  Bearing this in mind, it seems imperative to describe the {\em mathematical} cultures in which those involved with the theory of projective modules at its inception were enmeshed.\endnote{The plural cultures, used in Wilder's original sense, and perhaps more generally in the sense used by anthropologists of that era has also been criticized, even in more limited senses in which it may be applied to the sciences; see \cite{ChemlaKeller} for some discussion of this point, but I highlight one point that seems quite relevant to my discussion \cite[p. 8]{ChemlaKeller}: 
	\begin{quote}
		as the outcome of activity, a culture of scientific practice is subject to constant change—in relation to the problems actors address and the goals they pursue and in the ways they draw on the resources available to
		them to mold and remold their objects of research, their values, and so forth.
		Finally, we suggest that establishing bridges between cultures of scientific practice is also part of what actors do. Overcoming differences in knowledge and practice, constructing sameness (or even universality), and achieving consensus are not properties of scientific practice and knowledge that are given a priori but are outcomes of actors’ knowledge activities.
	\end{quote}
I will touch on this in the main body of the text, but in this terminology what will be important is actors perceptions of their own mathematical culture(s).	
	
More broadly, within spheres of scientific production one may appeal to different notions of culture.  For example, one may consider {\em epistemic cultures}; Knorr Cetina in \cite[p. 363]{KnorrCetina2005} defines these as follows:
\begin{quote}
	The notion of epistemic culture is designed to capture these interiorised processes of knowledge creation. It refers to those sets of practices, arrangements and mechanisms bound together by necessity, affinity and historical coincidence which, in a given area of professional expertise, make up how we know what we know. Epistemic cultures are cultures of creating and warranting knowledge.
	\end{quote}
Knorr Cetina's book \cite{KnorrCetina} produces, to my mathematical eyes, interesting analyses of several epistemic cultures around experimental physics and biology.  This notion of epistemic culture is, I believe, intended to be more narrowly applied; in the context of this document, one could imagine applying this notion to styles around the Bourbaki group itself.  This ``micro'' notion is then contrasted with a much broader ``macro'' counterpart called a {\em knowledge culture} \cite[p. 369]{KnorrCetina2005}.  For us I would argue that the ``macro'' conception means mathematics on the scale of international mathematics, broadly construed.}  

This all leads me to return to the question of the meaning of the word ``important'' in my initial premise.  I aim to argue that ``important'' mathematics is a dynamic and emergent concept within given mathematical cultures.  On the one hand, historically, mathematicians have aimed to do ``important'' mathematics.\footnote{Readers who have just been reminded of Gauss's construction of the 17gon will now turn back to Dieudonn\'e and ask how Gauss was trying to do ``important'' mathematics in this construction.  Here, I enquote ``important'' because it is being used vaguely for a constellation of valuative judgements for which a number of other words might be suitable, e.g., ``interesting'', ``good'', etc. One might try to distinguish ``important'' from ``interesting'' by suggesting ``important'' has an outward focus, e.g., it makes sense within a culture, whereas ``interesting'' has an inward perhaps more individual focus.  However, the social nature of the process of demarcating ``important'' or ``interesting'' ideas I propose to use blurs the distinction.}  For pre-modern mathematics, in the sense of Gray, external validation made it easier to describe importance, but the modernist turn gave rise to a subtextual anxiety. Unsurprisingly, mathematicians have even made efforts to characterize this and related concepts, for instance one can look at Saunders MacLane's article \cite{MacLaneExcellence} or Terry Tao's article \cite{Tao}.  More recently, Akshay Venkatesh \cite[\S 3]{Venkatesh} spent some time analyzing this concept in relation to the potential impact of general artificial intelligence on the course of mathematical research.  

For a long time I believed that importance of mathematical ideas had some ``intrinsic'' nature.  I had internalized a narrative that the universality of mathematics we work with was even {\em inevitable} in a certain sense.\endnote{Again, I want to distinguish {\em elementary} mathematics from modern mathematics here.  One might argue that the fact different basic mathematical notions arose in cultures with vast geographic separation is a testament to the inevitability of these mathematical ideas/structures.  At this level, I suggest that something akin to the idea of {\em parallel evolution} in biology actually suffices to explain those instances.  While I don't want to develop this argument here, distinctive features of the human brain and natural evolutionary pressures, e.g., around the development of agriculture seem sufficient to explain the rise of arithmetic in different cultures.  
	
Such a conception is certainly not new.  Writing in \cite[p. 425]{WilderGrowth} Wilder in reference to J.L. Coolidge's pronouncement that it was ``a curious fact in the history of mathematics that discoveries of the greatest importance were made simultaneously'' that there was nothing ``curious'' about it, nor was it confined to mathematics.  Nevertheless, even much later Mac Lane specifically highlights {\em inevitability} as a marker for excellence in mathematics \cite{MacLaneExcellence}.}  My feeling is that some kind of universality is a widely held belief among mathematicians,\endnote{This feeling is, admittedly, entirely based on anecdotal evidence and I should qualify it slightly.  Robert Wagner \cite[\S 2]{Wagner} writes:
\begin{quote}
Beauty, originality, importance and other measures of quality are not, I believe, more consensual among mathematicians than among other scientists. In fact, and contrary to common perceptions, even truth is not subject to consensus among mathematicians. Indeed, mathematicians have different conceptions of mathematical truth, even if these conceptions are not always as well-articulated as those developed by logicians and philosophers... 
\end{quote}
He does not offer any support for this belief either.  Note, however, that Inglis and Aberdein suggest that individual mathematicians {\em are} consistent in their aesthetic evaluations \cite{Inglis2015-INGBIN}, whereas aesthetic assessments of proofs between mathematicians {\em can} vary \cite{Inglis2016-INGDIP}.

Wagner does observe that mathematicians ``tend to agree much more than experts in other domains on whether a given argument or proof is valid.''  It something like this latter agreement, echoed in Jaffe's comments earlier on the role of proof that I imagine when speaking about universality.} perhaps stemming from the way in which we interact with modern mathematics itself.  To state things in a strong way: if history were to repeat itself, we would eventually come up with the same definitions and ``important'' definitions/theorems.\endnote{The careful reader will not have missed my Easter-egg use of the word ``meta-narrative'' early on.  The ``universality'' I just described goes to the heart of modernism as characterized by Lyotard \cite[p. xxii]{LyotardPostmodern}:
	\begin{quote}
		I will use the term {\em modern} to designate any science that legitimates itself with
		reference to a metadiscourse of this kind making an explicit appeal to some grand narrative, such as the dialectics of Spirit, the hermeneutics of meaning, the emancipation of the rational or working subject, or the creation of wealth. For example, the rule of consensus between the sender and addressee of a statement with truth-value is deemed acceptable if it is cast in terms of a possible unanimity between rational minds: this is the Enlightenment narrative...
	\end{quote}
Lyotard furthermore echoes my discomfort with universality and inevitability \cite[p. 54]{LyotardPostmodern}: 
\begin{quote}
	But what never fails to come and come again, with every new theory, new hypothesis, new statement, or new observation, is the question of legitimacy.
\end{quote}
Once again, it is the sentiment that I want to highlight, rather than any details of the argument.  Indeed, Lyotard refers to this book \cite[p. 26]{AndersonOriginsofPostmodernity} as the absolute worst of his books.  Moreover, his use of the notion of meta-narrative, especially in science, has been criticized \cite[p. 94]{bertens1997international}, e.g., for being self-refuting in a ``Russell paradox''-like way: that there are no meta-narratives is itself a meta-narrative.}  Reliance on {\em intrinsic importance} is still present in sociologically-inflected treatments from as late as the 1970s, e.g., Wilder, writing in \cite{WilderStress} highlights {\em capacity} as contributing to what he calls {\em hereditary stress}, and capacity is defined as ``The quantity and intrinsic interest of the results that the basic theory and methodology of a field are capable of yielding.''  There are echoes of this ``intrinsic importance'' also, for example, in Wigner's unreasonable effectiveness essay \cite{Wigner}.  As we will see, however, this sense of inevitability seems also to be an {\em a posteriori} construction, perhaps due to the ways in which we expose mathematics, and is largely invisible in mathematical practice.

At the Fields Medal symposium held in 2022, Tim Gowers studied a related question in his lecture entitled {\em Is mathematical interest just a matter of taste?} \cite{Gowers}.  In particular, Gowers claims that ``the extent to which a statement is interesting is related to objective features of how it fits into the corpus of mathematical knowledge.''  However, as pure mathematics has evolved away from its sources of inspiration, my point of view of importance has shifted and echoes Herrnstein-Smith \cite[p. 30]{SmithContingencies} who argues:
\begin{quote}
	All value is radically contingent, being neither a fixed attribute, an inherent quality, or an objective property of things but, rather, an effect of multiple, continuously changing, and continuously interacting variables.
\end{quote}
In my conception, ``important mathematics'' does not have any kind of objective or universal characterization. Instead, it arises organically from the ideas upon which  mathematical cultures choose their focus: we don't know it until {\em after} we see it, so to speak.  Moreover, this {\em choice} of focus evolves subject to myriad mechanisms including random chance, status hierarchies, personal aesthetics, historical narratives, etc. (some of which we'll highlight in our analysis proper).  How this choice plays out in real-time seems reflective of Pickering's ``mangle'' \cite{PickeringMangle}, which I will try to illustrate by way of a case-study.  

To make this case, I will first discuss the ``mathematical cultures'' in which the notion of projective module emerged, which will involve several interdependent narratives.  I will try to do this in a way that pays heed to Herbert Mehrtens's characterization of mathematical historiography \cite[pp. 258-259]{Mehrtens}:
\begin{quote}
	The mathematician's understanding of the history of his subject is frequently sharply anti-sociological. But as frequently it is a 'rational reconstruction' of a development governed by universal laws - those of mathematics itself.  In the unfolding of eternal mathematical laws through time there is no room for stories. They are only secondary flavouring, anecdotes of the lives of great men.\footnote{I wonder if Mehrtens has Andr\'e Weil's famous letter to his sister \cite[p. 244]{WeilI} in mind: \begin{quote}Je t'avertis que tout ce qui concerne l'histoire des math\'ematiques, dans ce qui suit, repose sur une \'erudition tout \`a fait insuffisante, que e'est pour une bonne part une reconstitution a priori, et que, m\^eme si e'est ainsi que les choses ont d\^u \^etre (ce qui n'est pas prouv\'e), je ne saurais affirmer que e'est ainsi qu'elles ont \'et\'e. En math\'ematiques, d'ailleurs, presque autant qu'en toute autre chose, la ligne de l'histoire a des tournants. 
	\end{quote}
	Moreover, as we will discuss, this view of history does seem to prevalent among we might call Bourbaki-aligned historical treatments.}
\end{quote}
Merhtens's comment does appear appropriate for the episodes that we will access.

A precise definition of the notion of projective module first appears in the landmark treatise of Henri Cartan and Samuel Eilenberg entitled {\em Homological Algebra} \cite{CE}, which retrospectively created the field of homological algebra out of whole cloth.  I aim to recreate some aspects of the mathematical scene prefiguring the emergence of the Cartan--Eilenberg text.  Why was this text written?  Who was thinking about what?  

I will focus my attention around several key developments, all originating in roughly the period 1920-1950. 
\begin{itemize}
\item[$\bullet$] {\em Fiber bundles} (See Section~\ref{ss:fiberbundles}), initially appearing as special examples of the then new notion of topological space studied by Seifert, penetrated into complex analysis and algebraic geometry through the efforts of Andr\'e Weil and Henri Cartan; these ideas percolated through and were internalized by many members of the famous Cartan seminar. Their use appeared, to Weil at least, as a unifying theme.
\item[$\bullet$] {\em Homological algebra} (See Section~\ref{ss:homologicalalgebra}) evolved from Poincar\'e's work on homology.  A proliferation of different approaches to homology were a challenge for the modernist mathematician in the 1940s leading to the axiomatic turn in the work of Eilenberg and Steenrod.  The transplant of the ``homological'' to other algebraic contexts, e.g., groups, associative algebras and Lie algebras presented yet another challenge for modernism.  
\item[$\bullet$]{\em Category theory} (See Section~\ref{ss:categories}) was a tool developed to understand the relations between formal structures appearing in topology and algebra.  Describing it as ``abstract nonsense'' betrays anxiety in the explicit use of its concepts, but the ``functorial point of view'' as invented by Eilenberg and Mac Lane created a framework for, among other things, transplanting analogies.
\item[$\bullet$]{\em Sheaf theory} (See Section~\ref{ss:sheaves}), invented by Leray, encapsulates the passage from ``local'' information to ``global'' information in topology.  Sheaf theory also emerged from homological investigations, this time taking inspiration from the work of Georges de Rham and presented an approach that wasn't easily amenable to the Eilenberg--Steenrod axiomatization.  Reconciling these two approaches to homology led beyond the work of Cartan--Eilenberg.
\end{itemize}

This selection may seem somewhat arbitrary and idiosyncratic, but I believe this cultural background presents a solid cross-section of the mathematics to which people in the Cartan--Eilenberg circle were exposed and serves to demarcate the ``mathematical mainstream'' of the era.  Granted this cultural background, I revisit the development of the notion of a projective module.  After describing the general reception of Bourbaki's algebra in Section~\ref{ss:bourbakialgebra}, which aims to situate how the Cartan--Eilenberg text may have been received by those just outside of their circle.  Then, Section \ref{ss:axiomaticprojectivity} aims to demonstrate that, in their initial conception, projective modules were mainly of ``technical'' interest within homological algebra.  More precisely, their very definition is motivated by what I will call a category-theoretic (or just categorical) aesthetic, as opposed to more mundane reasons for making definitions in mathematics.  More strongly, I argue that without the other influences mentioned above, projective modules would have remained an isolated ``technical'' tool, and likely would {\em not} become ``important'', in at least the limited sense that they were a topic of frequent consideration by influential mainstream mathematicians.  

Section~\ref{ss:fiberbundlesandprojectivemodules} then argues, in the context of the ideas circulating in the Cartan seminar, that the influence of the theory of fiber bundles and sheaf theory breathed new life into the notion of projective module: work of J.-P. Serre laid the groundwork for analogical use of the topological theory of fiber bundles, which he developed somewhat, but was taken up by many mathematicians shortly thereafter.  Section~\ref{ss:conclusion} argues that with some additional influences, e.g., from algebraic geometry via Grothendieck's proof of the Riemann--Roch theorem, the theory stood, near the end of the 1950s as a rich source of problems and questions giving rise to a more self-sustaining domain of mathematics.  This choice of end-point is somewhat arbitrary, but I hope sufficient to make my point.

The reader not interested in the influence of social factors in mathematics has probably long ago stopped reading this introduction.  However, for those that remain, the main body of the text aims to develop the narrative just expounded with an aim toward social context surrounding the invention of various mathematical notions making up ``mathematical cultures''.\footnote{At this point, the reader probably cannot help but to have observed that I have used the word {\em invented} rather than {\em discovered} several times; see the introduction to \cite{Hadamard} for a discussion of the difference between these two notions written around the time period under discussion.}  Having highlighted the arbitrary choices made in isolating the definition of projective module, I hope to emphasize that the modern notion was {\em in no sense} inevitable.  Rather, it was contingent upon interactions of specific mathematicians with particular styles of doing mathematics.  It depended on the particular climate of mathematics around the Cartan seminar, in conjunction with the general, ``modernist/axiomatic'' influence of Bourbaki ideology on French mathematics at the time, as well as more specifically the influence of Andr\'e Weil and Jean-Pierre Serre.  
The initial definition would not have grown to be ``important'' without constant influx of interest from what I call ``high-status'' mathematicians and mathematical domains, which themselves were sustained and elevated by the mathematical culture of the time, encompassing mathematics outside of Europe, but especially in the United States and centered at institutions including Princeton, Columbia and the University of Chicago. 

\subsection{Why now? On self-reflection in contemporary mathematics}
\label{ss:whynow}
\addtoendnotes{\vskip .2em\noindent {\bf Why now..} \vskip .2em}
But why now?  Is it simply reaching middle-age in mathematics that leads one to grapple with ``historical'' concerns (is this the mathematical equivalent of a mid-life crisis?).  Grappling with the sources of progress in mathematics (or science more generally) seemed a more common activity several generations ago.  Indeed, one need only look at the works of Poincar\'e, Hilbert, or Weyl who wrote at length on questions of value in science.  

Einstein, in his obituary for Mach, famously wrote (this translation is taken from \cite{Holton})
\begin{quote}
	Concepts which have proved useful for ordering things easily assume so great an authority over us, that we forget their terrestrial origin and accept them as unalterable facts. They then become labeled as ‘conceptual necessities,’ 'a priori situations,' etc. The road of scientific progress is frequently blocked for long periods by such errors.
\end{quote}
This leads into the oft quoted ``not an idle game'' motif, but also from our point of view to the question: where has modernist mathematics left us as regards progress in mathematics?\endnote{One investigation of this kind of question in a scientific context can be found in \cite{laudan1978progress} where on p. 3 he writes:
\begin{quote}
Some historians and philosophers of science (e.g., Kuhn and Feyerabend) have argued, not merely that certain decisions between theories in science {\em have been irrational}, but that choices between competing scientific theories, in the nature of the case, {\em must be irrational}.  They (especially Kuhn) have also suggested that every gain in our knowledge is accompanied by attendant losses, so that it is impossible to ascertain when, or even whether, we are progressing.
\end{quote}
I see no reason why something similar could not also apply to mathematics.  For that reason, it seems worth considering Laudan's reorientation of progress in science toward ``solutions of problems''.  As a {\em human} conception, widespread interest in a problem seems a good justification of worth, even if it can only be local and social agreement.  Undoubtedly we already use ``problem-solving effectiveness'' as a proxy for mathematical interest.  Nevertheless, I hope to indicate in the main body of the text that such a discussion is too simplistic for progress in a mathematical context; see especially Paragraph~\ref{par:pointsofview}}  As regards this text, I hope that a  mathematician will find in Section~\ref{s:culture} an exegesis\footnote{This word seems to me particularly apt in the context of algebraic geometry in light of the recurring use of religious terminology and references.  One need look no further than the introduction to \cite{SGA3.1} where the S\'eminaire Chevalley is literally referred to as [Bible].} of some ideas that have become, in fact, ``conceptual necessities'' in modern algebraic geometry/algebraic topology.  

Looking back at the era of Poincar\'e, Hilbert, Weyl, or Einstein, my feeling is that of a time that was extremely turbulent, but marked by an (Enlightenment?) ethos of universality and unification in science.  The current era seems, by contrast, to be marked by rising walls of disciplinary boundaries, increasing professional specialization, compartmentalization and exhortations to not overextend one's domain of expertise.  In contrast, the era which I approach in this text seems to be marked by forgetting and increasing fragmentation.  Here, I mean forgetting in the positive sense of Lewis Hyde's memoir {\em A primer of forgetting}, \endnote{Hyde, writes \cite[p. 264]{hyde2019primer} 
\begin{quote}	
	Myself, when writing poems, I practice revision by forgetting.  I write a draft of the poem, and then another and another, allowing the versions to pile up in a jumble--lines I am attached to, although they don't belong, lines that fit but go flat in the middle, words replaced and then reinserted, promising developments that never delivered--it all sits there, a shapeless pile, clammy with fatigue.  
	
	Then I set the mess aside and ignore it for at least one day.  Then I write the poem from memory.  Great chunks will have fallen into oblivion, while others will have returned clarified from the pool.  The double goddess attends, erasing as she records, drawing shape from shapelessness, dropping the discord to reveal the harmony.
\end{quote}

This practice seems to echo the perhaps apocryphal method of Jean-Pierre Serre for writing papers, which I have heard referred to as the $n+1$-method: write a draft of a paper, put it in a drawer.  Write another draft: if the two coincide, stop.  If they do not, throw away the first draft and repeat the process until it converges.  More abstractly, this echoes the description of the Bourbaki axiomatic approach: one aims to strip away the inessential.  My main question about all of this is: what factors influence how one decides what is inessential?  This seems just as fraught as deciding what is essential!  The conception of ``Losing knowledge'' as part of ``progress'' also appears in Peter Burke's sociology \cite[Chapter 5]{burke2012social}, but I will not develop this here.} which I will come back to momentarily.  

However one wants to characterize the contemporary era of mathematics, it is my feeling that self-reflection about questions of value and progress in mathematics has been pushed to the background, and is (arguably) even scorned.  For an indication of the latter, one look no further than Hardy's {\em A mathematician's apology} \cite{Hardy}.  In response to this, Hermann Weyl writes \cite[p. 163]{WeylMindandNature}:
\begin{quote}
	If I view the situation in which I find myself tonight a little less melancholically, it is not because I disagree with Hardy in that “mathematics is a young man’s game,” but because I do not quite share his scorn “of the man who makes for the man who explains.” It seems to me that in mathematics, as in all intellectual endeavors, both things are essential: the deed, the actual construction, on the one side; the reflection on what it means on the other. Creative construction unguarded by reflection is in danger of losing its way, while unbridled reflection is in danger of losing its substance.
\end{quote}
Weyl's statement also seems to reflect the sentiment contained in my epigraph from Paul Feyerabend's {\em Tyranny of Science} \cite{feyerabend2011tyranny} though focuses less on who will do the reflecting.  Another answer to the ``why now?'' question, for which this text might serve as background is an attempt to answer the question: what, if anything, has changed from the Poincar\'e era to the present day about the nature of self-reflection in mathematics?   


That being said, it is clear that reflection in general does exist, and I have discussed some examples in the preceding section.  To my mind, modern examples, especially as related to my field of interest, seem to have a methodological focus, e.g., concerns around the nature of proof.  This focus has some hallmarks of what is called ``boundary work'' in the sociology or history of science: how does one demarcate science from non-science?  Clear from the discussions above is that our conception of proof is involved in the process of demarcating mathematics from non-mathematics.  Thinking of the ``map'' of mathematics, one can take up Thomas Gieryn's cartographic analogy \cite{gieryn1999cultural}: where is the boundary between proof and non-proof and how do we draw the line?  Perhaps it was more widely believed 50 years ago that this line was easy to draw.  However, the complexity of modern proofs, sometimes distributed among thousands of pages written by many different authors makes finding the borders of the mathematical continents murkier.\endnote{One key example of what I mean here that has received prominent attention is the classification of finite simple groups.  Michael Aschbacher, a prime contributor to the proof, remarks on related questions in \cite{Aschbacher}.  He writes ``because of the complexity of the proof and the absence of a definitive treatment in the literature, one can ask if the theorem has really been proved.  After all, the probability of an error in the proof is one. Indeed, presumably any attempt to write down a proof of such a theorem must contain many mistakes.''  He then states ``I suspect most professional mathematicians feel that, after some (high) minimal standard of rigor has been met, it is more important that the proof convey understanding than that all formal details appear without error.''  This instance seems to be methodologically interesting to historians as well.  Alma Steingart analyzed this proof as well \cite{SteingartGroups}, writing ``Here, tacit knowledge and personal communication were not glossed over in publications, but instead were indispensable components of how mathematicians were able to approach, apprehend, and evaluate an unwieldy body of literature.

While group theory provides an example that is easy to wrap one's head around, modern mathematical domains frequently operate in a gray area as regards proof.  Examples of this kind of phenomenon range from use of statements made in talks for which complete write-ups are not available, to ``known-to-experts'' errors in published results.  In the latter case, sometimes fixes are known (to varying degrees!), but sometimes they are not, and mathematicians can operate under the belief that the proofs will {\em eventually} be fixed.  Numerous examples of some provisional form of truth are attached to episodes in the mathematical literature, ranging from the gap in Wiles' original proof of the modularity theorem, several episodes in symplectic geometry, through Mochizuki's claimed proof of the ABC conjecture.  My point is not to highlight errors, but simply to observe that vibrant fields can operate without any pretense to absolute certainty.  For a recent treatment of what one might call an ``epistemology of mathematical error'' I refer the reader to \cite{WeatherallWolfson}.}  Moreover, proofs for a paper written in a given domain frequently rely on expected cultural knowledge, leading to the idea, espoused by Andrew Granville that proof itself is a social compact \cite{Granville}.

Related directions include computer-aided proof, say around the four color theorem, or the Kepler conjecture, the nature of proof in the wake of the influence of external factors, e.g., physics, say as discussed in the Jaffe--Quinn article and its responses, or the nature of proof in the wake of large ``programmatic'' approaches to mathematics.  In the latter category I place the reflections of V. Voevodsky on proof-assistants\endnote{As I mentioned at the outset, this was supposed to be an introduction to applications of the Morel--Voevodsky motivic homotopy theory.  Voevodsky played a central role in the establishment of this theory, culminating in his proof of several conjectures including Milnor's conjecture on the mod $2$ norm residue isomorphism, and its generalization to other primes called the Bloch--Kato conjecture.  In the course of these proofs, Voevodsky developed his approach to ``motivic cohomology'', fulfilling part of a conjectural vision of Beilinson and Lichtenbaum.  
	
	In the course of lecturing on the Bloch--Kato conjecture, Voevodsky writes \cite{Voevodsky}: ``Only then did I discover that the proof of a key lemma in my paper contained a mistake and that the lemma, as stated, could not be salvaged.  Fortunately, I was able to prove a weaker and more complicated lemma, which turned out to be sufficient for all applications. A corrected sequence of arguments was published in 2006.''   Coincidentally, I attended these lectures whilst a graduate student, and they were formative in shaping my mathematical tastes.  This anxiety led Voeovdsky to his Univalent Foundations project and numerous discussions on ``certainty'' of proof in mathematics.} and formalization and the discussions of Kevin Buzzard.\todo{Jesse Wolfson reference.}  Of these discussions, only those of William Thurston and Michael Harris' memoir \cite{HarrisMWO} grapple with the nature of mathematical progress.  These are important questions, and I think having a clearer sense of how social factors affect evolution of mathematical ideas would be a useful tool in their analysis as well, but it is not my intention to address such issues here.



A second reason for ``why now'' stems from a certain reflexivity of method.  If we agree that mathematics is a human activity, then exposing the human aspects explicitly simply makes a better story, but then so too should analyzing the human aspects of my construction of the narrative.  Along the way I'd like to argue that this {\em human} aspect of mathematics has never really disappeared, but rather has been relegated to a background role.  In this brand of algebra/algebraic geometry/algebraic topology, the style of exposition that is highly valued, highlighted by, for example, Milnor's texts, the Cartan--Eilenberg text, or the numerous texts of Weil and Serre, has only served to increased the socially dependent nature of those texts.  The human and social aspect of mathematics plays out in private, in common-room conversations, letters, e-mails, snide remarks at seminar talks, apocryphal stories, but has been faithfully excised from much mathematical exposition and ``finished products''.  I claim that bringing it to the fore enriches our understanding.

To close, I want to highlight one final aspect of ``why now'' related to the notion of progress in mathematical knowledge.  To frame this discussion, let me recall what Yuri Manin writes about the ``modernist transformation'' of mathematics in his review of {\em Plato's Ghost} \cite{ManinPlato}:
\begin{quote}
	 ...“the modernist transformation of mathematics”...was one such periodically occurring refurbishing of basic vocabulary, grammar, and, yes, aesthetic requirements for mathematical thought and mathematical texts. 
\end{quote}
Manin simultaneously asserts that ``one remarkable feature of mathematical knowledge is this: we learn more and more about the same objects that ancient mathematicians already started to see with their mental eyes: integers and prime numbers, real numbers, polynomial equations in one or many variables, space and various space forms...''   

I'd like to propose a different view, somewhat antithetical to Manin's point of view and the cumulative view of mathematical knowledge.  We don't read mathematics so much as teach ourselves to perform it; we don't memorize proofs so much as teach ourselves to reconstruct them; we reconstruct ideas in our brain.  The ``refurbishing of basic vocabulary and grammar'' means that, in a sense, each generation invents/constructs mathematics anew, with all of the cultural change that portends.  For notions that are near to lived experiences or that can be shared by many people (say, integers and prime numbers), like some of the notions Manin describes, the resulting changes can seem cumulative, but for notions that are further away (say, our notion of space) the changes can be rather radical.\endnote{For a mathematical analogy: the micro-changes I am suggesting are analogous to analytic continuation and are to be distinguished from ``revolutions'' or ``paradigm shifts'' in any kind of K\"uhninan sense.}  Concretely: what (if any!) of modern arithmetic geometry would be recognizable to Dedekind?\todo{This is not an idle question: Bhatwadekar story.  How does the older generation respond to ``new'' proofs.}

Following Harold Bloom \cite{bloom2003map,BloomAnxiety}, I want to view Manin's ``refurbishing'' as acts of {\em misprision} or misreading.  In a different direction, Paraphrasing Grothendieck, to do new and creative mathematics frequently requires a certain naivet\'e \cite[p. 198]{GrothendieckRS}:\footnote{``Le petit enfant découvre le monde comme il respire — le flux et le reflux de sa respiration lui font accueillir le monde en son être délicat, et le font se projeter dans le monde qui l’accueille. L’adulte aussi découvre, en ces rares instants où il a oublié ses peurs et son savoir, quand il regarde les choses ou lui-même avec des yeux grands ouverts, avides de connaître, des yeux neufs — des yeux d’enfant.''}
\begin{quote}
	The child discovers the world like he breathes - the ebb and flow of his breathing make him inhale the world into his delicate being, and make him exhale into the world that welcomes him. The adult also discovers, in those rare moments when he forgets his fears and knowledge, when he looks at things or himself with wide-open eyes, eager to know, with new eyes--like those of a child.
\end{quote}
Putting these ideas together, doing mathematics, especially in the ``structural'' vein that has become the dominant paradigm, certainly around the types of mathematics I will discuss, consists of simultaneously focusing and forgetting, in the sense I mentioned earlier.  Reformulations made in the course of research emphasize some features, while deliberately obscuring or ignoring others.\endnote{I can take this idea even further, as it seems to apply to mathematicians renderings of the philosophies of their contemporaries and predecessors.  In Poincar\'e's review of Hilbert's Grundlagen he complains about the apparent diminishment of intuition, whose banishment is remarked upon repeatedly \cite{PoincareHilbert}.  Indeed, David Rowe, writing in \cite{RoweKleinHilbert} writes about Hilbert criticizing purely axiomatic mathematics.  He quotes Hilbert as saying:
\begin{quote}
Were this viewpoint correct, then mathematics should appear as nothing more than a series of logical arguments heaped one upon the other.  One would find nothing but an arbitrary series of conclusions driven by the power of logic alone.  But in reality nothing of the kind exists; indeed, the conceptual structure of mathematics is constantly led by intuition ({\em Anschauung}) and experience so that mathematics for the most part represents a closed structure free of arbitrariness.   
\end{quote}
This is evidently ``later'' Hilbert, from roughly 1919, but seems quite far from what Poincar\'e describes.  Nevertheless, it is a strong reading of earlier Hilbert that seems to support the view of Bourbaki's later projects in the popular mathematical consciousness of the day.}  We frequently do not remark upon that which has been forgotten, and it has to be unearthed.  I think it is important to pay attention to what has been forgotten: to understand that we have made many arbitrary choices to get to the mathematics of today.  Why did we make these choices?  Was it even conscious? 

\subsubsection*{Disclaimer}
Several words of warning are in order.  Evidently, I am not a professional sociologist (or philosopher, or historian, or...).\footnote{Barany and Kremakova write \cite[p. 2611]{BaranyKremakova}: ``Beyond their personal navigation of their professional worlds, some mathematicians pursue (with varying degrees of rigor and amateurism) social and humanistic studies of their subjects.''  I feel like this serves as an apt characterization of this work.  My disclaimer should also include the boilerplate phrase that it is ``susceptible to the biases and limitations'' of its author.}  As such my use of the literature in these domains is bound to be at best naive and at worst hopelessly wrong-headed.  Perhaps this work should be considered as a ``fan-fiction'' version of a sociological history of projective modules or maybe a first step towards what one might call ``mathematical criticism''.  Better: please view it as an invitation to a discussion.  

My referencing practices are more informed by my mathematical training than the examples I have seen in other literatures: I have attempted to give complete quotations when I use them, and I take responsibility for all errors of translation when that has been necessary.  As is perhaps clear from the preceding text, I have used endnotes to enclose long parenthetical discussions with additional references to literature that I think are relevant.  It is my belief that most of these comments are not essential to the argument contained in the main text, but I do use them to suggest questions I find interesting.  

This document has evolved rather organically, and in its current form is still somewhat tentative; for this I can only apologize to the prospective reader. I prefer to view it thus as closer to an alembic: I hope to eventually recast the collected distillate in a more definitive form.

\subsubsection*{Acknowledgments}
The seed for this work was probably planted when, in the middle of re-shelving books while working in a library in high school, I stumbled upon the ramblings of a disgruntled polymath Gynasium instructor turned cloistered scholar by the name of Oswald Spengler.  His massive two-volume tome with the provocative title ``The Decline of the West'' and imposing brutalist cover art \cite{SpenglerI,SpenglerII}, stuck out among the flashy detective novels.  Spengler's idea of cultures seemed bizarre to my burgeoning mathematical sensibilities, but I was too distracted by the diagnosis of cultural malaise to be too bothered.

When I first arrived at USC, through a web of connections, I was introduced to Martin Krieger who at the time was very focused on analogy in mathematics and physics; he had a particular fascination with the Langlands programme.  We lost touch during the pandemic, and I was saddened to learn that Martin passed away recently.  However, my various interactions with him over the years were extremely useful in circumscribing my ideas, and I would have really enjoyed sharing these thoughts with him. 

The actual argument in this text has been brewing for years, coming to life as a polemic that was tempered by long strings of conversations with two people to whom I evidently owe a huge intellectual debt.  On the one hand, Akshay Venkatesh has thought about interaction between AI and the practice of mathematics; discussions around these ideas played a formative role in my forays into ``inevitability'' in mathematics.  On the other hand, discussions with Mike Hopkins about history and influence in mathematics have greatly informed my conception of where ``interesting'' mathematics comes from.  I cannot imagine having written anything without these discussions.  In a related direction, the polemic to which I allude above had its source in the mathematical life lessons I have learned over the years with Jean Fasel; I thank him for his temperance and unflagging optimism.

In a less amorphous direction: I want to thank Jean-Pierre Serre and Hyman Bass for brief but extremely useful correspondence about their recollections of the history of projective modules.  While it probably goes without saying, their responses to my questions should not be viewed as an implicit endorsement of anything I am advocating (since I doubt they were aware of my project).  I thank Shane Kelly for criticism that helped me improve the structure of the argument.  I also thank Jesse Wolfson and Spencer Gerhardt for invitations to give lectures, preparations for which helped to clarify my thoughts.

\newpage
\section{A survey of pre-1951 mathematical culture}
\label{s:culture}
\addtoendnotes{\vskip 3.5em\noindent \begin{large}{\bf Mathematical culture pre-1951}\end{large} \vskip .2em}
The name ``homological algebra'' seems a pithy description of the collection of algebraic tools originating from the various definitions of homology.  The scope of this term has expanded considerably from its conception in Cartan--Eilenberg's {\em Homological algebra} with time; to see this, one need only open a more recent treatment like \cite{GelfandManin}.  Consequently, we attempt to reconstruct the mathematical cultures into which Cartan and Eilenberg's text was released.  

As we will recall in Section~\ref{ss:homologicalalgebra}, chronologically our story could begin with Poincar\'e's conception of homology, which gave generalizations of Betti numbers for ``suitably nice'' spaces. However, the story of the proliferation of homology, and its numerous variants will take us through various groups in the mathematical landscape.  By the time the work of Cartan--Eilenberg had appeared, homology had moved away from its original source, interacting fruitfully with, among other notions, group theory (initially the work of Eilenberg--Mac Lane), the theory of associative algebras (Hochschild), and the theory of Lie algebras (E. Cartan, G. de Rham and later Chevalley--Eilenberg).  The list of mathematical terms mentioned here serves only to indicate in modern terms the domains to which people interacting with homological ideas were exposed: algebra, topology, rounded out by complex and functional analysis.  While I indicate these distinctions here to serve as guideposts, disciplinary boundaries strike me as more porous then than now, so drawing precise lines is difficult.  

Instead, we choose to orient our discussion around people, rather than chronologically or by situating it around arbitrary mathematical domains.  Cartan and Eilenberg were both members of the Bourbaki collective, with Cartan a founding member.  Cartan ran a vibrant and very influential seminar in Paris, in the tradition of the S\'eminaire Hadamard and S\'eminaire Bourbaki; he was situated at the center of French academic hierarchy. Simultaneously, Eilenberg was a Jewish European emigrant, invariably enmeshed in the associated politics of refugee mathematicians in the United States.  He had a center of mathematical activity around himself at Columbia as well as numerous collaborators including Saunders Mac Lane (at Chicago) and Norman Steenrod (at Princeton at the time).

The ideas that led to up to the Cartan--Eilenberg book were developed over many years, by many people, interested in many kinds of problems.  Sections~\ref{ss:fiberbundles}-\ref{ss:sheaves} are devoted to sketching the  ``mathematical cultures'' around the initial release of the Cartan--Eilenberg book.  What did people know, and what were people interested in?  To answer this question, I will describe several scenes of mathematical activity before roughly 1951, focusing on ideas that, in my opinion, eventually informed that which was to give vital force to the theory of projective modules.  My aim is to describe, a cross-section of the mathematics under consideration that was viewed to be important by those around Cartan and Eilenberg.


\subsection{A bundle-theoretic unification theology}
\label{ss:fiberbundles}
\addtoendnotes{\vskip .2em\noindent {\bf A bundle-theoretic unification..} \vskip .2em}
By 1948, Henri Cartan was a well-established and extremely influential figure in the French mathematical community in particular and, more generally, in the international mathematical community.  He became a professor at the \'Ecole Normale Sup\'erieure (ENS) in 1940, following in his father Elie's footsteps.  While ENS initially was viewed as an institution to train teachers, it had slowly acquired the status of {\em the} place for training French academics and was invariably an elite institution by this point.\endnote{The ENS undoubtedly looms large in French academic society, and its culture provides background forces shaping the actions of our protagonists.  For an extra-mathematical example, Pierre Bourdieu's sociology, by his admission, was `` in reaction to the École Normale'' \cite{Riding}.  Speaking about pre-1968 French university society in general, Bourdieu \cite[p. 142]{BourdieuHA} writes:
	\begin{quote}
		It follows that appointment to the professorial body is subject to arbitrary decisions by the diverse authorities (and especially directors of research groups) whose choices are eventually validated and ratified by the body as a whole; and consequently that chances of appointment to research posts and increasingly, posts in higher education tend to depend at least as much on the scope, diversity and quality of academically profitable social relations (and thereby on place of residence and on social origins) as on academic capital.
	\end{quote}
Essentially all of the French mathematical figures in the period we discuss will be attached to the ENS.  Moreover, the over-representation of Normaliens in French mathematics has been analyzed by Bernard Zarca \cite{Zarca2006}, lending support to the idea that Bourdieu's analysis applies equally to French mathematics of the period we discuss.}

While there is much to be said about the cachet of the International Congress of Mathematicians (ICM) during the period 1930-1950, say by comparison to the status it currently enjoys within the mathematical community, Cartan's participation in this event serves as some measure of his international reputation: he had given a sectional lecture at the 1932 ICM held in Zurich, Switzerland and a plenary lecture in the 1950 ICM held in Cambridge, USA.\footnote{There is much to be said about the history of the ICM, I give only a few references here.  The 1932 ICM took place in the wake of the world-wide financial crisis, and we will have more to say about it momentarily.  Discussion of the 1950 ICM in the context of American Mathematics can be found in \cite[p. 496--509]{Parshall}.}

Cartan was a student of Paul Montel and, perhaps following Montel's interests, his early work was in complex analysis, specifically analysis of functions of one complex variable.  Complex function theory had deep historical roots in French mathematics stemming from the work of E. Borel and C. Picard.  The themes of Cartan's seminar seem a reasonable reflection of his research interests and the transcripts of the talks starting in 1948 provide a convenient starting point for our analysis. 


After the first world war, the mathematical enterprise in France was considerably diminished.\endnote{The scope of this diminishment has been explored in a number of places, especially in histories of Bourbaki, which is frequently described as a response.  Reinhard Siegmund-Schultze describes the situation thus \cite[p. 250]{SiegmundSchultzeIHP}:
\begin{quote}
	In 1920, when Émile Borel assumed the chair of mathematical physics and the theory of probability at the Sorbonne, it was generally understood that theoretical physics had fallen behind German developments by about a generation. The situation in pure and applied mathematics, as in several other disciplines, was not much better, even though France still had active and recognized researchers in some sub-disciplines, predominantly in analysis. The most basic problem was man-power. About half of the generation of young, promising French scientists had been killed on the battlefields of World War I....Historians agree on the relative decline of geometry in France; an almost total absence of algebra, number theory, and applied mathematics, and considerable weakness in mathematical physics..
\end{quote}}  Cartan explains that his interest was turned to higher-dimensional complex analysis \cite{CartanNotices}: 
\begin{quote}
	A little earlier, toward 1930, it is he who had oriented me toward the study of analytic functions of several variables by pointing out to me the work of Carathéodory on circled domains. One cannot overestimate what I owe to André Weil.
\end{quote}
Cartan and Weil were students together at the ENS in the 1920s.\endnote{For color, let us give an (incomplete) list of notable Normaliens from around the same time.  Outside of mathematics a list of people who finished around the time of Cartan and Weil includes: Georges Canguilhem (1924), Alfred Kastler (1921), Maurice Merleau-Ponty (1926), Louis Ne\'el (1924), Yves Rocard (1922), Jean-Paul Sartre (1924) and Andr\'e Weil's sister Simone Weil.  Other notable Normaliens in mathematics, many of whom will appear later in the story, include Claude Chevalley, Jean Delsarte, Jean Dieudonn\'e, Paul Dubreil, Marie-Louise Dubreil-Jacotin, Charles Ehresmann, Jacques Herbrand, Jean Leray, Charles Pisot, and Ren\'e de Possel.}  They were colleagues later in Strasbourg, and regular correspondents for many years; a number of the letters they exchanged may be found in \cite{CartanWeil}.  With the above comment about Weil's influence on Cartan as a jumping-off point, we want to describe Andr\'e Weil's point of view and its influence more generally. 

My goal in this section is describe and explore ideas of unification and analogy in mathematics and how they influenced the diffusion of ideas of the topological theory of fiber bundles into algebraic geometry.  I begin by exploring some aspects of Andr\'e Weil's biography, providing some background for his position in the mathematical community more generally, and highlighting some incidents that seem to bear on the nature of Bourbaki's mathematical practice.  I then recall some facts about the Bourbaki movement and its associated ideological convictions.  With this background in place, I review some aspects of the theory of fiber bundles in topology, then passing to Andr\'e Weil's efforts to transplant these topological ideas into algebraic geometry.  I argue that Weil's efforts were not valued in the broader algebro-geometric community.  Nevertheless Weil's point of view is important: I mention briefly in conclusion that it is this point of view that {\em vitalized} and broadened the conception of projective modules; this theme will be explored in greater detail in Section~\ref{ss:fiberbundlesandprojectivemodules}.

\begin{entry}[Andr\'e Weil]
	\label{par:Weil}
Undoubtedly, Weil was held in high esteem by Cartan, but Cartan suggests this was an opinion shared by many \cite{CartanNotices}: ``I knew that Weil was considered as a rather exceptional being, but I did not know that during his first year he had read Riemann and finished all his examinations.'' Weil spoke English, and German at home, had taught himself Greek, Latin and Sanskrit by the age of 16 \cite[Chapter 1]{WeilAuto}.\footnote{One can also find discussion in the breathless treatment of Weil in Aczel's biography of Bourbaki \cite{Aczel}.  This popular account of Bourbaki seemingly contributes to the mythology of Bourbaki.}  Undoubtedly he was viewed as a consummate scholar by those of his generation, but was well-respected by members of the previous generation as well. 

After finishing at the ENS in 1925, Weil writes
\begin{quote}
	almost all of my classmates left to do their military service as second lieutenants. Our class was the last one entitled to this rank without having to pass any examinations.  I was still too young to do likewise. 
\end{quote}
Unable to enter military service on account of his age, Weil applied for and received a modest scholarship, which while initially intended to be used to stay in Paris, he managed to modify to allow him to travel to Rome. Weil arrived in Rome in October 1925 to a gracious welcome extended by Vito Volterra.\endnote{To add some depth to the scene around this visit, we turn to \cite[Chapters 6-7]{GuerraggioPaoloni} about Volterra's life. Mussolini had come to power in 1922 and was prime minister by November of that year. Volterra worried greatly about the political situation in the wake of Mussolini's rise partly because of his (Volterra's) role in Italian and international science and mathematics, but also, Volterra being Jewish, because of rising tides of Fascist anti-semitism.  Volterra was elected to the Accademia dei Lincei one of the oldest and most prestigious European scientific institutions having Galileo as a member in 1923.  Around the same time he was president of the National Research Council and was the European referent for international arm of the Rockefeller foundation (the vaguely named International Education Board) \cite[p. 113]{GuerraggioPaoloni}.  
	
The political situation was tense by the time of Weil's visit.  Weil mentions the assasination of Giacomo Matteoti by Mussonlini's blackshirts, which took place in June 1924 and Benedetto Croce.  At this time, Giovanni Gentile who was the first Minister of Education in the first Mussolini fascist government wrote the {\em Manifesto Gentile}, which was an attempt to solicit intellectuals' adherence to the fascist movement.  Croce, an Italian senator, drafted the intellectuals response to this document.  

By the time of Weil's visit and Volterra had ``carried out his final duties as president of the Accademia dei Lincei and the National Research Council'' by June 1925 \cite[p. 125]{GuerraggioPaoloni}.  Against this backdrop, Weil arrives and reminisces of Volterra:
\begin{quote}
		Though probably a less universal mathematician than Hadamard, he was an admirable man in all respects. The king had named him Senator for life; he and Croce were 
		the two senators who, until the bitter end, come what may, voted against the Fascists. 
\end{quote}
Weil also recalls spending considerable time with Volterra's son Edoardo, who was roughly his age and about Volterra himself listening to his ``shapeless'' ideas about mathematics \cite[pp. 47-8]{WeilAuto}.  During this period, Weil reminisces also about meetings with O. Zariski, and F. Severi, as well as about Severi's views on Lefschetz: Severi compared Poincar\'e to an eagle and Lefschetz to a Hawk (Lefschetz's famous Analysis Situs paper had just been published in 1924).  Severi's interactions with politics at precisely this period led to fascist investigation ``for alleged administrative irregularities'' in his role as Rector of the University of Rome \cite[p. 65]{Guerragio} (he had been appointed by Gentile).  Famously, Severi later became a fascist in the early 1930s \cite{Guerragio}, though Volterra refused to take the fascist oath, a fact even remarked upon in BAMS \cite[p. 337]{VolterraBAMS}}  By his own admission, this trip had considerable effect on Weil's life.  After being pointed to work by Louis Mordell during this trip, which would eventually lead to his thesis (i.e., the theorem that is now known as the Mordell--Weil theorem), Volterra organized a Rockefeller grant for Weil to travel to Germany, where he visited G\"ottingen, Berlin and Frankfurt as well as to Sweden, where he visited Djursholm and Stockholm.  In this period, he met Richard Courant, Hans Lewy (Courant's assistant), David Hilbert, Emmy Noether, Bartel van der Waerden, Luitzen E.J. Brouwer, Paul Alexandroff, Gosta Mittag--Leffler and Max Dehn.  Cartan writes of this period that Weil ``became friends with the great mathematicians of the time, despite his being younger than they'' \cite[p. 633]{CartanNotices} undoubtedly served by his ability to perform the role of a scholar.  

Weil writes quite a bit about time in Frankfurt in 1926, especially discussions with Max Dehn and brief overlap with C.L. Siegel both of whom he held in extremely high regard (he refers to the latter as ``already a legend''.)  Reflecting attitudes about publication, Dehn already worried that ``mathematics was in danger of drowning in the endless streams of publications'', a worry still, nearly 100 years later.  Dehn professed a solution:
\begin{quote}
	...but this flood had its source in a small number of original ideas, each of which could be exploited only up to a certain point. If the originators of such ideas stopped publishing them, the streams would run dry; then a fresh start could be made. 
\end{quote}  
But mathematicians kept producing work
\begin{quote}
	at a time when the attitude reflected in the 
	American motto ``Publish or perish'' had invaded German universities and increasingly filled Siegel and others like him with disgust. 
\end{quote}
Weil also alludes to a text of Siegel's on transcendental numbers in Dehn's possession that the latter would only allow him to see in his home; he gives the sense of a world of mathematical secrets, literally kept under lock and key, revealed only to those in the ``know''.

As the 1930s opened, Weil took up a position in India that would last two years.  He returned to Europe in the summers, traveling widely, including back to Germany and England and then to a position in Marseilles in 1932.  He briefly mentions the International Congress in Zurich in September, glossing over the ``crisis'' then unfolding in Germany (the Nazis gained a plurality of seats in Germany and many communists come to power).

In November 1933, Weil moved to a position in Strasbourg, soon after Cartan had arrived there (Weil makes no mention of this fact in his autobiography, but we learn from Cartan that they had initially competed for a position that was won by Cartan in 1932).  In this period, in discussion with Cartan, the Bourbaki project was conceived, and that conception will play a background role in much of what I say below.\footnote{Harry Paul \cite[Chapter 9]{PaulKnowledgetoPower} gives some discussion of scientific research infrastructure in France, say around funding in the 1930s that is helpful in understanding what kind of national resources were devoted to research.} After a tumultuous period in the late 1930s and early 1940s, Weil arrived in the United States in March 1941 \cite[Chapter VII]{WeilAuto}.  To call this period ``formative'' for him is perhaps an understatement as it makes up the bulk of his autobiography.  After a period in San Paolo beginning in 1945, that is marked equally by social isolation and a flurry of correspondence, which I will discuss at greater length in Section~\ref{ss:sheaves}, Weil eventually settled at the University of Chicago in 1947, where he remained until he moved to the Institute for Advanced Study in 1958.  
\end{entry}


\begin{entry}[Interlude: Bourbaki and ``structuralist mathematics'']
	\label{par:bourbaki}
	Since I have mentioned Nicholas Bourbaki several times so far, before I go any further, I want to disentangle some things.  I have already mentioned ``modernism'' in the introduction, but for the forthcoming discussion the reader should be familiar with Hilbert's ``axiomatic method'' and have some ideas of the associated ``structuralism'' in mathematics. 

	As discussed in Section~\ref{s:prolegomena}, Jeremy Gray's book Plato's Ghost \cite{Gray}, presents an analysis of the transition from ``pre-structural'' mathematics to ``structural'' mathematics in the period 1870 to 1930.  This might be viewed as a first stage in the structural transition; the end of this period, which we mark roughly by the appearance van der Waerden's text ``Moderne algebra'' sees axiomatic approaches to mathematics appear as the dominant paradigm \cite[\S 1.3]{Corry} in mathematics;   This shift will be discussed in Paragraph~\ref{par:Noetherandgroups}.  
	
	The Bourbaki project provides a lens through which to view the next period, roughly 1930-1950.  As regards structuralism in mathematics: from a naive mathematical point of view, structuralist mathematics is the idea that mathematics is the study of structures.  This view, concordant with views espoused by Hilbert perhaps first crystallized in algebra, leaves a glaring hole: what is a structure?  One answer is that it is an attempt to axiomatize the numerous classes of objects that were being studied in algebra: groups, rings, fields, vector spaces, etc.  Bourbaki spent considerable time trying to come up with a precise definition of structure, encompassing all of these examples.  \todo{Tegmark siezes on this naive notion of culture.}  Structuralism in this sense is part of a style that Alma Steingart identfies in her beautiful book {\em Axiomatics} \cite[p. 18]{Steingart} as ``high modernism'' :
	\begin{quote}
		Starting in the 1930s and expanding in the postwar period, axiomatics came to be associated with universalizing theories and an emphasis on structure.
	\end{quote}
	But Steingart emphasizes that ``high modernism'' is broader than Bourbaki.

	Bourbaki was (and is), what I think is best described as a mathematical secret society.  Though retrospectively it was not so secret and much has been written about the group and its early members.\endnote{In \cite{Beaulieu}, Liliane Beaulieu aims to dispel some myths about the early history of Bourbaki, which have some bearing on what I say in what follows.  I will quote some of Dieudonn\'e's discussion on the aims of Bourbaki below in the famous ``Architecture'' \cite{Bourbakimanifesto}, but Beaulieu argues \cite[p. 242]{Beaulieu} ``this vision, however, did not prevail at the outset of the Bourbaki venture.  One may say that the often quoted ``Architecture'' article distorts the historical perspective to some extent by projecting Bourbaki's later hallmark onto the aims and concerns of the team at its inception.''  Nevertheless, I still believe that it is this {\em stronger} interpretation of the Bourbaki ideology that is reflected in the attitudes of the generation of mathematicians in the late 1940s and early 1950s that are crucial to the narrative.}  Traces of Dehn's attitude toward mathematical production and the gifting of mathematics seem visible in Bourbaki's practices and early distribution methods.  Initially, Bourbaki was conceived by Weil and Cartan (around the time of their joint stay in Strasbourg), purportedly to write a new textbook treatise on analysis.  The group consisted largely (entirely?) of Normaliens who had the lofty goal of lifting French mathematics back up to their retrospective view of its pre-WWI heights.\endnote{This goal hides numerous simmering tensions.  Before World War I, Poincar\'e was a huge figure in French mathematics, and French intellectual society more generally.  Poincar\'e, having attended and taught at \'Ecole Polytechnique,  self-identified repeatedly as a ``Polytechnicien''.  Other prominent French mathematicians of this period who were Polytechniciens included Charles Hermite and Paul Painlev\'e.  The \'Ecole normale was slowly becoming more prominent during the period just before World War I, but this change portends a tension between Polytechniciens and Normaliens.  Furthermore, the seeds of an ``applied'' vs. ``pure'' battle are seemingly planted in this era.  In the figure of Poincar\'e in particular, as compared to Bourbaki more generally, I can see a tension between ``intuitive'' and ``formalist'' mathematics as well.}  This initial project expanded to an effectively encyclopedic scope: Bourbaki aimed to write a definitive treatment of all of modern mathematics.  The Bourbaki attitude was unabashedly dogmatic, and there were many reactions to their approach.  The mathematicians involved in the Bourbaki project would grow to have tremendous influence within the mathematical establishment, especially around the fields we aim to discuss.
	
	Leo Corry describes \cite[p. 320]{CorryBourbakiStruct} the Bourbaki project as follows.
	\begin{quote}
		Like van der Waerden had done for thitherto disparate disciplines, which since then were included under the heading of ``modern algebra'', Bourbaki undertook the task of presenting the whole picture of mathematical knowledge in a systematic and unified fashion, within a standard system of notation, addressing similar questions, and using similar conceptual tools and methods in the different branches. But Bourbaki went a step beyond van der Waerden and attempted to provide a formal theory of structures affording a common foundation for all the other theories considered in his treatise. Bourbaki's work was originally motivated by the idea that the whole of mathematics may be presented in a comprehensive treatise from a unified, single best point of view the concept of structure was to play a pivotal role within it. This  conception, however, proved overconfident and Bourbaki soon realized that he must limit himself to include in his treatise only a portion of mathematics.
	\end{quote}
	
	By 1950, Bourbaki had different members and aims, summarized in their ``Architecture of mathematics'', an ideological manifesto of sorts \cite{Bourbakimanifesto}.  The Bourbaki ideology included various features that I will rely upon later, spelled out here.  First and foremost, Bourbaki ``shall not undertake to examine the relations of mathematics to reality or to the great categories of thought''.  Bourbaki thus isolates some core of mathematics from perceived external features reflecting, in part, the pure/applied dualism.  
	\begin{quote}
		Today, we believe however that the internal evolution of mathematical science has, in spite of appearance, brought about a closer unity among its different parts, so as to create something like a central nucleus that is more coherent than it has ever been.
	\end{quote}
	Bourbaki intends ``to remain within the field of mathematics'' and to proceed by analyzing ``the procedures of mathematics'' itself.  
	
	The ``procedures of mathematics'' to which they refer are the axiomatic method, which ``has its cornerstone in the conviction that, not only is mathematics not a randomly developing concatenation of syllogisms, but neither is it a collection of more or less ``astute'' tricks, arrived at by lucky combinations, in which purely technical cleverness wins the day''.  One sees in this a quest to isolate something like natural laws of mathematical evolution.  Inessential ingredients should be stripped away, so there is a ``minimalist'' sensibility to the approach.  But chiefly Bourbaki symbolizes unification and structure.   
	
	Of course, this summary itself of Bourbaki's ideology was written at near the end of the period I will consider in my discussion of ``mathematical cultures'' and evolved from perhaps less strongly-worded convictions present in group members themselves at the inception of the project.  Thus, I will try to be careful to distinguish the attitudes of members of Bourbaki from the collective: temperaments of individual members reflected the ``unification and structure'' paradigm to greater or lesser extents.  Another complication swept under the rug in the above description is posed by the Bourbaki seminars, which were organized to expose ``current'' mathematics of central interest (to those in the Bourbaki group).  This pursuit was arguably independent of the knowledge re-organization projects described above, but reflected similar expositional goals and thus its status would rise alongside Bourbaki.  
\end{entry}

In 1948, the University of Chicago, was teeming with new life: the newly arrived Weil was organizing/participating in a seminar that included several notable and also newly arrived mathematicians: the differential geometer Shing-Shen Chern had just returned from a stay in China to take up a position, and the topologist Edwin Spanier had just finished his degree with Norman Steenrod at the University of Michigan and taken up a position.  Hans Samelson, then a professor at Michigan also participated \cite[p. 259]{CartanWeil}.  Weil was already friends with Chern, having met him earlier in the 1940s when they overlapped during a stay in Princeton.  In a ``Topology caucus'' held in Madison, WI on Sep. 8, 1948, Weil, Samuel Eilenberg and Norman Steenrod drafted a plan for Nicolas Bourbaki to approach algebraic topology \cite[p. 260]{CartanWeil}: the first topic on the list was fiber bundles.  How did the Bourbaki vision interact with the forefront of mathematical research? 

\begin{entry}[Fiber bundles]
	\label{par:fiberbundles}
The preface of Norman Steenrod's 1951 text {\em The topology of fiber bundles} \cite{SteenrodFB} opens with the claim:
\begin{quote}
	The recognition of the domain of mathematics called fibre bundles took place in the period 1935-1940. 
\end{quote}	
Steenrod singles out Hassler Whitney as having given the first ``general'' definitions and the work of H. Hopf and E. Stiefel for applications to differential geometry.  Here, Steenrod is referring to \cite{Whitneyspherespace} where Whitney introduces the notion of sphere-space, eventually renamed ``sphere bundles''.  

After asserting that sphere spaces occur in nature, Whitney's stated goal is to define invariants to distinguish sphere spaces.  Particular examples highlighted in this paper arise from differentiable manifolds (the ``tangent'' sphere space) and embeddings of manifolds (the ``normal'' sphere space).  Whitney sketches the construction of homology classes playing the role of invariants. Whitney makes no claim of originality here, mentioning precursors in the work of Seifert-Threlfall \cite{SeifertThrelfall} and simultaneously mentions that examples of his sphere-spaces have been constructed by Hotelling \cite{Hotelling}.  Why was Whitney interested in sphere spaces?

To modern eyes, these sources might seem rather disparate: Hotelling was interested in dynamical problems proposed by Birkhoff, while the discussion of ``fiber spaces'' in Seifert--Threlfall is relegated to some remarks in an appendix, which themselves reference earlier work of Seifert \cite{Seifert}.  Tracing through the latter we find a common source.  In modern terminology, Seifert was interested in construction of three-dimensional manifolds, and was focused on the homeomorphism classification problem; already by this point the Poincar\'e conjecture was a focal question of the field then called combinatorial topology.  Of course, one source of Poincar\'e's interest in topology in general and three-manifolds in particular arose from his investigations of the dynamical three-body problem.  Seifert gives a long discussion \cite[\S 3]{Seifert} of an explicit continuous map $S^3 \to S^2$ with fibers that are circles, parenthetically remarking that the example he is describing is nothing other than Hopf's now eponymous map.  In modern terminology, Seifert is interested in locally trivial fibrations with closed curve fibers, especially circle fibers.  
 
Stiefel's thesis, written under the direction of Hopf simultaneously introduced the Stiefel ``manifold of frames'' and, extending work of Hopf on the index theorem, homology classes bearing obstructing existence of linearly independent vector fields.  In coordination, the work of Stiefel and Whitney formed a link between the problem of existence of linearly independent vector fields on manifolds and homology theory.  Thus, one sees fiber bundles can be used to systematically build manifolds, but simultaneously fiber bundles are used to build invariants; both conceptions are related to the homeomorphism classification problem.  Furthermore, these characteristic (co)homology classes are themselves becoming objects of independent interest.  

A more detailed account of the pre-1948 development of fiber bundles can be found in a note of McCleary \cite{McClearyFB}. By this point, Steenrod had given lectures at the University of Michigan and Princeton on fiber bundles in 1947-1948 (see the introduction to \cite{SteenrodFB}).  Progressively more general notions of fiber bundle were considered in this period because of applications (we will touch on some of these notions in the sequel), but Steenrod summarized the situation in 1951 thus:  
\begin{quote}
	The subject has attracted general interest, for it contains some of the finest applications of topology to other fields, and gives promise of many more. It also marks a return of algebraic topology to its origin; and, after many years of introspective development, a revitalization of the subject from its roots in the study of classical manifolds. 	
\end{quote}
\end{entry}

Eilenberg and Weil prepared a report on the theory of fiber bundles for Bourbaki \cite[p. 425]{BorelNotices} based on their interactions with Steenrod, but this document, as well as much of the material that Bourbaki discussed around algebraic topology never saw the light of day.  John McCleary has argued \cite{McClearyBBkiAT}, that Bourbaki's leaving discussion of the theory of fiber bundles unwritten was a result of the rapid pace of advance of algebraic topology at the time: presentation of this area was better suited to a more fleet-footed medium, such as a seminar. 

Weil was intimately familiar with theory of characteristic classes.  He wrote to Cartan in 1944:\footnote{From \cite[p. 94]{CartanWeil} \begin{quote}
		Ce qui interessera ton p\`ere, et vous tous, c'est que Chern (le m\^eme chinois qui a travaill\'e avec ton p\`ere in 1937) est \`a Princeton pour cette ann\'ee, et fait du tr\`es brillant travail sur les rapports entre topologie (particuli\`erement espaces fibr\'es) et g\'eom\'etrie diff\'erentielle.  Il a d\'ebut\'e par une d\'emonstration intrins\`eque, par les m\'ethodes de ton p\`ere, du th\'eor\`eme du Gauss--Bonnet pour poly\`edres riemanniens de dimension quelconque, que j'avais obtenu en collaboration avec un certain Allendoerfer en 1941, et publi\'ee aux Transactions.  Il est maintenant all\'e beaucoup plus loin, et on peut dire qu'il a une methode generale pour traiter tous les questions de ce genre.  
\end{quote}	}
\begin{quote}
	What will interest your father, and all of you, is that Chern (the same Chinese who worked with your father in 1937) is at Princeton for this year, and is doing very brilliant work on the relationships between topology (particularly fiber spaces) and differential geometry.  It began with an intrinsic demonstration, using your father's methods, of the Gauss--Bonnet theorem for Riemannian polyhedra of any dimension, which I had obtained in collaboration with a certain Allendoerfer in 1941, and published in the Transactions.  He has now gone much further, and it may be said that he has a general method for dealing with all questions of this kind.
\end{quote} 

Weil internalized and developed Chern's work going on to formulate what eventually became known as Chern--Weil theory \cite[1947b, p. 374]{WeilI} (submitted 8 Nov 1946 and appearing in 1947).\footnote{Chern's famous paper \cite{Chern} on characteristic classes of Hermitian manifolds appeared in 1946 and explicitly mentions complex analytic fiber bundles.  As Weil remarked in his letter, it in part generalizes some results of Allendoerfer--Weil on the Gauss--Bonnet theorem in the Riemannian setting, but the final paper does much more than Weil advertises.}  Chern--Weil theory provides a link between, on the one hand differential geometric ideas including the theory of connections and curvature, studied in great detail by E. Cartan, and G. de Rham's differential approach to cohomology and, on the other hand, the topological approach to bundle theory.  In brief, Chern--Weil theory describes characteristic classes of Hermitian vector bundles on complex manifolds in terms of E. Cartan's theory of curvature and connections.  Already by 1948, students around H. Cartan were learning this approach to characteristic classes and fiber bundles.  For example, Borel remarks about his attempts to understand the proofs of Weil's approach to Chern classes to Cartan \cite{BorelNotices}.\footnote{It's worth pointing out here that Weil also interacted with the anthropologist Levi-Strauss, and around this time he solved a combinatorial problem posed by Levi-Strauss \cite[p. ?]{WeilI}.  Levi-Strauss and Weil were classmates at the ENS and it's perhaps not surprising given the breadth of his interests that he had some awareness of anthropology.} 

\begin{entry}[Foundations of algebraic geometry revisited]
	\label{weil:foundations}
	Weil's ``Foundations of Algebraic geometry'' was published in 1946 \cite{WeilFoundations}, roughly the same time period as just discussed, and I would like to investigate its reception in the algebro-geometric community.  To understand Weil's viewpoint in this text, it will be useful to have a sense of how the notion of ``manifold'' was understood in this era; a slightly more detailed discussion of manifolds according to Poincar\'e can be found in Paragraph~\ref{par:Poincarehomology}.  Nowadays, I think it is standard to frame a dichotomy between ``embedded'' and ``abstract'' points of view on manifolds.  Briefly, the embedded point of view is that manifolds can be cut out of Euclidean space by equations (and inequalities).  In the ``abstract'' conception, one prescribes ``local models'' for manifolds, say open subsets of Euclidean space, that are glued together in a controlled fashion.  I think it is also now standard to view the ``embedded'' point of view as ``classical'', or at the very least older, and the ``abstract'' point of view as modern, the latter being codified by Veblen and Whitehead.  
	
	It probably comes as no surprise that this distinction is an oversimplification of matters.  To the extent that people thought about ``spaces'' in general, both conceptions seem present in the literature essentially from the start.  Erhard Scholz argues that Riemann was one of the first people to expound the concept of manifold.  Riemann's ideas about analytic continuation certainly influenced his conception of space, and Scholz \cite[p. 26]{Scholz} writes ``the distinction between local simplicity of manifolds, because of the presupposition of local coordinate systems, and globally involved behaviour was indicated by Riemann'' but only emphasized in some of his works.\footnote{Riemann in his habilitation also placed Gauss alongside the psychologist/philosopher J. Hebrart in providing the groundwork for his conception of manifold \cite[p. 288]{Laugwitz}.}  It is nevertheless embedded in modern ideas about analytic continuation of analytic functions.  A related discussion of Riemann's ideas about space and, in particular, the classification of surfaces, can be found in \cite[2.1]{Volkert}.  Poincar\'e also rather explicitly uses both the ``embedded'' and ``abstract'' conceptions of manifold: he introduces the former in \cite[\S 1]{PoincareI}, but then describes the latter in \cite[\S 3]{PoincareI}, before attempting to reconcile the two points of view.  In fact, in his constructions of manifolds, say via identification spaces, Poincar\'e certainly used both points of view.  
	
	In a sense, Weil's Foundations of algebraic geometry puts this dichotomy between ``abstract'' and ``embedded'' for manifolds at the forefront.   Algebraic varieties in the 19th century were also primarily thought of in ``embedded terms'', either as subsets of affine or projective space. Bowdlerizing some aspects of Hodge's description of algebraic geometry \cite{HodgeLefschetz} in that period, there were several schools broken down along roughly national lines: a German school of Abel, Riemann and Weierstrass.\footnote{My discussion of Riemann above indicates why I say ``primarily'': certainly Riemann surfaces gave particularly interesting examples of algebraic varieties, and the ``abstract'' point of view was useful here.}  A more geometric point of view was envisioned by Noether, and eventually taken up by Castelnuovo and Enriques.  A French school, directed by Picard and Simart introduced a third function-theoretic school.  The algebraic point of view was further developed by Dedekind and Krull.  
	
	None of these authors tried to introduce any sort of ``abstract'' point of view on algebraic varieties and Weil was one of the first to do so; he was rather explicitly inspired by the notion of abstract manifold: one would like to describe algebraic varieties by ``gluing or patching together'' suitable local models.  In Weil's words: ``this idea, inspired by the usual definition of a  topological manifold by means of overlapping neighborhoods, leads to the definition of an ``abstract variety'' \cite[p. xi]{WeilFoundations}. Transplanting the ``abstract'' notion of manifold to algebraic varieties was complicated by several factors.  
	
	First, while algebraic varieties over the complex numbers were frequently studied by classical topological techniques, what happens if one works over other fields? Owing to his interest in finite fields and the Riemann hypothesis for curves over finite fields, Weil certainly had such applications in mind (see Paragraph~\ref{entry:weilanalogy} for further discussion).  
	
	Second, what should one take as local models, i.e., what was a suitable analog in algebraic geometry of an open subset of Euclidean space?  In what to me seems a telling remark regarding the primacy of topological thinking in his approach Weil writes \cite{WeilFoundations}:
	\begin{quote}
		Similarly, however grateful we algebraic geometers should be to the modern algebraic school for lending us temporary accommodation, makeshift constructions full of rings, ideals and valuations, in which some of us feel in constant danger of getting lost, our wish and aim must be to return at the earliest possible moment to the palaces which are ours by birthright, to consolidate shaky foundations, to provide roofs where they are missing, to finish, in harmony with the portions already existing, what has been left undone. 
	\end{quote}
	
	Some sense of the value of Weil's treatment and the value of this generalization of the notion of variety can be gleaned from Zariski's review of Weil's book for the Bulletin of the AMS \cite{ZariskiWeil}.  About Weil's eschewing of algebra, Zariski writes:
	\begin{quote}
		It is a remarkable feature of the book that—with one exception (Chap. III)—no use is made of the higher methods of modern
		algebra. The author has made up his mind not to assume or use modern algebra ``beyond the simplest facts about abstract fields and
		their extensions and the bare rudiments of the theory of ideals.''  Weil's faithful realization of this program of strict mathematical economy is an achievement in itself. In some cases this leads to the ``best possible'' proofs. However, on the whole one may question the wisdom of this self-imposed régime of austerity.
	\end{quote}
	In ``questioning the wisdom'', Zariski seems to be alluding to whether the resulting text is usable or rather an exercise in style.  His opinion on this becomes clearer when one reads his commentary on the presentation, which is in the axiomatic style:
	\begin{quote}
		As a result, the reader finds himself very much in the position of a man who must collect a large amount of cash most of which is in pennies...Most readers will find it difficult to follow the author through the seemingly endless series of propositions, theorems, lemmas and corollaries (their total must be close to 300). A further obstacle to continuous reading of the book are the numerous crossreferences in the proofs...
	\end{quote}
	Zariski's complaint about this style of writing, which is a hallmark of Bourbaki texts and was imitated later by others, applies equally well to more modern treatments of algebraic geometry.  I claim it increases the socially dependent nature of the text.  Indeed, by simply reading such a text one does not get a sense of which aspects of the arguments are ``routine'', i.e., part of the culture (in the sense I used it in Section~\ref{ss:whathappenshere}), and which are ``novel'', which is a complaint I have heard repeated by many.  One is thus led to seek out competent guides who can signpost the text or, lacking that, to develop one's own signposts, however idiosyncratic. 
	
	Finally, Zariski comments about the value of the ``abstract'' notion of variety:
	\begin{quote}
		The long Chap. VII is dedicated entirely to what Weil calls ``abstract varieties,'' or Varieties with a capital V...We find it very difficult to estimate the necessity or the permanence value of this new concept. 
	\end{quote}
	He highlights that Weil's notion of ``completeness'' has no examples beyond classical ``projective'' examples.  In his 1950 ICM address \cite{ZariskiICM}, Zariski mentions Weil's abstract varieties only briefly, writing:
	\begin{quote}
		An even more radical revision of the concept of a variety has been offered by André Weil. His so-called abstract varieties are not defined as subsets of the projective space, but are built out of pieces of ordinary varieties, pieces that must hang together in some well-defined fashion. It is still an open question whether the varieties of Weil can be embedded in the projective space. 
	\end{quote}
	Thus, while the distinction between ``embedded'' and ``abstract'' seems present in the theory of manifolds essentially from the outset, it apparently strikes Zariski, even in 1950, as artifice, and perhaps irrelevant to his conception of ``mainstream'' algebraic geometry.  Weil was no doubt aware of this possibility writing \cite[p. viii]{WeilFoundations}: ``Of course every mathematician has a right to his own language — at the risk of not being understood; and the use sometimes made of this right by our contemporaries almost suggests that the same fate is being prepared for mathematics as once befell, at Babel, another of man's great achievements.'' 
\end{entry}

\begin{entry}[Weil to Cartan on fiber bundles]
	\label{paragraph:WeilFB}
Weil apparently did not take Zariski's criticism to heart.  In January 1949,	Weil introduced, in a talk at the University of Chicago, a version of the theory of fiber bundles in algebraic geometry \cite[p. 411-12]{WeilI}.  By several accounts (see the autocritique later for further discussion of this), his was the first attempt to make such a definition.  Weil admits he was unsure about the use of the Zariski topology in his Foundations of Algebraic geometry (a fact he elliptically remarks upon in the introduction) \cite{WeilFoundations}, but remarks \cite[p. 566]{WeilI}:\footnote{``C'est bien, je crois, la th\'eorie des espaces fibr\'es qui me convertit d\'efinitivement \`a l'usage de cette topologie, malgr\'e quelques h\'esitations initiales.}
\begin{quote}
	It is, I believe, the theory of fiber spaces which definitively converted me to the use of this topology, despite some initial hesitations.
\end{quote}
Note that Zariski's eponymous topology was introduced around 1944 \cite[p. 683]{Zariskitopology}, and appears in \cite[\S 2]{ZariskiICM}.

Nevertheless, Weil gives a series of examples intended to demonstrate that theory of fiber bundles in algebraic geometry was rich and allowed the unification of a number of disparate strands of problems.  He studies fiber spaces in algebraic geometry where the structure group is the multiplicative group and observes that isomorphism classes of such bundles coincide with divisor classes modulo linear equivalence, a notion that had long been studied in arithmetic situations, say going back to Dedekind--Weber.  He mentions the case where the structure group is the group $PGL_2$ and observes that in this case the theory of fiber bundles is connected with the theory of ruled surfaces, which he marks ``as a problem of considerable difficulty and interest.''  He describes a theory of the tangent bundle and remarks on corresponding definitions of characteristic classes of algebraic varieties and invariants of Stiefel--Whitney--Chern type.  He then goes on to speculate about definitions of fiber bundles in arithmetic contexts closing with the bold assertion:
\begin{quote}
	We thus stand at the threshold of a very promising new field of research, which it will certainly take many years to develop fully. 	
\end{quote}
\end{entry}

\begin{entry}[Interlude: Weil and analogy]
	\label{entry:weilanalogy}
At this point it is natural to wonder if there is anyone but Weil who could have initiated a discussion like this: the unique collection of historical referents he uses in his justification reads as a summary of his research interests and experiences.  Weil's oeuvre, as well as his informal writings demonstrates an intense fascination with certain kinds of generalization and analogy.  

For example, his 1938 paper \cite[p. 185]{WeilI} proposed a ``matrix'' generalization of theta functions.  Atiyah references this in \cite{Atiyah} as an early approach to the theory of fiber bundles in algebraic geometry, which I will discuss more momentarily.  Looking back much later, Weil himself, in his commentaries on this paper \cite[p. 537]{WeilI} views it as this and more:\footnote{``Avec le recul dont on dispose \`a pr\'esent, on peut dire que l'essentiel dans [1938a] est d'avoir inaugur\'e l'\'etude des fibr\'es vectoriels sur une courbe alge\'brique, \`a fibres de dimension quelconque $r > 1$; subconsciemment je m'efforçais de construire pour ceux-ci des ``varietes de modules'' que, faute de notions claires sur ce sujet, j'aurais \'et\'e bien en peine de d\'efinir.''}
\begin{quote}
With hindsight, we can say that the essential thing in [1938a] is to have inaugurated the study of vector bundles on an algebraic curve, with fibers of any dimension $r > 1$; subconsciously I tried to construct for these ``moduli varieties'' that, for lack of clear notions on this subject, I would have been hard pressed to define.
\end{quote}
The theory of moduli of vector bundles on curves was to become much more active in the 1960s.  

Weil had been, one might say, obsessed about the analogy between number fields, function fields, and complex function theory, a.k.a., the theory of Riemann surfaces for some time.  While this analogy itself goes back at least to Dedekind,\footnote{See for example \cite[Chapter 5 IV p. 242]{Krieger} for another discussion of these points and, in particular, commentary on F. Klein's summary of the ``dictionary'' in 1926.  See \cite{OortSchappacher} for a discussion of points of history related to the Riemann hypothesis in positive characteristic.}  Weil explicitly discussed meta-mathematics of analogy in several places, but perhaps most famously in the ``Rosetta Stone'' that he laid out in correspondence with his sister \cite[p. 236]{WeilI}.  For the sake of the chronology, let me mention that Weil gave a proof of the Riemann hypothesis for function fields of curves near this time as well \cite{WeilRH}, eventually leading, around the time of his note on fiber bundles, to the Weil conjectures (cf \cite{WeilConjectures}).
\end{entry}

I think it's impossible not to view the transplant of the notion of fiber bundles from topology to algebraic geometry from this context.  All of these results justified Weil's vaunted ``intuition''.  However, the definition of fiber bundle in algebraic geometry was still at this stage extremely ad hoc.   Weil's definition of abstract algebraic variety was conceived to imitate the definition of abstract manifold. In contrast to the definition of an abstract manifold, where one only needed to mention existence of covering charts, in Weil's algebro-geometric setting the covering charts were part of the data, complicating the presentation.  The resulting complexity spilled over into the theory of fiber bundles over algebraic varieties where once again covering charts were part of the data.\footnote{Perhaps one inadvertent side-effect of the explicit nature of charts in Weil's definition of algebraic fiber bundle, which will reappear when we discuss cohomology and sheaves later, is that one may readily provide an interpretation in terms of Čech cohomology. Cognizant of the fact that it may appear at this stage like I am projecting contemporary biases onto old mathematics, our later discussion of the Cartan--Weil correspondence around this period, especially as regards sheaf-theory will suggest that this interpretation seems plausible.}  

From the standpoint of aesthetics of the axiomatic approach in the hands of Bourbaki (especially in light of the ``austerity'' that we observed in Zariski's characterization of Weil, though one might also call it minimalism), these additional choices can be viewed as ``inelegant''.  As regards this accusation, Weil characterizes as ``subtle'' the difference between his point of view and the differential geometric point of view:\footnote{Chez eux, il s'agit de d\'efinir une structure sur un ensemble pre\'existant, en le recouvrant par des voisinages; chez moi cet ensemble n'appara\^it qu'{\em a posteriori} par recollement de morceaux donn\'es individuellement \`a l'avance. }  
\begin{quote}
	For them, it is a question of defining a structure on a pre-existing set by covering it with neighborhoods; for me this set only appears {\em a posteriori} by gluing together pieces given individually in advance.
\end{quote}
To me, another ``subtle'' difference can be described as follows: in the differential geometric conception of manifolds, local models are open subsets of Euclidean space, and therefore ``all of the same type'', whereas for Weil, local models were open subsets of affine varieties, so there are {\em different types} of local models, giving Weil's algebraic varieties a chimeric nature.

The resulting theory was, nevertheless, still rich enough to address Hilbert's 15th problem of putting enumerative algebraic geometry of Schubert on firm footing, so good enough for ``practical purposes'' in some sense\footnote{Weil highlights this in the introduction to \cite{WeilFoundations}, writing ``Our results include all that is required for a rigorous treatment of so-called ``enumerative geometry'', thus providing a complete solution of Hilbert's fifteenth problem.} but it certainly does not appear to rise to the Bourbaki standard then espoused.  Indeed, Borel remarks \cite[p. 423]{BorelNotices} the state of the theory was such ``for lack of suitable concepts'' and the resulting theory is ``rather unwieldy and requires a somewhat discouraging amount of algebra to be worked with. Weil himself writes \cite[p. 556]{WeilI} \footnote{``Pour superficielle qu'elle soit, cette diff\'erence n'en a pas moins contribu\'e \`a compliquer inutilement les d\'efinitions relatives aux vari\'et\'es abstraites au Chap. VII de mes Foundations, et peut-\^etre \`a rebuter bien des lecteurs. Mais en ce temps je n'avais guere le choix.''}
\begin{quote}
	Superficial as it may be, this difference has nevertheless contributed to unnecessarily complicating the definitions relating to abstract varieties in Chap. VII of my Foundations, and perhaps to putting off many readers. But at that time I had little choice.
\end{quote}
Presented by another mathematician, not enmeshed in a group like Bourbaki, one wonders whether the resulting theory and its evident complexities would have been ignored.  Instead, the language of Weil's Foundations was, eventually, widely used, and the inherent complexities gave birth to a simplifying impulse; we will come back to this when we discuss \cite{SerreFAC} which serves as a partial response.

Undoubtedly the theory of fiber bundles was well-represented on both sides of the Atlantic.  By 1949-50, the Cartan seminar had focused on fiber bundles as well, but their interests appear rather different.  The seminar was thinking about bundle theory and characteristic classes in topology, and the Cartan--Serre method of killing homotopy was being developed.  In contrast, Weil was taken with fiber bundles as a unifying theme, and a letter from Weil to Cartan \cite[p. 311]{CartanWeil} seemingly pushes for a shift away from thinking about fiber bundles only in topological/differential-geometric contexts.  

Weil opens the letter with an entreaty for Cartan to take its contents into account in his preparations for the 1950 ICM.  Using a complex analytic theory of fiber bundles, Weil formulates the Cousin problems in terms of fiber bundles, and observes that, following work of Oka, for ``domains of holomorphy'' (later generalized by the notion of Stein manifold) the only obstructions to holomorphic solutions are topological.  In other words, he reformulates what would later be called the Oka principle in terms of principal bundles under the group of non-zero complex numbers $\cplx^{\times}$.  He closes the letter by writing:\footnote{``
	Une fois qu'on a pris l'habitude de voir les espaces fibr{\'e}s dans ces questions, on ne tarde pas a se convaincre qu'ils apparaisant partout (ou tout au moins \guillemetleft presque partout \guillemetright) et qu'on gagn'e \'enorm\'ement... Il me para\^it certain que presque tous les probl\`emes o\`u il s'agit de recoller des donn\'ees locales pour passes \`a du global sont de la nature \guillemetleft espace fibr\'e \guillemetright.''}
\begin{quote}
	Once we get used to seeing fiber spaces in these questions, we quickly become convinced that they appear everywhere (or at least ``almost everywhere'') and that we gain enormously... It seems certain to me that almost all the problems where it is a question of gluing together local data to produce global data are of ``fiber space'' nature.
\end{quote}
Borel cites part of this passage in \cite[p. 426]{BorelNotices} observing that only a fragment of the letter was available to him, noting that this fragment alone ``makes one strongly wish to see the rest''.  Cartan responds with interest, requesting clarifications and along the way notes:\footnote{``	Mais je n'avais pas vu clairement le lien avec les espaces fibres.''}
\begin{quote}
  But I didn't clearly see the link with fiber spaces.
\end{quote}	
In his collected works, Cartan discusses these ideas much less, writing only \cite[p. XIV]{CartanI}:\footnote{Les premières indications relatives à l'utilisation de la théorie des faisceaux pour l'étude des fibrés holomorphes remontent à une conférence que j'ai faite au Séminaire BOURBAKI (Décembre 1950).} 
\begin{quote}
	The first indications relating to the use of sheaf theory for the study of holomorphic bundles date back to a presentation that I gave in the BOURBAKI Seminar (December 1950).
\end{quote}
The Bourbaki talk to which Cartan refers \cite{CartanBourbaki} will be discussed much more later, but it presents a nicely distilled version of the discussions of Cartan and Weil.

\subsubsection*{Auto-critique}
When I began to write this narrative of fiber bundles, I couldn't help but fall more and more in the thrall of Weil.  How magnificent his intuition!  His vision was amazing!  But Zariski's criticisms poked a hole in this idea suggesting that Weil's conception might not have been widely accepted.  Zariski's 1950 ICM address makes clear that, while algebraic geometers were aware of Weil's Foundations, many could proceed as if they effectively did not exist.  

Nevertheless, I want to take a step back and try to investigate my point of view.  Three separate questions arise.  First, ``But why Andr\'e Weil''?  Second, am I just subscribing to another ``great-man'' narrative of mathematics?\endnote{Weil himself subscribed to some kind of ``great-man'' view of history, writing \cite{WeilAuto}:
	\begin{quote}
		I had become convinced that what really counts in the history of humanity are the truly great minds, and that the only way to get to know these minds was through direct contact with their works. I have since learned to modify this judgment quite a bit, though I have never really let go of it completely. 
	\end{quote}
	In H. Spencer's description \cite[p. 31]{Spencer}:  
	\begin{quote}
		The lessons given to every civilized child tacitly imply, like the traditions of the uncivilized and semi-civilized, that throughout the past of the human race, the doings of conspicuous persons have been the only things worthy of remembrance. 
	\end{quote}
	However, Spencer counters with \cite[p. 34-35]{Spencer}
	\begin{quote}
		then you must admit that the genesis of the great man depends on the long series of complex influences which has produced the race in which he appears, and
		the social state into which that race has slowly grown.
	\end{quote}
	One wonders if this admission of social factors is what Weil means when he speaks of ``modified judgment.''
}  Third: because historians of mathematics have been very critical of Bourbaki's historiography, and my only sources are largely Bourbaki members themselves, is the treatment above just repeating a Bourbaki-approved version of the events?

\begin{entry}[But why Andr\'e Weil?]
	Maybe a more precise question here is: why did Andr\'e Weil feel compelled to introduce a definition fiber bundles in algebraic geometry, and would someone else have done the same thing had he failed to do so?  In asking a question like this about ``inevitability'' I follow Mac Lane who asks precisely this question about category theory (see the discussion of Paragraph~\ref{par:categoricalinevitability}).  
	
	To approach this question, I want to investigate at least one other mathematician interested in applying ideas of topology to algebraic geometry: Solomon Lefschetz (though there are others; we will very briefly touch on Hodge and what one might call the English school).  Solomon Lefschetz famously wrote \cite{LefschetzAuto}:
	\begin{quote}
		As I see it at last it was my lot to plant the harpoon of algebraic topology into the body of the whale of algebraic geometry. But I must not push the metaphor too far. 
	\end{quote}
	Indeed, many of the tools Lefschetz developed in algebraic topology (to which we will return shortly) were created precisely to answer questions about topology of complex algebraic varieties.  
	
	Lefschetz was also educated in France before emigrating to the United States.  Like Weil, he studied with Picard and was very interested in complex function theory.  Eventually, Lefschetz developed his fixed point formula about which Steenrod \cite[p. 24]{LefschetzSymposium} writes:
	\begin{quote}
		The fixed-point problem seems to have dominated nearly all of Lefschetz's work in topology. Some dozen papers appearing during
		the period 1923-38 were concerned directly with the problem. 
	\end{quote}
	But there is a related direction that later became important in algebraic geometry: Lefschetz studied what is now known as the Picard--Lefschetz formula \cite{PicardLefschetz}. Moreover, Lefschetz was interested in even the foundations of algebraic geometry through the writing of his book \cite{LefschetzAG}, which appeared in 1953.  Could Lefschetz have invented the theory of algebraic fiber bundles?
	
	In modern terminology, Picard--Lefschetz theory studies the following situation: one takes a suitable proper holomorphic function on a complex manifold and aims to study the action of the fundamental group of the complement of the critical locus on the homology of the fibers; the situation of interest is now called a Lefschetz fibration.  Once the concept of fiber bundle was invented by Whitney, it is hard to imagine that Lefschetz would not have been immediately aware that Lefschetz fibrations were examples of fiber bundles.  
	
	There are two distinctions to be made here.  First, as I explained above, Weil seems to be the person who really pushed studying algebraic varieties using the Zariski topology.  In this direction: ``interesting'' Lefschetz fibrations are essentially never locally trivial in the Zariski topology and would have required an even more radical conception of ``topology'' in algebraic geometry (a suitable such notion was invented later by Grothendieck, but was not seemingly inspired by this analysis).  Second, there is an issue of mathematical style.  Given what is written about Lefschetz's mathematics, he was less concerned about ``unifying structures'' and themes than he was about his own idea of ``originality''.  As such, it seems unlikely that he would have focused attention on these kinds of structures in particular.  In other words, it seems unlikely that Lefschetz would have isolated ``algebraic fiber bundles''.  Instead, the Bourbaki conception of ``structure'' and Weil's romantic associated with analogy seem to conspire in Weil's definition.
\end{entry}

\begin{entry}[And\'re Weil and great-man narratives]
As regards the second point, there have been many histories of mathematics that paint mathematics as a sequence of events performed by ``great men.''  Let's try to be a little more precise about this: we need to admit that Weil's status within the mathematical community played an important role in acceptance of his ideas.  Here, status is used in the Weberian sense discussed in Economy and Society.  In fact, use of these concepts with regards to mathematical cultures is not original, for example, Michael Harris discusses charisma in this context in \cite{HarrisMWO}. Undoubtedly Weil possessed charisma in this sense.\endnote{By 1948-49, Weil was arguably one of the most important mathematicians on the scene.  One testament to this appears in the discussions around the award of the 1950 Fields Medals as discussed in \cite{BaranyFM195058}.  Among other things, Barany argues that the age-limit now more or less formally imposed on the awards of the fields medal turned on whether the 1950 award would or would not be awarded to Weil.  He remarks: `` Bohr cryptically suggested that a cut-off of 42 “would be a rather natural limit of age”.''  He citing some strong opinions of other committee members that came to light during the discussion: Alfhors worried that to give a medal to Weil would be “maybe even disastrous” because “it would make the impression that the Committee has tried to designate the greatest mathematical genius.”  Damodar Kosambi thought it would be “ridiculous” to deny him a medal. Hodge worried “whether we might be shirking our duty” if Weil did not win.}   To say it differently, if Weil had not actively pushed his point of view, it might have led to completely different research priorities.  I think it's safe to say that we would not have had the same perception of fiber bundles in algebraic geometry without Weil's influence.

A related question is {\em to whom} was this influence directed?  Cartan's comments suggest that Weil was esteemed among colleagues and older mathematicians.  Furthermore, his sojourn in India undoubtedly left a mark on Indian mathematics, and a conduit by which Indian mathematicians could begin to have contact with the undoubtedly more well-established European mathematical scene. However, Weil, in contrast to Cartan, seems much less successful in an advisorial capacity.  If one believes the information at the mathematics genealogy project (which is notoriously unreliable for older mathematicians), then Weil's only students around this time were Harley Flanders and Arnold Shapiro. Shapiro himself wrote a thesis on the cohomology of fiber bundles \cite{ShapiroThesis}, but focused more on topological directions, perhaps in part due to his close friendship with Raoul Bott.  Bott refers to Shapiro \cite[p. 532]{Duren2} as his ``personal remedial tutor''.\footnote{Bott writes: ``But these people — together with Kodaira and Spencer — and my more or less ``personal remedial tutor,'' Arnold Shapiro, were the ones I had most mathematical contact with.''  Bott also refers to his ``good friend'' Shapiro in \cite[p. 378]{JacksonBott}	where he details Shapiro's contributions to the proof of Bott periodicity.  Later Shapiro gave a Bourbaki seminar about his joint work with Bott on Clifford algebras and the periodicity theorem \cite{ShapiroBourbaki}.  Shapiro died in 1962, and the results exposed in the Bourbaki seminar became what is now known as the Atiyah--Bott--Shapiro theorem \cite{AtiyahBottShapiro}} What is certainly clear is that none in his academic lineage contributed to the theory of algebraic fiber bundles.  In other words, we have to assume that it was only through his colleagues that theory of algebraic fiber bundles gained credence.
\end{entry}


\begin{entry}[Bourbaki self-mythologizing]
Bourbaki as a group, and Weil and Dieudonné in particular, have been harshly criticized for their views on history of mathematics; see for example \cite[pp. 329-338]{Corry} and the references therein.  In particular, Corry writes \cite[p. 331]{Corry}:
\begin{quote}
	Bourbaki's historiography, as manifest in the Elements d'histoire as well as in the individual writings of Dieudonné and of Weil, has been strongly connected with their overall conception of mathematics. In particular, they have applied similar criteria to differentiate important from unimportant ideas in both present mathematical research and past mathematical theories.
\end{quote}
Related discussion appears in \cite[Chapter 6]{Steingart}.  The MathSciNet review of \cite{BourbakiHistoire} includes the description
\begin{quote}
In their presentation the authors have been guided by their view of modern mathematics. It is this retrospective view of the development of mathematics (outweighing the usual genetic description) which is probably the most important and unusual feature of this torso. It also accounts for the heavy emphasis given to the mathematics of the 19th and 20th century.
\end{quote}
Note, however, that these descriptions seems to match many other (but not all!) historical criticisms of histories of mathematics written by mathematicians themselves.\endnote{Cartan, writing in \cite{CartanonBourbakiEng}, has this to say about Bourbaki's historical interludes:
\begin{quote}
	The ``Notes Historiques'' and ``Fascicules de R\'esultats'' deserve special mention. Bourbaki often places an historical report at the end of a chapter. Some of them are quite brief, while others are detailed commentaries. Each pertains to the whole matter treated in the chapter. There are never any historical references in the text itself, for Bourbaki never allowed the slightest deviation from the logical organization of the work. It is only in the historical report that Bourbaki explains the connection between
	his text and traditional mathematics and such explanations often reach far back into the past. It is interesting to note that the style of the ``Notes Historiques'' is vastly different from that of the rigorous canon of the rest of Bourbaki's text. I can imagine that the historians of the future will be hard put to explain the reasons for these stylisticdeviations.
	\end{quote}
This note is translated from German, but what is one to make of Cartan's implicit distinction between Bourbaki's presentation and ``traditional mathematics''?  At its inception, Bourbaki was certainly an avant-garde institution, but by the time this was written, many aspects of Bourbaki's terminology and choice of presentation had become widely accepted.}

Later, Dieudonn\'e wrote a ``Panorama of mathematics'' as seen by N. Bourbaki \cite{DieudonnePanorama}.  One aim of this book is to adumbrate the ``major parts'' of mathematics.  Dieudonn\'e begins by pronouncing that ``the history of mathematics shows that a theory almost always originates in efforts to solve a specific problem'' before going on to classify the kinds of problems that can appear: stillborn problems, problems without issue, problems that beget a method, problems that belong to an active and fertile general theory, theories in decline, and theories in a state of dilution.  His tone is unapologetically judgmental about the first pair and last pair in his classification, and situates Bourbaki's interests among the third and fourth classes asserting that ``this is, I believe, as objective an opinion as I can form, and I shall abstain from further comment.''  Each chapter of the book contains at its end ``a list of the mathematicians who have made significant contributions to the theories described'' \cite[p. 3]{DieudonnePanorama}.  It is perhaps unsurprising that Bourbaki-aligned mathematicians feature prominently in these lists, with Weil appearing most regularly.  Landsburg in his reminiscences on Weil \cite{Landsburg} writes ``with 11 major contributions to his credit, Weil's name appeared more often than any other.''\footnote{I do not begrudge Landsburg his evident reverence for Weil, about whom he writes: ``Weil's presence was enhanced, as is the case with many great geniuses, by his personal eccentricities and the legends they inspired---the strangely guttural French accent, the acerbic wit, the exacting standards, the complete inability to tolerate any form of stupidity (quite a burden for a man compared to whom almost everyone else in the world was basically a dunce), and the mischievous vanity. These traits live on in his writings and in the oral history that is lovingly preserved by mathematicians worldwide.''}

Grattan-Guinness \cite{GrattanGuinness} writes more broadly of histories produced by mathematicians, characterizing them as ``an account of how a particular modern theory arose out of older theories instead of an account of those older theories in their own right''. Admittedly, my focus is {\em precisely} on ``how we got here'' and where a particular mathematical idea arose from past mathematical theories.\endnote{Weil had his own strong opinions about how the history of science should be written.  He wrote \cite[p. 231-32]{WeilICM1978}:
	\begin{quote}
		How much mathematical knowledge should one possess in order to deal with mathematical history? According to some, little more is required than what was known to the authors one plans to write about; some go so far as to say that the less one knows, the better one is prepared to read those authors with an open mind and avoid anachronisms. Actually the opposite is true. An understanding in depth of the mathematics of any given period is hardly ever to be achieved without knowledge extending far beyond its ostensible subject-matter. More often than not, what makes it interesting is precisely the early occurrence of concepts and methods destined to emerge only later into the conscious mind of mathematicians; the historian's task is to disengage them and trace their influence or lack of influence on subsequent
		developments.
	\end{quote}
	When Weil writes ``interesting'' here, he means it, I believe, in the same sense I mean it: interesting to other mathematicians.  Corry writes about Weil's ICM discussion: ``Clearly, his main point was not to discuss the “why” or “how,” as his title had it, but rather the “who”'' \cite[p. 200]{CorryPoetics}.}  There is undoubtedly a question of audience here and Rowe is slightly more circumspect in his discussion:
	\begin{quote}
		No one, it seems to me, ought to dismiss the practical advantages of Weil's approach to history. True, it produces a highly rationalized picture of past achievements, not to mention a contracted image of how mathematics developed, but it also addresses the interests of the audience for which it is intended. Moreover, the historical studies of Weil, Dieudonn\'e and others in the Bourbaki tradition represent serious and very valuable contributions to scholarship. 
	\end{quote}
	Grattan--Guinness introduces a nice terminological difference between two styles that he calls {\em history} and {\em heritage} \cite{GrattanGuinnessHistvsHeri}: 
	\begin{quote}
		History addresses the question “what happened in the past?” and gives descriptions; maybe it also attempts explanations of some kinds, in order to answer the corresponding “why?” question...History should also address the dual questions “what did not happen in the past?” and “why not?”; false starts, missed opportunities, sleepers, and repeats are noted and maybe explained...Heritage addresses the question “how did we get here?,” and often the answer reads like “the royal road to me.” The modern notions are inserted into N when appropriate, and thereby [a notion] is unveiled (a nice word proposed to me by Henk Bos): similarities between [a notion] and its more modern notions are likely to be emphasized; the present is photocopied onto the past.''
	\end{quote}
	He then states that: ``{\em Both kinds of activity are quite legitimate}, and indeed important in their own right; in particular, mathematical research often seems to be conducted in a heritage-like way'' (emphasis in the original), before asserting that ``taking heritage to be history'' is frequently the mathematician's view).  This does seem an apt description of mathematical research: we choose to read what we want into the mathematics of the past, sometimes simplifying, sometimes highlighting aspects we deem important before putting these notions to work in the sense we need.

As an amusing sidenote, Corry also remarks that many glowing reviews of the Bourbaki texts were written by members of Bourbaki themselves \cite[p. 299]{Corry}.  Pierre Samuel's review of Bourbaki's algebra text, included the following line: ``As Thucydides said about his {\em History of the Peloponnesian War}, this is...a treasure valuable for all times (I.22).''  Samuel was a member of Bourbaki at the time of writing.  One can only wonder if this is irony.\endnote{Thucydides himself has been criticized as ``an artist who responds to, selects and skillfully arranges his material, and develops its symbolic and emotional potential'' \cite[pp. 231-32]{ConnorThucydides}.  Connor also writes  \cite[p. 233]{ConnorThucydides}: ``Thucydides authority can intimidate, especially at this remove when so much evidence has disappeared, when alternative and dissenting versions of events have been lost, when many controversies of his age have been forgotten.'' This sentiment seems almost equally applicable to Bourbaki's history project.}   \endnote{One might also lob the criticism of Thucydides above at Weil himself; Weil mentions Thucydides twice in his autobiography \cite[p. 27 and 54]{WeilAuto}. At the beginning of his 1950 ICM address \cite{WeilICM1950}, Weil writes:
\begin{quote}
 There appears to have been a certain feeling of rivalry, both scientific and personal, between Dedekind and Kronecker during their life-time; this developed into a feud between their followers, which was carried on until the partisans of Dedekind, fighting under the banner of the ``purity of algebra'', seemed to have won the field, and to have exterminated or converted their foes.  Thus many of Kronecker's far-reaching ideas and fruitful results now lie buried in the impressive but seldom opened volumes of his Complete Works. While each line of Dedekind's XIth Supplement, in its three successive and increasingly ``pure'' versions, has been scanned and analyzed, axiomatized and generalized, Kronecker's once famous Grundzuge are either forgotten, or are thought of merely as presenting an inferior (and less pure) method for achieving part of the same results, viz., the foundation of ideal-theory and of the theory of 
 algebraic number-fields. 	
\end{quote}
Is there any doubt that this is written to elicit sympathy for a historical figure or with whom Weil's sympathies lie?}  

To the extent that mathematicians read the history of their subject at all, I think this kind of rational reconstruction is problematic {\em even for mathematicians}. Indeed, this view of history retroactively legitimizes some points of view while de-legitimizing and circumscribing others.  Undoubtedly it is useful to contemporary mathematicians to have some {\em organized} form of the knowledge of the past that is easily accessible to them, but having made those choices what are to we make of {\em using} that organized knowledge for purposes of justification of current mathematical work?  
\end{entry}

As regards Weil's primacy in theory of algebraic fiber bundles, most of the historical sources available to us seem like breadcrumbs laid out by Weil himself (or, perhaps, those who are closely allied with his views).\footnote{Weil certainly viewed himself as a protagonist in mathematical history.  For a discussion of Weil's role in shaping the history around the Riemann hypothesis for function fields of curves over finite fields, see ``La guerre de recensions'' by \cite{Audin}.}  However, there are a few other leads we may chase: we can start looking through MathSciNet or Zentralblatt for articles mentioning terms like analytic fiber bundle, algebraic fiber bundle, algebraic bundle, or their French or German equivalents in the period before 1955 or so.  There are very few such papers.  The first few I could find all appear in the 1950s. Besides Chern's theory of characteristic classes \cite{Chern}, the following papers seem notable:
\begin{enumerate}[noitemsep,topsep=1pt]
\item N. Hawley wrote \cite{Hawley}, a student of S. Bochner wrote an article entitled ``Complex fiber bundles'' delivered to the National Academy in 1952; this paper mentions analytic line bundles on Riemann surfaces and proposes an ``analytic homotopy'' classification along the lines of Steenrod.  He also mentions connections with some results of Hirzebruch on projective line bundles over the projective line on what are now called Hirzebruch surfaces.  His analysis was mentioned in work of 
\item Atiyah \cite{Atiyah}, where it is remarked that these results are in error; this paper explicitly mentions Weil's work in its discussion of complex analytic bundles.  Atiyah's paper also points back to the work of Weil, via the theory of non-abelian theta functions mentioned above.
\item Kodaira and Spencer \cite{KodairaSpencer} directly mention the Cartan seminar in their treatment of complex analytic bundles in terms of sheaf cohomology and the work of Weil on Picard varieties; these papers will be discussed in detail later.     
\item Nakano \cite{Nakano} discusses complex analytic fiber bundles as well, but he explicitly makes mention of Weil's 1949 Chicago colloquium presentation.
\end{enumerate}

\subsubsection*{Conclusion}
Putting all this together, I think we can attribute the idea of studying fiber bundles in algebraic geometry to Weil, and assert that Weil influenced the treatment of complex analytic fiber bundles by means of his highlighting problems of interest to Cartan in this language.  Weil was very well-respected among mathematicians of his age and those of an older generation, but for his ideas to gain a foothold, those ideas would have to be adopted by the next generation of researchers; this distribution seems to have been achieved through the Cartan seminar.  

By 1952, Weil had given a course on algebraic fiber bundles notes from which were widely cited \cite{WeilFBCourse}, the Cartan seminar had, after discussing cohomology in various contexts, moved back to the theory of analytic spaces and had considered fiber spaces, in particular from a cohomological point of view.  The merged Cartan--Weil point of view on fiber bundles and fiber spaces was thus internalized by the participants of the S\'eminaire Cartan.  I think it's safe to say that Weil's point of view on analytic and algebraic fiber bundles was an integral part of the mathematical culture around the Cartan seminar by 1952/3, and by May of 1953, Serre gave a summary of this work in a S\'eminaire Bourbaki talk \cite{SerreFSBBKI}.  In the penultimate section of this note, concerning examples, Serre explains that ``the most important are vector bundles''\footnote{Serre writes: ``Les plus importants sont les espaces fibrés à fibre vectorielle.''}, a point to which we will return.  


\subsection{The homological style in algebraic topology}
\label{ss:homologicalalgebra}
\addtoendnotes{\vskip .2em\noindent {\bf The homological style..} \vskip .2em}
Eilenberg grew up in the Polish school of mathematics centered in Warsaw where he was a student of Karol Borsuk \cite[p. 1346]{EilenbergObit}.  I'd like to spend a little bit of time tracing his role through the foundations of algebraic topology leading up to Cartan--Eilenberg's {\em Homological algebra}.  A great many things have been written about the history of homological algebra (for example, see \cite{WeibelHHA,McClearySS} etc.) so I'd like to weave a path that avoids much that has been written (or, at least those things of which I'm aware).

The Warsaw school, which flourished between the two World Wars, published in the journal {\em Fundamenta Mathematica} founded in the 1920s.  Led by Sierpiński, and including many notable figures such as Kuratowski and Borsuk, the interests of the Warsaw school centered around point-set topology, real analysis, logic and set theory.  This school, which was firmly grounded in axiomatic approach to mathematics, took inspiration from the work of Felix Hausdorff as well as a ``Russian modernist'' tradition as argued by Corry \cite[p. 102]{CorryRussia}.  Point-set topology in style here was an undoubtedly modern field in the sense of the word used in the introduction, growing from the work of P. Alexandroff and A. Urysohn who were interested in studying ``general'' topological spaces and had analyzed various corresponding ``classification'' questions (the list of these names brings to mind, e.g., separation axioms).  Sierpiński was trained by Luzin in Moscow, bringing these ideas to Poland and the Warsaw school and the Russian school had close contact up until the 1930s.  

In this section, we aim to situate the vast numbers of points of view that were introduced into algebraic or combinatorial topology in the late 1930s and 1940s in response to Poicnar\'e's early work.  Of particular interest to us for our later discussion will be how this material is reflected in Eilenberg's aesthetics and mathematical methodology as this provides a frame for Cartan and Eilenberg's text.  Moreover, I aim to indicate in this section and partially the next one, how certain algebraic ideas were adopted by the core of the group of combinatorial topologists.



\begin{entry}
	\label{par:Eilenberg}
Samuel Eilenberg, frequently referred to as Sammy, received his degree in 1936, but had already published a number of papers about e.g., continuous maps to spheres (in the direction of Borsuk's cohomotopy), Lusternik--Schnirelmann category, etc.  These papers show increasing interest in the objects of combinatorial (now algebraic) topology, as opposed to the general topology more stereotypically associated with the Polish school.  Eilenberg came to the United States in 1939, first to Princeton where Oswald Veblen and Solomon Lefschetz were welcoming and helping to place refugee mathematicians fleeing WWII.  He settled in 1940 at the University of Michigan with its vibrant topology group built by R. Wilder (see \cite[p. 196-7]{Duren3} for Wilder's recollections of Eilenberg's hiring).  Eilenberg's stay in Michigan overlapped with that of Norman Steenrod who arrived in 1942. This period was very fertile, laying the groundwork for Eilenberg's work with Steenrod, as well as Mac Lane (note Wilder writes in \cite[p. 196]{Duren3} that this was around 1947, contradicting recollections of others including Mac Lane).  


Much has been written about Eilenberg's mathematical style, and aspects of this style will play a significant role in the subsequent discussion.  To frame this style we begin by looking at Eilenberg's early interaction with Solomon Lefschetz, especially around the foundations of singular homology; Mac Lane characterizes this interaction as an ``argument''.  According to Mac Lane, Eilenberg found Lefschetz's book \cite{LefschetzAT} ``obscure''.  Others have been more explicit in their description of the obscurity, characterizing it in terms of purported defects in Lefschetz's approach to  singular homology (specifically, the chain groups defining the theory were not free).  Gian-Carlo Rota describes \cite[pp. 18-19]{Rota} Lefschetz's book in the following way:
\begin{quote}
	This book, whose influence on the further development of the subject was decisive, hardly contains one completely correct proof.  It was rumored that it had been written during one of Lefschetz's sabbaticals away from Princeton, when his students did not have the opportunity to revise it and eliminate the numerous errors, as they did with all of their teacher's other writings.
\end{quote}
Whitney's review of \cite{LefschetzAT} for Math Reviews includes the phrase: ``In each proof the necessary elements are brought in, but sometimes in a somewhat displaced order, making the logical structure of the proof more difficult to ascertain.''

Eilenberg's style contrasts rather starkly with that of Lefschetz who is frequently described as a ``purely intuitive mathematician.''  Perhaps Mac Lane's characterization of this interaction as an argument stems from another description of Lefschetz: ``He despised mathematicians who spent their time giving rigorous or elegant proofs for arguments which he considered obvious'' \cite[p. 19]{Rota}.  Eilenberg, cognizant of deficiencies in Lefschetz's treatment of singular homology prepared his own, appearing as \cite{EilenbergSing}.  Even if this paper was initially not well received, the appearance of this paper in Annals of Mathematics can be viewed as Lefschetz's tacit endorsement (see the discussion in Section~\ref{ss:categories} for further explanation).


According to Mac Lane \cite[p. p. 1346]{EilenbergObit} Eilenberg operated under a general principle to ``dig deep and deeper till he got to the bottom of each issue''.  Hyman Bass made more explicit comparisons \cite[p. 1352]{EilenbergObit} about Eilenberg's methodology:
\begin{quote}
	He fit squarely into the tradition of Hilbert, E. Artin, E. Noether, and Bourbaki; he was a champion of the axiomatic unification that so dominated the early postwar mathematics. His philosophy was that the aims of mathematics are to find and articulate with clarity and economy the underlying principles that govern mathematical phenomena.
\end{quote}

To explore the principles that informed the above ``clarity'' (clear to whom and in what sense?), let us visit the preface of \cite[p. ix-x]{EilenbergSteenrod} where one finds:
\begin{quote}
We were faced with the problem of presenting two parallel lines of thought. One was the rigorous and abstract development of the homology groups of a space in the manner of Lefschetz or Čech, a procedure which lacks apparent motivation, and is noneffective so far as calculation is concerned. The other was the non-rigorous, partly intuitive, and computable method of assigning homology groups which marked the early historical development of the subject. In addition the two lines had to be merged eventually so as to justify the various computations. These difficulties made clear the need of an axiomatic approach.
\end{quote}
Axiomatics are justified as a way to tame the plethora of different approaches to (co)homology at the time, a topic which we will revisit shortly.  They describe their approach to axiomatization in the following way: 
\begin{quote}
	No motivation is offered for the axioms themselves. The beginning student is asked to take these on faith until the completion of the first three chapters. This should not be difficult, for most of the axioms are quite natural, and their totality possesses sufficient internal beauty to inspire trust in the least credulous.
\end{quote}
Here ``clear'' seems to be used in opposition to ``non-rigorous'' but, more importantly, also in partial opposition to ``intuitive''.  The axioms are asserted to be ``natural'', but ``natural'' with respect to what assumed background?  This is a thoroughly modernist position in the sense used by Corry \cite{Corry}.
\end{entry}

Andr\'e Weil first met Eilenberg in 1944, around the time of the publication of the singular homology treatise, but apparently after Eilenberg and Steenrod had worked out part of their axiomatic approach to homology.  Weil writes \cite[p. 526]{WeilII}:\footnote{par exemple Eilenberg me fit part en 1944 de la theorie axiomatique qu'il venait de batir avec Steenrod...} ``for example Eilenberg told me in 1944 of the axiomatic theory that he had just built with Steenrod...''  Weil was already familiar with some of Eilenberg's work of this period as he writes to Cartan about Eilenberg \cite[p. 93]{CartanWeil}:\footnote{``Le volume III de la topologie nous para{\^i}t bon, mais nous avons r\'edig\'e un appendice sur la topologie du plan, d'apr\`es les m\'ethodes d'Eilenberg, qui y trouverait bien sa place et que nous vous ferons parvenir bient{\^o}t.''}
\begin{quote}
  Volume III of the topology seems fine, but we have written an appendix on the topology of the plane according to the methods of Eilenberg, which would find its place there and which we will send to you soon.
\end{quote}
Later in the same letter \cite[p. 94]{CartanWeil} he writes, alluding to Eilenberg's treatment of singular homology:\footnote{Eilenberg a aussi une th\'eorie de l'homologie que nous sera bien utile pour nos volumes de topologie combinatoire.}
\begin{quote}
Eilenberg also has a theory of homology which will be very useful for our volumes on combinatorial topology.
\end{quote}
Weil and Cartan would eventually recruit Eilenberg to the Bourbaki fold, but Weil makes first mention of a disjunction that would play an important role in the development of homology: his disappointment that the Eilenberg--Steenrod axiomatization does not apply to de Rham theory.

\subsubsection*{(Co)homology from Poincar\'e to Eilenberg--Steenrod}
To put ourselves in the right frame of mind, let's look at the state of (co)homology in 1942.  I would like to visit the story as it is presented by Eilenberg--Steenrod, highlighting aspects of the outline they provide \cite[p. viii]{EilenbergSteenrod}.  This treatment will, I think, look rather different from many treatments in the mathematical literature, e.g., those appearing in \cite{James} or \cite{DieudonneHAG} since my goal will be to highlight how mathematical historical presentations can depend on dogmatic views.  The same structure is employed by Steingart \cite[p. 44]{Steingart}.\footnote{Steingart refers to the outline presented in the Introduction to \cite{EilenbergSteenrod} as a ``{\em historical} development'' of algebraic topology (emphasis in the original).  I quibble with this assessment as they make no such claim, which raises the question: who gets to determine what is and is not history or a historical development?}

\begin{entry}[Poincar\'e: from manifolds to homology]
	\label{par:Poincarehomology}
By his admission, Poincar\'e was interested in ``manifolds'' \cite{PoincareI}.  Because it will be useful to situate our discussion, let us recall some ways in which Poincar\'e thought about manifolds; some points discussed here were already mentioned in Paragraph~\ref{par:fiberbundles}.  We refer the reader to \cite{Scholz} for more general history on the evolution of the manifold concept and \cite[\S 3.2]{LefschetzEarly} specifically of Poincar\'e's treatment.

Poincar\'e used at least two different conceptions of manifolds within {\em Analysis Situs}: the description from \cite[\S 1]{PoincareI} corresponds, in modern terminology, to something like an embedded $C^1$-manifold, possibly with boundary, but the examples he gives are all specified in terms of equations and inequalities.  Later in the paper, he uses a version in terms of something like an open covering by overlapping parametric $n$-cells.  He asserts that the second definition is equivalent to the first one, but frequently falls back upon the first; the treatment of the second definition is undoubtedly influenced by Poincar\'e's understanding of analytic continuation of functions in the complex plane as this point of view is prominent in Riemann's conception of surface.

Poincar\'e also gave various constructions of manifolds, say in terms of ``identification spaces (think of the torus as built from the square by identifying opposing edges).  The word homeomorphism is also introduced for the first time in \cite[\S 2]{PoincareI}; unlike the modern use of the word which relies only continuity properties, Poincar\'e's version implicitly contains some differentiablity hypotheses.  Generalizing what is now called the genus of a compact orientable surface, Poincar\'e first defined numbers along the lines of Betti.  These numbers were themselves extracted from suitable ``incidence'' matrices that were themselves built out of some kind of ``combinatorialization'' of a manifold.  Poincar\'e's combinatorializations were obtained by decomposing his manifolds into something like cells homeomorphic to simplices while keeping track of orientations of some sort.  The incidence matrices Poincar\'e wrote down corresponded to linear equations from whose solution sets Betti numbers could be extracted.  By introducing and considering the ``dual'' of a given decomposition, Poincar\'e established a duality result: ``the Betti numbers equally distant from the extremes are equal'', this is the primordial version of what is now called Poincar\'e duality.



If the description given sounds vague to a modern reader, perhaps it strikes the right tone as the ``combinatorialization''  seems baroque to modern eyes.  The ``pieces'' into which one decomposes a manifold were themselves manifolds: Poincar\'e writes \cite[\S 5]{PoincareTrans}
\begin{quote}
	Consider a manifold $V$ of $p$ dimensions; now let $W$ be a manifold of $q$ dimensions ($q \leq p$) which is a part of $V$. We suppose that the boundary of $W$ is composed of $\lambda$ manifolds of $q - 1$ dimensions..
\end{quote}
He is thus {\em not} constructing implicitly or explicitly, what we might call a triangulation or simplicial approximation.

Poincar\'e's description reads (to me) somewhat like a cookbook, or perhaps more generously, an algorithm.  Readers at the time found it vague as well.  Heegard, in his thesis, observed some mistakes in Poincar\'e's discussion of orientations leading to an error in his proof of the duality theorem.  To respond to this and subsequent concerns about his treatment, Poincar\'e would eventually publish five supplements to his initial treatise.  In \cite{PoincareII}, Poincar\'e's linear equations were formulated in terms of integer-valued ``incidence matrices'', with better treatment of orientations.  Extracting normal forms for these matrices, e.g., Smith normal form, Poincar\'e's second formulation led to two types of numbers: Betti numbers and so-called torsion numbers.\footnote{My treatment here has been informed by \cite{PoincareTrans} and \cite{Stillwell}.  This translation of the computational problem to one of normal forms of matrices remains in modern computational approaches to homology.}  

Nevertheless, Poincar\'e's analysis leads to many questions.  The discussion of ``combinatorialization'' given by Poincar\'e demonstrates that it was something of an art.  Looking back one can see a number of questions.  For identification spaces that Poincar\'e considered the procedure for extracting Betti numbers was intuitive.  But even for spaces defined by equations, how could one know that a given ``combinatorialization'' was fine enough to reflect the geometry of the space under consideration?  If one chose a different combinatorialization, why didn't the Betti or torsion numbers change?  As notions of ``space'' evolved and became more general, could one still attach Betti numbers and torsion numbers? 
\end{entry}


Because it has some bearing on the timeline, we mention that Poincar\'e was specifically interested in three-manifolds and his interest here came from analysis of differential equations; it was investigations of closed curves in the context of differential equations that led to Poincar\'e's introduction of the fundamental group.  In Poincar\'e's conception: Betti numbers and fundamental groups were of secondary importance to what we would now call the homeomorphism problem for three manifolds.  Poincar\'e's treatment of homology is roughly the state of the art for the next 20-25 years. 

Instead, this period sees intense development of notions of space: bridging metric conceptions of geometry and the manifold conception, the notion of general topological space is formulated in this period.  Brouwer introduces the notion of homotopy around 1912 in the course of his analysis of classification of continuous maps.  There is a slow refinement of the process of combinatorialization as can be seen by consulting Veblen's treatment of Analysis Situs \cite{VeblenAS}.  Indeed, even in 1923/24 K\"unneth published his famous formulas for the homology of product spaces essentially using Poincar\'e's formulations \cite{Kunneth1,Kunneth2} and Lefschetz was still using these formulations in 1925.  

Volkert has described the importance of the homeomorphism classification problem for $3$-manifolds in particular, beginning in the mid 1880s \cite[\S 2.3.1]{Volkert}.  Indeed, Betti numbers and fundamental groups, and later torsion numbers were secondary and built arguably in service to the classification problem.  Volkert goes on to argue that beginning in the 1920s the invariants themselves became the focus.  

\begin{entry}[Enter abstract algebra]
	\label{par:Noetherandgroups}
Around 1925, the character of the discussion of homology changes and much of the subsequent discussion centers around the following note of E. Noether \cite{Noether1925} establishing that ``elementary divisor theory'' could be interpreted in terms of group theory.  Noether wrote:\footnote{
``Die Elementarteilertheorie gibt bekanntlich fur Moduln aus ganzzahligen Linearformen eine Normalbasis von der Form ($e_1 y_1$, $e_2 y_2$, ... , $e_r y_r$), wo jedes $e$ durch das folgende teilbar ist; die $e$ sind dadurch bis aufs Vorzeichen eindeutig festgelegt. Da jede Abelsche Gruppe mit endlich vielen Erzeugenden dem Restklassensystem nach einem Solchen Modul isomorph ist, ist dadurch der Zerlegungssatz dieser Gruppen als direkte Summe gro.Bter zyklischer mitbewiesen. Es wird nun umgekehrt der Zerlegungssatz rein gruppentheoretisch direkt gewonnen, in Verallgemeinerung des fur endliche Gruppen iiblichen Beweises, und daraus durch Ubergang vom Restklassensystem zum Modul selbst die Elementarteilertheorie abgeleitet. Der Gruppensatz erweist sich so als der einfachere Satz; in den Anwendungen des Gruppensatzes - z. B. Bettische und Torsionszahlen in der Topologie - ist somit ein Zuriickgehen auf die Elementarteilertheorie nicht erforderlich.''}
\begin{quote}
	As is well known, elementary divisor theory gives modules of integer linear forms a normal basis of the form ($e_1 y_1$, $e_2 y_2$, ... , $e_r y_r$), where each $e$ is divisible by the following; The $e$ are thus uniquely defined up to the sign. Since every Abelian group with a finite number of generators is isomorphic in the residue class system to such a module, the decomposition theorem of these groups is thereby proven as a direct sum of cyclic ones. Conversely, the decomposition theorem is now obtained directly using purely group theory, in a generalization of the usual proof for finite groups, the elementary divisor theory is derived from it by moving from the residue class system to the module itself. The group theorem turns out to be the simpler theorem; in the applications of the group theorem - e.g. Betti numbers and torsion numbers in topology - it is therefore not necessary to go back to elementary divisor theory.
\end{quote}
This comment has served as a turning point in the discussion of group theory in topology, and its significance is perhaps generational in a sense that I will now explain.  

Several authors have pointed to Noether as the source for the use of group theory in homology.\footnote{Perhaps we should be even more precise here and discuss the use of group theory in homology {\em in Germany}.  Indeed Tucker writes in his review of Seifert and Threlfall's famous topology book \cite{TuckerSFBR}: ``The influence of the German algebraic school is reflected in a frank use of group theory.''}  Lefschetz, writing in \cite[\S 3.1]{LefschetzEarly} qualifies this by writing that group theory was introduced ``by Emmy Noether (through Alexandroff).''  This period is also marked by the appearance of the first edition of van der Waerden's text \cite{vdWI,vdWII}.  The preface of van der Waerden's text cites a number of sources attesting to the prevalence of algebraic ideas; these ideas are fleshed out in van der Waerden's own recollections of the writing of his textbook \cite{vdWModAlghistory}.  Indeed, van der Waerden cites E. Steinitz's work from 1910, augmented by a treatise written with R. Baer and H. Hasse, as well as lectures of Emil Artin in Hamburg on Algebra, a seminar by Artin, W. Blaschke, O. Schreier and himself on the theory of ideals, and lectures of Emmy Noether in G\"ottingen in 1924/25 on group theory and hypercomplex numbers as key sources used in the preparation of his landmark treatise.  All of this is a testament to the prevalence of group-theoretic and, more generally, algebraic ideas in Germany during this period. 

Dieudonn\'e \cite{DieudonneNoether} wrote that Noether had persuaded Hopf to use group theory in his discussion of homology using the above ideas and thus credits the transition from matrices to group theory to Noether.  This version of events is repeated in other ``histories'' including \cite[\S 1.3]{WeibelHHA}.  However, the situation appears more complicated.  Mac Lane argues that by 1926 Vietoris \cite{MacLaneVietoris} was already using groups in homology, independently of Noether's influence.\footnote{Mac Lane qualifies this slightly.  He first notes that Dieudonn\'e's history makes the claim that ``Mathematical ideas originated in G\"ottingen and then spread to lesser places. It was not always that simple.''  After discussing the work of Vietoris in more detail he writes: ``Probably ideas passed back and forth in both directions between Vienna and Göttingen (and Berlin, Moscow and Paris).''  Granted the frequent interactions and letters that passed back and forth and the general interconnectedness of the era, this seems more plausible to me than unidirectionality.}
  
We won't quibble over dates, but the remarks of Vietoris on Mac Lane's manuscript seem rather more important in this discussion.  Vietoris writes (in Mac Lane's translation)
\begin{quote}
	Without doubt H. Poincar\'e and his contemporaries knew that the Betti
	numbers and the torsion coefficients were invariants of groups, whose elements were (classes of) cycles under the operation of addition...Then one worked with the numerical invariants rather than with the invariant groups. That was a matter of ‘taste’.
\end{quote}
Dieudonn\'e's MathSciNet review of Mac Lane's paper contains this response:
\begin{quote}
	The author quotes a letter from Vietoris, who attributes to Poincaré the knowledge of the relation of his ``incidence matrices'' with abelian groups. The reviewer has never seen any statement from Poincaré which would support this idea; Poincaré was dealing with groups in almost every part of mathematics; if he had seen groups in homology, he certainly would have said so.
\end{quote}
Even though authors of the time only used the terminology abelian groups, Dieudonn\'e insists in \cite[p. 38]{DieudonneHAG} on calling them $\Z$-modules...certainly non-standard.\footnote{The systematic adoption of module-theoretic terminology was definitely part of Bourbaki's approach to algebra and sped up by the influence of these texts.}  He furthermore writes that Noether was: ``engaged in the process of liberating linear algebra from matrices and determinants.''  Noether's approach was undoubtedly spiritually closer to Dieudonn\'e's aesthetics, and arguably the Bourbaki abstraction aesthetic than Poincar\'e's approach, which we might now-a-days call ``down to earth.''  Certainly the lack of rigor which we now associate with Poincar\'e seems antithetical to the Bourbakiste ideal. 

One can see in this discussion something of a generational difference about acceptance of new mathematical terminology or ideas. Indeed, Vietoris finished his degree in 1920, while Hopf finished his degree in 1925.  One might argue that Hopf ``grew up'' with the new kinds of abstraction and accepted and internalized the inclusion of abstract groups, making little mention of its use going forward.\footnote{Mac Lane reports \cite{MacLaneVietoris} that Hopf explicitly comments that Noether's ideas provide a ``simplification'' of the proof.  As usual, this simplification comes at the cost: that of familiarity with and facility with use of a mathematical structure, i.e., the notion of group.}  In contrast, Vietoris explicitly refers to the use of abstract groups as a matter of taste.\footnote{See \cite[\S 6]{Volkert} for a more detailed timeline and much additional discussion of the algebraisation of topology in this period.}

Werner Mayer (the same Mayer of the Mayer--Vietoris sequence) seems to be first person to take the next leap in two senses: (i) he wrote down the notion of a (chain) complex \cite{MayerI,MayerII}, though note that at this stage the notion was still viewed as provisional, and (ii) provided an axiom system for defining homology \cite[I. Abschnitt]{MayerI}.\footnote{V. Katz points out in \cite[p. 120]{Katz} that this kind of axiomatic development of homology long predates that of Eilenberg--Steenrod, whose axiomatic treatment of homology wasn't published until February 1945 \cite{EilenbergSteenrodPNAS}.  Certainly the use of ``Axioms'' was ``in the air'' and ``a matter of taste''.  Nevertheless, Dieudonné's treatment in \cite[IV \S 6 B p. 107]{DieudonneHAG} makes no mention of Mayer's paper.}  Mayer explicitly thanks Vietoris for help with topology.  Dieudonn\'e writes \cite{DieudonneNoether} that Mayer ``did not mention Emmy Noether at all in this paper. However, by that time the spirit of ‘modern algebra’ had spread to many German universities..”''  Mac Lane argues \cite{MacLaneVietoris}, contra the view espoused by Dieudonn\'e, that Mayer's conception of group theory arose from interactions with Vietoris.  Once again, Mayer's Ph.D. (granted in 1912) predates 1925, and this seems to support the view stated by Vietoris that ``groups were in the air'' and their use was ``a matter of taste.''

By the end of the 1920s, something about language seems to have shifted: most papers using homology are phrased in terms of groups.  Moreover, this change  was {\em not} viewed as inevitable even by practitioners; Noether's point of view comes out the winner, and it seems that as one of the most prominent practitioners, retrospectively she is identified with the theory.  The name complex itself seems to crystallize from the various combinatorial approximations used to build spaces.  
\end{entry}

\begin{entry}[Mathematical dead-ends]
	\label{par:deadends}
Somewhat later, Mayer analyzed variants of what would become our modern definition of chain complex.  In a chain complex, the differential squares to zero, and Mayer considers situations where these might not be the case \cite{Mayer}.  Instead, he considers sequences of abelian groups together with sequences of maps between them having the property that the composite of $p$ sequential terms is zero, as opposed to just $2$.  Eventually, Spanier showed \cite{SpanierMHT} that this notion did not give ``new invariants'' relative to the old method, and not much is said about the theory subsequently.  This gives an example of what one might call a mathematical dead-end.  I have a vague recollection that in the years since, Mayer's definition has been ``reinvented'' several times, though I can no longer recall precise sources.  
\end{entry}
  
\begin{entry}[A mangle of theories]
	\label{par:manyhomologies}
	Another strand of development after the general acceptance of group theory arose from the various studies of ``combinatorialization'' of spaces.  I will not try to give any precedence to ideas here, but I want to highlight some points of discussion that will be useful later.  
	
	Poincar\'e's own procedure for computing homology gave rise to what one might call primordial cell-decompositions, and first developments of combinatorialization followed this avenue of thought exploring variants of this notion leading to various notions of what one would eventually call a ``cell complex''.  As regards the word ``complex'' Dieudonn\'e writes \cite[p. 37]{DieudonneHAG}:
	\begin{quote}
		Unfortunately that word is given different meanings by the mathematicians who use it; for the sake of clarity we shall use a terminology that distinguishes these meanings, even if it does not coincide with the one used in the papers we describe. 
	\end{quote}
	Dieudonn\'e frames the development differently: he begins with ``decompositions of manifolds'' into what we might now call finite cell complexes via ``triangulations''.   Note that Poincar\'e himself never mentions ``triangulations'', which to a modern reader mean something else, but instead speaks of ``polyhedra''.  The notion of polyhedral subdivision changes as one can see in \cite{VeblenAS}.  All of these constructions were finite in nature.  
	
	In his 1932 thesis, Albert W. Tucker had created a combinatorialization of a manifold \cite{TuckerI} with the goal of more clearly formulating Poincar\'e's duality theory.  Later, he introduced a more general notion of cell space in \cite{Tucker2} which greatly expanded the kinds of topological spaces to which homological techniques could be applied: these developments permitted analysis of complexes that were not necessarily finite.\footnote{Since he has already been mentioned once and especially since his writings will reappear at numerous points throughout the narrative, let us mention that Albert W. Tucker was a student of Solomon Lefschetz who after his thesis joined the Princeton Mathematics Faculty in 1933 and remained there until 1974.  He was department chairman during the 1950s and 1960s and was described as the ``intellectual soul of Princeton's mathematics department'' \cite{NasarTucker}.  While his initial work was in topology, Tucker later laid the foundations for the theory of linear programming and analyzed game theory; it is these later results that are more widely known.  For example the Kuhn--Tucker conditions frequently make their way into multi-variable calculus textbooks.  Later in life, he organized an oral history project focusing on the Princeton Mathematics department which was an extremely useful resource.}

	Another variant of Poincar\'e's approach led to simplicial complexes as studied by Lefschetz \cite{LefschetzTop1930}.  This led eventually to the concept of singular simplices and singular cohomology, one treatment of which appears in \cite{LefschetzAT}; Lefschetz's particular presentation of singular homology will be revisited in Section~\ref{ss:categories}.
	
	Around the same time, Alexandroff and Čech introduce another combinatorialization by way of coverings and their intersections: this leads to the notion of the ``nerve'' \cite{CechI} and to yet other forms of homology.  The sense we want to give is of an explosion of activity and a huge number of different approaches to homology and later cohomology; we will revisit this discussion in Paragraph~\ref{par:cechhomology}.
	
	By 1935, Alexandroff and Hopf had published their famous topology book \cite{AlexandroffHopf}.  In the introduction, it is explained that preliminary versions benefited immensely from, among a whole list of mathematicians, discussions with the Princeton topology group consisting of Alexander, Lefschetz and Veblen. 
\end{entry}

\begin{entry}[Points of view]
	\label{par:pointsofview}
The year 1935 was also marked by the first international conference devoted to topology, held in Moscow.  Whitney writes \cite[p. 97]{Duren1} that the conference was notable in several ways:
\begin{quote} 
 To start, it was the first truly international conference in a specialized part of mathematics, on a broad scale. Next, there were three major breakthroughs toward future methods in topology of great import for the future of the subject. And, more striking yet, in each of these the first presenter turned out not to be alone: At least one other had been working up the same material. 
\end{quote}
The activities of the conference were summarized as well by Tucker \cite{TuckerMoscow} thus:
\begin{quote}
	The papers showed a broad conception of topology and its relations with	other subjects. The discussions which followed the papers were keen and informal; they revealed many instances of overlapping investigations, and led to an invaluable exchange of ideas. In fact one could hardly conceive of a more successful scientific gathering.
\end{quote}
The mention of the fact that a number of people were thinking about similar ideas seems notable, as does the fact that the ``successfulness'' of the conference is characterized in terms of the exchange of ideas, not, for example, in terms of what we might call a ``deliverable'' (say the solution to a problem).  

Whitney frames some aspects of the conference around the book of Alexandroff--Hopf \cite{AlexandroffHopf} that was about to appear.  The book is famously titled ``Topologie I'', with the idea that it would be the first of a three part series on the state of algebraic topology.  But \cite[p. 97]{Duren1}:
\begin{quote}
	the conference was so explosive in character that the authors soon realized that their volume was already badly out of date; and with the impossibility of doing a very great revision, the last two volumes were abandoned. 
\end{quote}
Tucker reviewed the book of Alexandroff and Hopf for the AMS \cite{TuckerBR} in 1936, writing:
\begin{quote}
	Topology is not a young subject; this year may be described as the two hundredth anniversary of its birth if we agree that it had its beginning in the problem of the seven bridges of Königsberg settled by Euler in 1736. But the systematic development of topology is new; it has only come since the work of Poincaré at the turn of the century. The International Topological Conference held at Moscow last September showed that the subject has attained a definite measure of maturity and a wide range of influence on other branches of mathematics, but that it is still undergoing rapid growth and flux. Just when topological activity seems to be slackening some new point of view sets it seething again; in the past year we have had an example of this in the ``dual cycles'' of Alexander and Kolmogoroff.
\end{quote}
The ``dual cycles'' of Alexander and Kolmogoroff are now viewed as the beginning of cohomology as a notion to be investigated systematically and distinct from homology (Poincar\'e had investigated ``dual cycles'' in his duality investigations, but this was within the confines of homology).  Striking to me about this paragraph is Tucker's phrase ``some new point of view sets it seething again''.  More explicitly: it is not problems, conjectures or theorems that are fundamental in this new ``systematic development of topology'', but ``points of view.'' Could anything sound more social than a mixing pot of ideas, where ``points of view'' are the most important thing?\endnote{Tucker closes his review with the comment:
\begin{quote}
	The book contains no real bibliography; it has merely a list of topological texts and one of works which bear directly on individual sections of the book and influenced their composition. References to the literature through footnotes have been reduced to a minimum. Certain concepts and proofs have been designated by the names of their originators, but otherwise the authors have resolutely refrained from attempting to trace notions back to their sources.  The introduction to the book embraces a short history of the development of topology and a survey of the relations of topology to neighboring branches of mathematics, as well as a discussion of the authors' program.
	\end{quote}
I get the sense that Alexandroff and Hopf just want to get back to ``proving theorems''.  Developments were too rapid and a history could be written later.  For a book that has retrospectively been so influential, this undoubtedly had the effect of burying sources of points of view.}  What informed these ``points of view''?  Why were some more readily accepted than others?  Certainly, concrete developments contributed to the conception of importance, but as I hope I have illustrated with Mayer's work, that does not appear to me to be the only factor. \todo{We need to tease out what the important ``points of view'' were here.  There were numerous approaches to topology: Tucker, Vietoris, Čech...} 

Formulations of Poincar\'e duality also led to the concept of homology or cohomology with coefficients.  But this development had other sources also.  Hopf had shown, extending ideas of Poincar\'e that an orientable manifold admits a nowhere vanishing tangent vector field if and only if its Euler characteristic is $0$.  Stiefel had extended this result to a finite number of (linearly independent) tangent vector fields using (characteristic) cohomology classes.  In an initially parallel direction, Whitney had reformulated \cite{WhitneyMaps} the Hopf classification theorem of maps from an $n$-complex to an $n$-sphere in cohomological terms, inventing the notion of cohomology with coefficients in a homotopy group along the way.  Eilenberg had extended Whitney's analysis studying obstructions to extensions of maps (up to homotopy) using cohomology with coefficients in a homotopy group \cite{EilenbergObst}.  By 1941, Steenrod extended some of Eilenberg's results and used cohomology with coefficients to analyze Stiefel's problem about linearly independent vector fields \cite{SteenrodTF}.  Eventually, Steenrod took another step and studied homology with local coefficients \cite{SteenrodLoc} more generally.  
\end{entry}

\subsubsection*{Conclusion}
We have already seen that ``points of view'' were driving forces in the early development of combinatorial/algebraic topology and it appears that problems and conjectures were important, but perhaps served more as touchstones of progress than driving factors.  Indeed, the problems of the era seem like a string of evolutions/refinements, depending on the terminology in vogue.  One begins with the Euler characteristic and its interpretations/computations: Hopf's theorem about zeros of vector fields and Euler characteristics is eventually subsumed by characteristic cohomology classes in several ways, e.g., the Gauss--Bonnet theorem in terms of curvature.    

Riemann's genus, also connected to the Euler characteristic, led in different directions: first to Betti numbers, homology groups and eventually problems around classification of manifolds (see, e.g., \cite{AlexanderVeblen}).  The classification of surfaces was ``known'' to Riemann, though see \cite[\S 2]{Volkert} for further clarification regarding the use of the word ``known'' here.  We remarked earlier on the importance of the manifold classification problem in Poincar\'e's work on three manifolds, and subsequent work.  Problems such as the Poincar\'e conjecture were undoubtedly of great interest and topologists such as J.H.C. Whitehead were quite interested in such problems (Whitehead's purported, but false, proof of the Poincar\'e conjecture appearing in 1934 \cite{WhiteheadPoincare} attests to this).  

Links between combinatorial group theory and topology were developed significantly in this period.  General discussion of these ideas of these links can be found in \cite[Chapter II]{Chandler}; particularly interesting is the work of the German school, e.g., Artin as related to fundamental groups of knot complements.  Volkert describes \cite[Chapter 6]{Volkert} the gradual shift in the period 1920-1950 away from the classification problem and toward general invariant computations; in a sense Steenrod's discussion of fiber bundles from Paragraph~\ref{par:fiberbundles} shows the pendulum starting to swing backwards. 

It's difficult, with modern prejudices, to look upon the ``mood'' of the times as regards the manifold classification problem, but it's important to note that insolubility of the word problem for finitely generated groups was still many years in the future, and the concept of algorithmic unsolvability didn't exist in any formal sense yet.  Volkert argues that the idea of classification of three manifolds was an important precursor to the development of topological ideas \cite[Chapter 3]{Volkert}.  In the 1920s and 1930s it perhaps seemed likely that a classification of ``all'' manifolds, mirroring the case of surfaces where some small collection of topological invariants could be used to enumerate all homeomorphism types, was plausible and perhaps even eventually possible.  

At the same time, higher dimensional Riemannian geometry was ``in the air'' because of Einstein's relativity theory, with which Hilbert and Poincar\'e were conversant (and Weil discusses his learning special relativity in his autobiography).  Mayer was a privatdozent with Einstein in Vienna in the 1920s before his work in topology, so people walked back and forth between these worlds: some higher dimensional notions were in the popular academic conscious at the time.  

Each of these ideas appears to follow a historical arc, but changing ``points of view'' lead the development in different directions.  On the other hand, the work of Mayer on modifications of chain complexes seems to fall outside this pattern,  It appears to be generalization for its own sake: an idea is developed with the potential of giving new invariants, but eventually does not succeed, even on its own terms.\endnote{A related direction appears in Tucker's review of two books \cite{TuckerSFBR}.  Tucker refers to the contents of Steinitz's {\em Vorlesungen Uber die Theorie der Polyeder unter Einschluss der Elemente der Topologie} \cite{Steinitz} as: 
\begin{quote}
 an isolated chapter of geometry which bears little relation to the main stream of contemporary topological research, but which stands by itself, firm in its own intrinsic worth. It deals with a question almost as old as analysis situs itself--the combinatorial classification of ordinary polyhedra. 
\end{quote}
This strikes me as encapsulating the fashions of topology at the time, at least in Princeton, rather explicitly.}  The most successful branches described above seem central in the sociological sense.  A clear case is Chern's work as it linked characteristic classes, interests of Stiefel, Whitney, Steenrod, with differential-geometric ideas of Cartan.   

The ``axiomatic'' approach of Eilenberg--Steenrod to homology thus hardly seems inevitable.  Its appearance is predicated on the vast number of different homology theories that had arisen.  By the time of its appearance, homology and cohomology had been, in some form, contained in the mathematical literature for almost 50 years.  It was also certainly not a foregone conclusion that this method would be widely accepted.  Indeed, the fact that ``combinatorial topology'' was amenable to different points of view seems to be one feature of the subject that allowed it to be receptive to the axiomatic approach in the first place.  In the next section, we will explore the limits of this receptivity.

\subsection{A cartographic view of functoriality in the mathematical landscape}
\label{ss:categories}
\addtoendnotes{\vskip .2em\noindent {\bf A cartographic view..} \vskip .2em}
By 1942 Eilenberg and Mac Lane were investigating homology of infinite complexes having moved away from manifolds to more purely ``combinatorial'' geometries and using Tucker's theory of complexes (see the references in \cite{EMinfinitecycles}).  In this paper, Eilenberg and Mac Lane analyze infinite complexes by successive ``finite'' approximations via limiting constructions: one sees the notion of direct limits appear, as well as the beginnings of universal coefficient theorems.  This paper rested on a fortuitous accident that Mac Lane has recounted with slightly varying details in several places, e.g., \cite{EilenbergObit,MacLaneEilplusCat}; this version is from \cite{EilenbergObit}.\endnote{The source of the Eilenberg--Mac Lane collaboration is recounted in other sources, with varying degrees of detail and emphasis.  For example, the treatment in \cite{McLarty} summarizes the episode thus:
\begin{quote}
	He published on several topics in his early career including logic but focused on technical problems in algebra aimed at number theory. His solution to one of these was a strange family of groups. Samuel Eilenberg knew these same groups solve a problem in topology. When Eilenberg (who, by the way, liked philosophy a great deal less than Mac Lane did) learned of Mac Lane’s result, the two of them agreed this could not be a coincidence. They set out to find the connection. 
	\end{quote}
And views this as the primordial source of category theory.  For our treatement, it will be important to separate the two conceptions.}
\begin{quote}
	 I had calculated an example of the group of group extensions for an interesting factor group involving a prime number p. When I told Sammy this result, he immediately saw that it answered a question of Steenrod about the regular cycles of the p-adic solenoid (inside a solid torus, wrap another one p times around, and so on, ad infinitum). So Sammy and I stayed up all night to find out the reason for this unexpected appearance of group extensions. We found out more: it rested on a “universal coefficient theorem” which gave cohomology with any coefficient group G in terms of homology and an exact sequence involving Ext, the group of group extensions. 
\end{quote}
Mac Lane writes that this bit of serendipity had tremendous mathematical consequences.  While others, e.g., Steenrod could likely have made the ``homological'' observation (though Steenrod was not yet at Michigan when this happened), Eilenberg's specific aesthetic tastes undoubtedly shaped the form of the final product. 

In his autobiography, Mac Lane explains the role Lefschetz played in the immediate product of the above computation \cite{MacLaneAuto}: ``Of course, we told Lefschetz about it; he at once asked us to write it up as an appendix (of five pages) to his 1942 book Topology.''  The final version of this paper appeared in Annals of Mathematics in 1942 \cite{EMExtensions}, whose transformation in the hands of Lefschetz we will explore momentarily.

The systematic analysis of properties of Ext groups and limiting constructions used above, together with conceptions from homotopy, e.g., homology and homotopy groups led to another invention: the theory of categories and functors \cite{EMNatTransf}.  This observation is underscored by Weil in his review of \cite{EMNatTransf}
\begin{quote}
	A vague idea of covariance and contravariance is often met with in group-theory, topology, etc.; that is, one feels that the character-group is contravariant to the group, that the homology and co-homology groups of a complex are, respectively, covariant and contravariant to the complex. This is of special importance in the building up of limits of direct and inverse systems (``projective'' and ``inductive'' limits) of groups, spaces, etc. The authors have succeeded in finding for this a precise definition, which is likely to be helpful in classifying and systematizing known results and also in looking for new relations between groups.
\end{quote}
Mac Lane in his review of Cartan and Eilenberg \cite[p. 622]{MacLaneCEReview} gives a rather more prosaic view: 
\begin{quote}
The authors' approach in this book can best be described in philosophical terms and as monistic: everything is unified.\endnote{Einstein was one of the first to use the terminology ``unified'' in physics viz his search, beginning in the 1920s for a theory that unified general relativity and electro-magnetism.  Higher-dimensional approaches to this problem appeared around the same time: the so-called Kaluza-Klein models.  Where did this terminology come from?  Did it create ``unification'' talk in mathematics, did it come from elsewhere?}
\end{quote}
Granted everything else that has been written about Eilenberg's point of view, and Cartan's admission that Cartan--Eilenberg was largely (entirely?) written by ``Sammy'', this description seems an apt summary.\endnote{While this has undoubtedly been discussed elsewhere, it seems to me that given the time period in which all of this mathematics was invented, this beautiful, orderly, unified structure provided a stark counterpoint to the external chaos of the world.}  

In the previous section we saw that algebraic/combinatorial topology was receptive to many different ``points of view''.  But how does a subject draw the lines that determine what work should be ``promoted'' and what should be ``ignored''?  In other words, what were the boundaries of receptivity?  Were there ideas that should be ignored?  Such question seems to have some features in common with what is usually called boundary work \cite{gieryn1999cultural}, and is the source of my use of cartography in the section heading.  

I proceed to investigate the question of drawing boundaries of mathematical disciplines; my analysis focuses on two confluent streams: the reception of Eilenberg and Mac Lane's ideas about category theory within mathematics, and mathematical publication practices in the late 1940s and early 1950s.  Among other things, I aim to argue that the propagation of category-theoretic ideas was far from inevitable.  Moreover, just as the scientific idea of mapmaking reveals evidence of social and political influence under scrutiny,\footnote{J.B. Hartley writes in \cite{Hartley}: ``...the scientific
		rules of mapping are, in any case, influenced by a quite different set of rules, those
		governing the cultural production of the map. To discover these rules, we have to read
		between the lines of technical procedures or of the map’s topographic content. They are
		related to values, such as those of ethnicity, politics, religion or social class, and they are also embedded in the map producing society at large.''} I aim to illustrate how social and political factors contributed to the visibility of category theory, not the least of which is the insistence of Eilenberg and Mac Lane and their place within American mathematical society in particular, and on the international mathematical scene more generally.  

\begin{entry}[Adopting categories and functors]
The introduction to \cite{EMNatTransf} reflects the key themes in the above discussion by beginning with a philosophical point:
\begin{quote}
	Frequently in modern mathematics there occur phenomena of ``naturality'': a ``natural'' isomorphism between two groups or between two complexes, · a ``natural'' homeomorphism of two spaces and the like. We here propose a precise definition of the ``naturality'' of such correspondences, as a basis for an appropriate general theory. In this preliminary report we restrict ourselves to the natural isomorphisms of group theory; with this limitation we can present the basic concepts of our theory without developing the axiomatic approach necessary for a general treatment applicable to various branches of mathematics.
\end{quote}
And the philosophical was evidently in the air as regards category theory.  Mac Lane writes in \cite{MacLaneCategories} that the name category was taken from Aristotle and Kant, while the notion of functor was adapted from Carnap.  

While the notion of tensor product of vector spaces was older, Whitney had only recently invented the tensor product of abelian groups \cite{Whitneytensor} (see the beginning of Section~\ref{ss:bourbakialgebra} for a much more involved treatment of the evolution of the notion of tensor product).  Weil closes his review by mentioning:
\begin{quote}
	In particular, the authors use it to derive some interesting relations concerning Whitney's ``tensor-product'' of groups, and clarify the nature of the latter.
\end{quote}
It's worth pointing out that while the tensor product of abelian groups was invented here, it was Bourbaki that led the shift to usage of modules (as opposed to the notion of groups with operators used in van der Waerden), and the tensor product of modules did not appear in print until roughly 1948 (though once again, was known and used by members of Bourbaki).

By 1944, Mac Lane had moved to Columbia as part of the US War effort to study problems in applied mathematics.  He writes: ``I took the occasion to hire a number of mathematicians well-known to me, including Samuel Eilenberg, Paul Smith, and Hassler Whitney.''  While Smith was already at Columbia in the 1940s, this seems to be the occasion upon which Eilenberg moved there.  

The more general treatment promised in the above announcement appeared in 1945 as \cite{EMCategories}.  Eckmann argues \cite[7.5]{EckmannFB} that the development of the theory of fiber bundles played a decisive role in the creation of the theory of categories and functors. While the PNAS note of Eilenberg and Mac Lane focuses on group theory, given Eilenberg's extensive work in homotopy theory before, this seems plausible.\endnote{The introduction to \cite{EMCategories} closes with the following (bold?) self-assessment of the theory:
\begin{quote}
	The invariant character of a mathematical discipline can be formulated in these terms. Thus, in group theory all the basic constructions can be regarded as the definitions of co- or contravariant functors, so we may formulate the dictum: The subject of group theory is essentially the study of those constructions of groups which behave in a covariant or contravariant manner under induced homomorphisms. More precisely, group theory studies functors defined on well specified categories of groups, with values in another such category.
	
	This may be regarded as a continuation of the Klein Erlanger Programm, in the sense that a geometrical space with its group of transformations is generalized to a category with its algebra of mappings.
\end{quote}
McLarty argues \cite[p. 217]{McLarty} that the later conceptions of categories and functors competed with Bourbaki's ideas of structure, and the latter were to some extent, abandoned.}

The story of the publication of Eilenberg--Mac Lane's category theory treatise is particularly fascinating, and versions of the anecdote that follows have been recorded in various places with varying details.  Mac Lane recalls \cite{MacLanePNAS}:
\begin{quote}
	The first of these papers is a more striking case; it introduced the very abstract idea of a “category”—a subject then called “general abstract nonsense”! When Eilenberg and I submitted a full presentation in 1945 (to the Transactions of the American Mathematical Society), we feared that the editor would turn it down as too “far out,” not really mathematics. So Eilenberg, who knew the editor well, persuaded him to choose as referee a young mathematician—one whom we could influence because he was then a junior member of the Applied Mathematics Group at Columbia University (war research), where Eilenberg and I were then also members, and I was Director.
\end{quote}
Writing in \cite[p. 130]{MacLaneDevelopments} Mac Lane adds: ``the editor of the Journal (The Transactions, AMS) was quite skeptical of its content.''  The names of both the editor and referee are revealed in \cite[p. 125]{MacLaneAuto}: the editor was none other than Paul Smith, though the ``persuasion'' aspect described above is toned down to ``discussion'', while the referee was a young student of Marshall Stone by the name of George Mackey.\footnote{This claim is contradicted in \cite[p. 62 2.3.2.1]{Kromer} where it is asserted that Smith was {\em not} an editor of Transactions at the time.} 

I want to reflect briefly on this episode.  Evidently Mac Lane has mentioned this statement in numerous, rather public, venues.  Mac Lane thus appears unashamed of the act of editorial intervention.  Perhaps this is accounted for by a certain sheepishness around the theory (``the material is trivial'' as we saw mentioned above) counteracted by a desire to publish what has since been argued to be perspective masquerading as mathematics (for example the refrain that category theory is ``abstraction for abstraction's sake'').\footnote{Mac Lane writes \cite[p. 30-31]{MacLaneCategories} ``Now the discovery of ideas as general as these is chiefly the willingness to make a brash or speculative abstraction...''}   In support of this viewpoint, I offer Mac Lane's discussion of the importance of publishing an announcement in PNAS:
\begin{quote}
	So in this case publication in the Proceedings was perhaps vital at the start; Category Theory is now accepted. In other words, without the Proceedings, this idea might well have been buried, unpublished.
\end{quote}
At the time of its invention, category theory was a sidelight to Mac Lane's mathematical work around the war effort.  In any case, I leave it to the reader to determine the comparative ethics of this episode, but the mathematical community was small and plausible readers of such cutting edge mathematics were probably hard to find. 
\end{entry}


\begin{entry}[Categorical inevitability]
	\label{par:categoricalinevitability}
Mac Lane's specific choice of language also deserves comment; he writes in \cite[p. 30]{MacLaneCategories}:
\begin{quote}
	Categories, functors, and natural transformations themselves were discovered by Eilenberg-Mac Lane [1942a] in their study of limits (via natural transformations) for universal coefficient theorems in tech cohomology. In this paper commutative diagrams appeared in print (probably for the first time).
\end{quote}
Mac Lane was a trained philosopher; he for example reviewed Carnap's work for the Bulletin \cite{MacLaneCarnap} and goes on to mention the philosophical pedigree of the word choices in category theory\footnote{In \cite[p. 30]{MacLaneCategories} he writes ``in this case supported by the pleasure of purloining words from the philosophers...''}, so I highlight his use of the word {\em discovered} rather than invented.  

Mac Lane elsewhere supports his use of the word discovered by arguing for the inevitability of category theory \cite[p. 131]{MacLaneDevelopments}.\footnote{Here, he remarks further on the philosophical pedigree of the words ``category and functor'' writing: \begin{quote}
		Since the philosopher Kant had made ample use of general categories, the term was borrowed from him for its present mathematical use, while Carnap, in his book on Die Logische Syntax der Sprachen had talked of functors in a different sense and made some corresponding mistakes. It seemed in order to take over that word for a better and less philosophical purpose.\end{quote} In his use of the word ``better'', he seems to betray his view of the respective status of philosophy vs. mathematics.}  
Explicitly: ``If Eilenberg and Mac Lane had not formulated category theory, who would have done so - or might it never have appeared?'' He goes on to mention Chevalley, Hopf, Steenrod, Cartan, Ehresmann, von Neumann and Ulam before discussing each of these figures.  

The one sentence justifications given by Mac Lane as to why this list of people could have invented category theory seem at best tenuous to me.  Undoubtedly, Chevalley made a number of category-theoretic contributions later, and I believe one can argue that category theory fits firmly within his mathematical aesthetic, about which I will have more to say later, but Chevalley's early interactions with Eilenberg, both as a member of Bourbaki and as a collaborator, call into question the source of his interest.  Moreover, as discussed amply in \cite[\S 8.5 p. 372]{Corry}, while some members of Bourbaki used category theory, there were others who were significantly less so-inclined, e.g., Weil. Bourbaki itself never adopted this formalism for its own foundational projects.  

Steenrod remarks about the paper of Eilenberg and Mac Lane that ``no paper had ever influenced his thinking more'', so it seems difficult at best and impossible at worst to disentangle Eilenberg's influence on Steenrod's papers.  However, if we look at Steenrod's published works only his joint work with Eilenberg comes close to the level of abstraction of the Eilenberg--Mac Lane category theory.  I think one can make similar arguments about Cartan and Ehresmann.  MacLane himself appears to dismiss Ulam as a plausible inventor, and as regards von Neumann, it's worth remembering that he wrote in 1947 \cite{vonNeumann}: ``at a great distance from its empirical source, or after much ``abstract'' inbreeding, a mathematical subject is in danger of degeneration'', which seems quite far from the aesthetic that other authors (even at the time) associated with category theory.  Curiously, Mac Lane leaves Grothendieck off his list, but that is perhaps another discussion.
\end{entry}



\subsubsection*{On the politics of academic publishing}
The episode involving Eilenberg--MacLane's category theory paper just recounted brings to mind another public discussion of academic publishing.  Writing in ``A beautiful mind'' Sylvia Nasar narrates \cite[p. 58-59]{NasarABM} the following about the status of the journal Annals of Mathematics under the direction of Lefschetz. 
\begin{quote}
	Entrepreneurial and energetic, Lefschetz was the supercharged human locomotive that ... pulled the Princeton department out of genteel mediocrity right to the top. He recruited mathematicians with only one criterion in mind: research. His high-handed and idiosyncratic editorial policies made the Annals of Mathematics, Princeton's once-tired monthly, into the most revered mathematical journal in the world. He was sometimes accused of caving in to anti-Antisemitism for refusing to admit many Jewish students (his rationale being that nobody would hire them when they completed their degrees), but no one denies that he had brilliant snap judgment. He exhorted, bossed, and bullied, but with the aim of making the department great and turning his students into real mathematicians, tough like himself.
\end{quote}
This discussion is relevant to us because Cartan and Eilenberg's homological algebra treatise appeared in ``Annals of Mathematical Studies'' which is viewed as the sister ``monograph'' series to Annals of Math.\endnote{The separation of the publication of Annals of Math Studies from Annals of Mathematics is discussed in some interviews with Al Tucker \cite{TuckerInterviews}.  He recalls: \begin{quote}The Annals Studies was started in a rather strange way. At that time the Annals of Mathematics had a surplus of papers, and the editors felt that they were plagued especially by long papers, papers of a hundred pages or so. At that time the Annals had a total page count for the year of perhaps 700 or 800 pages, and so two or three 100-page papers took up almost half of a year’s production. So it was decided, largely by Lefschetz, that the formalizing of the Princeton Mathematical Notes could be combined with a means of publishing long papers or perhaps monographs consisting of several papers on a single topic. And this was the reason for the name Annals of Mathematics Studies, to enable the editors of the Annals of Mathematics to transfer long papers or groups of papers to the Studies. That’s the reason for the title.\end{quote}} 

Undoubtedly, {\em now} Annals of Mathematics is widely viewed as one of the most high-status mathematical journals (though one's view of its relative rank among the most high-status journals is colored by nationality and field). Nasar's comment suggests this was not always the case, a point which is fleshed out in Karen Parshall's book \cite[pp. 241-246]{Parshall} so let's take this at face value and discuss the evolution of prestige around Annals of math from the 1920s through the 1950s, which is also tied to the rise of the United States and Princeton in particular as a center of international mathematical power.  

Both Nasar and Parshall center their discussion of the rise of prestige around Lefschetz and there are several episodes that we commingle to demonstrate Lefschetz's editorial hand.\footnote{Nasar's analysis of Lefschetz comes from various biographical sources, as well as first-hand interviews with various mathematicians who knew Lefschetz.}  Many of the papers we have discussed, e.g., the work of Lefschetz himself, Mayer, Tucker, Eilenberg, Steenrod, as well as the work of Serre (both in algebraic topology and later algebraic geometry) that we will discuss soon, all appeared in Annals of Math.  Exploring the sociology of publishing in Annals specifically and broader questions about academic publishing therefore seem natural topics for our understanding of the evolution of the theory of projective modules.

\begin{entry}[On the proliferation of mathematical journals]
Before moving forward, let us briefly recall a few things about the evolution of the culture of mathematical publication, following the discussion of Bartle \cite{Bartle}, but especially to highlight differences from present-day publication practices.  By the mid 19th century, mathematical publication had become specialized enough and there were so many journals publishing mathematics that bibliographic journals were created to organize, say thematically, what was being written, but also sometimes to critically review the contents.  Bartle suggests the most important such journal at the time was the Jahrbuch uber die Fortschritte der Mathematik (henceforth, Jahrbuch).  By the turn of the century, these review journals were already around 1000 pages themselves.  Over the next 20 years, due to numerous factors including journal delivery lag times, slowness of reviewers and external political factors such as World War I, there was increasing lag time in publication of reviews, which led to dissatisfaction on the part of users.  

In the late 1920s, to rectify the problems with Jahrbuch mentioned above, Otto Neugebauer in collaboration with the Springer publishing house conceived ZentralBlatt, the first issue of which appeared in April 1931.  By 1934, Neugebauer had emigrated to Denmark and he resigned from the editorial board of Zentralblatt in 1938.  Problems with Zentralblatt around the War, including the mass emigration of mathematicians to the United States led to the desire to recreate a reviewing publication in the US.  Neugebauer helped found Math Reviews in conjunction with the efforts of a number of mathematicians, including Veblen and Warren Weaver, who at that time was the director of the division of Natural Sciences at the Rockefeller foundation.  The fascinating history of the formation of Math Reviews is documented in \cite[pp. 69-89]{Pitcher}, but it's worth remarking that amount of mathematics research being done increased considerably during the period 1920-1950.  

All of this mathematical production led to further proliferation of journals and, by the time of the Cold War, further influence from commercial publishing houses on mathematical publication.  Organically, status hierarchies formed among the journals, and in the wake of the two World Wars, journals themselves sometimes reflected nationalistic sentiments.  Much has been written about changes in academic publishing post cold war, but I will largely not discuss this period.  Even with the proliferation of journals, mathematical publishing was (and still is) a relatively slow procedure.  
\end{entry}

\begin{entry}[The size of the mathematical establishment]
The proliferation of papers just described might lead to the impression that the academic establishment was huge.  Indeed, this is true by comparison to the size of the establishment in the preceding century, but the role of Colleges and Universities in this period was certainly not established in the way it is now (see, e.g., \cite[p. 21]{riesman2017academic}).  Already by the 1940s and 1950s, the vast majority of people working in what we might call today ``pure mathematics research'' were housed in universities.  In order to obtain teaching positions at an American university, already by this stage having a Ph.D. or equivalent qualification was expected.  From that point of view, it is useful to compare the number of Ph.D.s being produced in this era with the modern era.  

The first thing that to observe is that Academia was almost incomparably smaller than today.  From \cite[p. 79]{AmericanDoctoratesNSF} one learns that of the $39,806$ doctoral degrees granted by U.S. universities between 1920 and 1999, $35,592$ were granted after 1960.  This document provides subject-specific information, but it is unclear from the provided data whether this information is useful for 1920-1960 mathematics as the grouping ``Other mathematics'' accounts for almost 2000 of the 4214 stated doctorates.

To get more refined information about the period 1930-1960, we turn to \cite[p. 12]{AmericanDoctoratesSurvey}.  The number of mathematical doctorates conferred in 5 year blocks in this period is summarized 
\begin{center}
\begin{tabular}{|c|c|c|c|c|c|c|}
	\hline
	5-year block  &  1930-34 & 1935-39 & 1940-44 & 1945-49 & 1950-54 & 1955-59\\ 
	\# Math Ph.D.s & 398 & 380 & 362  & 470  & 1056  & 1265 \\
	\hline
\end{tabular}
\end{center}
What should be clear is that academic mathematics was considerably smaller than it is now.  I venture that mathematical centers such as Princeton had less diluted influence than they do now. 
\end{entry}

\begin{entry}[Peer-review and academic judgment]
	Histories of peer-review in science make it clear that publishing in journals was an important part of the legitimation of discoveries: while ideas could be communicated in talks and seminars, written treatments were still essential.  The culture of academic publishing, more specifically in science and even more specifically in mathematics during the period 1920-1950 bears little resemblance to system that a modern reader might know, so I take some time to highlight differences.\footnote{The interested reader can find further discussion regarding peer-review in science in \cite{ChubinHackett} or academia more generally, e.g., \cite{lamont2009professors}.  I aim to focus the discussion here on peer-review in mathematical journals.}
	
	To set the stage, we recall a now famous episode as regards a paper of Einstein--Rosen on gravitation waves sent to Physical review around 1935; this episode is analyzed by Kennefick \cite{EinsteinRosen}.  In what from a modern point of view seems like standard procedure, the editor sent the paper to a reviewer (identified as H.P. Robertson at Princeton in the paper just referenced).  The reviewer returned ten pages of comments which cast doubts on many claims made in the paper, to which Einstein replied:
	\begin{quote}
	We (Mr. Rosen and I) had sent you our manuscript for publication and had not authorised you to show it to specialists before it is printed. I see no reason to address the – in any case erroneous – comments of your anonymous expert. On the basis of this incident I prefer to publish the paper elsewhere.	
	\end{quote}
	
	Kennefick argues that editorial practices in German journals in the early 20th century were considerably less fastidious.  He goes on to discuss publication in the Prussian academy of sciences of which Einstein was a member.  From 1914 onward, Einstein was regularly called on to review submissions for this journal, and he frequently used the word ``worthless'' in his reviews.  More importantly for our discussion here, as a member of the academy, Einstein's papers were published here without question or revision.  In the US where practices were slightly different, professional status played a role in the review process.  Editors didn't want to slow the publication of work they deemed important, and frequently relied on their own instincts rather than careful review. 
	
	There are various places where one can find histories of peer review devoted to the period we are considering here, I refer the reader to \cite[Chapter 3]{GouldPR}.  The discussion of Einstein's work and journal publication sketched above seems relatively common through at least roughly 1950: ``peer review'' as such was simply a gatekeeper, and correctness of results was of secondary interest.  Editors of journals had considerable power and latitude in judgment.  Depending on the editor, peer review might consist of the editor themselves determining a paper was suitable for a journal, or it might consist of consulting other individuals for opinions. 
\end{entry}

With this background in place, we returning to the case of Lefschetz and Annals of Mathematics.  To start, we will analyze the dynamics of the editorial board of Annals of Mathematics.\footnote{There were ideas of ``top journals'': Acta Mathematica, Mathematische Annalen, and Crelle's Journal had strong reputations reflecting some of the Swedish and German schools.  In France, CRAS published many announcements, ASENS had existed for some time.  In the United states, PNAS was a venue for publishing announcements, and the AMS journals including Transactions and the Bulletin reflected some level of quality.  At the top, I believe were Annals of Math and American Journal of Math.  I am unsure what status was afforded to Duke Mathematical journal, but Whitney had published a number of papers in Duke as well.  In this landscape, Inventiones, PMIHES, JAMS didn't exist yet.}




\begin{entry}[Dynamics of mathematical editorial boards]
Lefschetz began to edit Annals of Mathematics in 1927 and continued to do so until 1958.  One can see that the composition of the editorial board changed rather significantly over the next few years, with Lefschetz as one of the few constants.  Table \ref{table:annalsofmatheditorialboard} describes the editorial board for Annals of Mathematics from just before Lefschetz became an editor to the year before the submission of Cartan--Eilenberg.  In the table, I use the following abbreviations: OSt $=$ Ormond Stone, Ale $=$ J.W. Alexander, Be $=$ A.A. Bennett, Bl $=$ H. Blumberg, Ei $=$ L.P. Eisenhart, Hi $=$ Einer Hille, Gr $=$ T.H.Gronwall, Pf $=$ G.H. Pfeiffer, Ve $=$ Oswald Veblen, Wh $=$ J.K. Whitemore, Fo $=$ W.B. Ford, Gi $=$ D.C. Gillespie, Ha $=$ W.L. Hart, He $=$ E.R. Hedrick, La $=$ R.E. Langer, P-W $=$ A. Pell-Wheeler, Ta $=$ J.D. Tamarkin, Ri $=$ J.F. Ritt, Ba $=$ H. Bateman, Va $=$ H.S. Vandiver,  Bi $=$ G.D. Birkhoff, We $=$ J.M.H. Wedderburn, Or $=$ O. Ore, Wi $=$ N. Wiener, Alb $=$ A.A. Albert, Mo $=$ M. Morse, MSt $=$ M. Stone, Za $=$ O. Zariski, Ar $=$ E. Artin, Hi $=$ T.H. Hildebrandt, Ja $=$ N. Jacobson, Mac $=$ S. Mac Lane, Why $=$ G.T. Whyburn, He $=$ M.R. Hestenes, Ste $=$ N. Steenrod, Wh $=$ H. Whitney, Bo $=$ F. Bohnenblust, Ch $=$ S.S. Chern, McS $=$ E. J. McShane, Eil $=$ S. Eilenberg, Smi $=$ P.A. Smith, Fr $=$ K.O. Friedrichs, Le $=$ N. Levinson, Sz $=$ O. Szasz, Wil $=$ R.L. Wilder, Mon $=$ D. Montgomery, Spe $=$ D.C. Spencer, Do $=$ J.L. Doob, Che $=$ C. Chevalley, Sze $=$ G. Szego, GWM $=$ G. W. Mackey, Kak $=$ S. Kakutani, Sch $=$ M. M. Schiffer, Zyg $=$ A. Zygmund, JHC = J.H.C. Whitehead, Wei $=$ A. Weil. 
\begin{table}
	\label{table:annalsofmatheditorialboard}
\begin{center}
{\bf Annals of Math: editorial board}\\
\begin{footnotesize}
	\begin{tabular}{|c|c|c|}
		\hline
		Year & Managing Editors & Assistant editors \\
		\hline
		1925-26 & Wedderburn & OSt, Ale, Be, Bl, Ei, Hi, Gr, Pf, Ve, Wh  \\
		1926-27  & Wedderburn & OSt, Ale, Be, Bl, Ei, Hi, Gr, Pf, Ve, Wh  \\
		1927-28  & Lefschetz, Wedderburn & OSt, Fo, Gi, Ha, He, Hi, La, P-W  \\
		1928-29	 & Lefschetz, Wedderburn & OSt, Fo, Gi, Ha, He, Hi, La, Ta, P-W \\
		1930  & Hille, Lefschetz &  OSt, Ri, Al, Ta, Ba, Va, Bi, Ve, Ei, We, Or, P-W, Wi\\
		1931  & Hille, Lefschetz & OSt, Ri, Al, Ta, Ba, Va, Bi, Ve, Ei, We, Or, P-W, Wi\\
		1932 & Hille, Lefschetz & OSt, Ri, Ba, Ta, Bi, Va, Or, P-W, Wi \\
		1933 & Couldn't find & -- \\
		1934 & Lefschetz, von Neumann & Alb, Bi, Hi, Mo, Or, Ri, MSt, Ta, Wi, Va, Za\\
		1935 & Lefschetz, von Neumann & Alb, Mo, Va, Ba, Or, P-W, Bi, Ri, Wi, Hi, MSt, Za, Ta \\
		1936 & Lefschetz, von Neumann & Alb, Or, Va, Ba, Ri, P-W, Bi, MSt, Wi, Hi, Ta, Za \\
		1937 & Lefschetz, von Neumann, & Alb, Mo, Va, Ba, Or, P-W, Bi, Ri, Wi, Hi, MSt, Za, Ta \\
		     & Bohnenblust   & \\
		1938 & Lefschetz, von Neumann, & Alb, Mo, Va, Ba, Or, P-W, Bi, Ri, Wi, Hi, MSt, Za, Ta \\
		     & Bohnenblust & \\
		1939 & Lefschetz, von Neumann, & Alb, Mo, Va, Ba, Or, P-W, Bi, Ri, Wi, Hi, MSt, Za, Ta \\
		& Bohnenblust & \\
		1940 & Lefschetz, von Neumann, & Alb, Ar, Hi, MSt, Ja, P-W, Bi, Mac, Why, He, Ste, Za\\
		 & Bohnenblust & \\
		1941 & Lefschetz, von Neumann, & Alb, Hi, MSt, Ar, Ja, P-W, Bi, Mac, Why, He, Ste, Za \\
		 & Bohnenblust & \\
		1942 & Lefschetz, von Neumann, & Alb, Hi, MSt, Ar, Ja, P-W, Bi, Mac*, Why, Hes, Ste, Za \\
		 & Bohnenblust & \\
		1943 & Lefschetz, von Neumann, & Alb, Hi, P-W, Ar, Ja, Wh, Bi, Ste, Why, He, MSt, Za\\
		 & Bohnenblust & \\
		1944 & Lefschetz, von Neumann, & Alb, Hi, P-W, Ar, Ja, Wh, Bi, Ste, Why, He, MSt, Za \\
		& Bohnenblust & \\
		1945 & Lefschetz, von Neumann & Ar, He, Ste, Bi, Hi, Ta, Bo, Ja, P-W*, Ch, McS, Wh, Eil, Smi, Why \\
		1946 & Lefschetz, von Neumann & Ar, Fr, Smi, Bi, Ja, Ste, Bo, Le**, Sz, Ch, Mac, Wil, Eil, McS  \\
		1947 & Lefschetz, von Neumann & Bo, Hi, Sz, Ch, Le, Wh, Eil, Mac, Wil, Fr, McS, Za, Mon, Smi \\
		1948 & Lefschetz, von Neumann & Bo, Le, Spe, Ch, Mac, Sz, Eil, McS, Whi, Fr, Mon, Wil, Hi, Smi, Za \\
		1949 & Lefschetz, von Neumann, & Hi, Bo, Le, Sz, Ch, Mac, Whi, Do, McS, Wil, Eil, Smi, Za, Spe \\
		& Montgomery, Steenrod & \\
		1950 & Lefschetz, von Neumann, & Hi, Bo, Le, Sz, Ch, Mac, Whi, Do, McS, Wil, Eil, Smi, Za, Spe\\
		 & Montgomery, Steenrod &  \\
		1951 & Lefscehtz, von Neumann, & Che, Le, Sze, Do, GWM, Whi, Hi, Mac, Wil, Kak, Sch, Za, Sz \\
		& Montgomery, Steenrod & \\
		1952(1) & Lefschetz, von Neumann, & Che, Le, Whi, Do, GWM, Wil, hi, Sch, Za, Kak, Sze, Zyg, JHC \\
		 & Montgomery, Steenrod &  \\
		1952(2) & Lefschetz, von Neumann, & Che, GWM, JHC, Do, Sch, Whi, Hi, Sze, Wil, Kak, Wei, Za, Le, Zyg\\
		& Montgomery, Steenrod & \\
		\hline
	\end{tabular}
\end{footnotesize}
\end{center}
\end{table}

Of course, we can look back and see that there are many historically excellent mathematicians present on the editorial board.  However, I want to highlight two things.  First, the composition of an editorial board is a socio-political choice: managing editors appoint subsequent managing editors and associate editors.  Second, an editor's mathematical tastes are reflected in a journal, and from this point of view we should compare the size of the editorial board to the size of the mathematical establishment.  Regarding the second point, observe that the editorial board does in fact grow somewhat as the size of the mathematical establishment grows.  

Let us try to analyze these choices through what one might call the submission life-cycle of a paper.  A paper is submitted and received by a editor.  The editor looks at the paper and might make an initial value judgment: is the paper ``good'' or ``interesting'' by whatever metrics the editor uses?  Of course, reasonable people can disagree at this stage, especially when making a quick judgment, as we will see in a case involving Lefschetz below.  Perhaps the editor is already familiar with the paper from a talk given by the author (or authors) or has already had detailed correspondence about the topic of the paper, or perhaps not.  Perhaps, as we have already seen can happen, the editor has solicited the paper.  

If the editor deems the paper of sufficient value, then they can send it to be refereed.  We discussed the historical standards of peer-review earlier, but checking correctness is something that has only evolved recently.  Instead, even if a paper was refereed, the editor was likely only seeking an opinion as to whether the paper is ``suitable''.  As we have seen above, the shape of a referee report can also itself be influenced by social pressures.  How does the editor weigh the opinion of the report: is it written by a high-status individual or someone young and untested?  Is it produced pro-forma only to be dismissed?  

It is important to remember that the status of a journal reflects {\em reciprocally} the status both of the editors and the authors whose papers appear in the journal.  The composition of the editorial board has a huge impact on the selection of papers that appear in the journal.  In the case of Annals of Math, the editorial board composition moves toward strong representation by ``abstract'' algebraists/topologists of the modernist stripe: many of the names mentioned by the late 1940s and 1950s are precisely the people whose papers we have been discussing.  
\end{entry}

With this in mind, let us explore Lefschetz's influence more directly in two episodes.  First, Parshall frames Lefschetz's early interventions at Annals in terms of the style of topology that was prioritized.  Alma Steingart, in her book ``Axiomatics'' \cite{Steingart} makes a number of remarks about Lefschetz's personality based on correspondence between N. Steenrod and R.L. Wilder.  Steenrod reported to Wilder ``He discussed the matter of young mathematicians acquiring the habit of publishing numerous papers on trivial problems.  It appears that the true Princetonian method is to work only on {\em general} problems and to publish only when some step in the theory has been accomplished'' \cite[p. 35]{Steingart}.  Lefschetz also ridiculed the Polish and R.L. Moore school of point-set topology, both of which had axiomatic leanings.  The Moore school was, in particular, lambasted as the ``concerning school'' \cite[p. 36]{Steingart} because so many papers appeared to be titled ``Concerning $X$'' for some topic $X$.

Gordon Whyburn was a student of R.L. Moore and Whyburn's student Beatrice Aitchison published an abstract of a forthcoming paper ``Concerning regular accessibility'' in the Bulletin of the AMS \cite{AitchisonBull}.  A more complete treatment of this work was submitted to the Annals of Math.  Having just presented a characterization of Lefschetz's view of the Moore school, the outcome was predictable: Aitchison's paper was rejected.  Parshall argues \cite[p. 242]{Parshall} that Whyburn did not take this rejection lightly: he viewed is a tacit rejection of his own mathematics by Annals.  Einar Hille delivered the news of rejection to Whyburn, framing it in terms of lack of sufficient enthusiasm and constraints imposed by the fiscal situation of the journal.  Whyburn was not content with this rejection responding: why not send the paper ``to someone thoroughly familiar with and active in {\em our field} such as Moore, Kline, Wilder, Ayres and ask them to referee it before passing final judgment.'' Evidently Whyburn did not feel the paper was given fair trial, but his third-party appeal to reverse the political rejection was, unsurprisingly ignored.  Aitchison's paper never appeared in Annals, and instead only appeared in 1933 in the journal of the Polish school Fundamenta Mathematica \cite{Aitchison}.  Even though we have described the Polish school of mathematics as a modernist school (see the discussion near the beginning of Section~\ref{ss:homologicalalgebra}), we see that modernist mathematics is not a monolith.  Whyburn did eventually become an editor at Annals of Mathematics in the 1940s as a representative of the AMS, but by this point combinatorial topology dominated the mathematical landscape and ``general'' topology as a pursuit in itself had receded somewhat (or its development had been absorbed within other subfields of mathematics).

To set the stage for our second episode, we describe comments by W.V.D. Hodge on Lefschetz.  Hodge writes:
\begin{quote}
	But possibly his main contribution to mathematics during the thirties lay in his powerful influence on others: he worked very hard to keep himself informed on what his students and associates were doing, and was a vigorous critic of anything he did not approve of. He asserted (with much truth) that ``he made up his mind in a flash, and then found his reasons''. Naturally he made mistakes this way, but once he was really convinced that he was wrong he could be extremely generous. 
\end{quote}
As regards the consequences of his actions, Hodge writes of Lefschetz:
\begin{quote}
	In this he was so successful that it is no surprise to find so many of his pupils in leading positions in the mathematical world.
\end{quote}
Putting these statements side-by-side, it hard not to draw the conclusion that Lefschetz consciously remade the landscape of mathematics in the way he wanted.  His student Albert Tucker was chair at Princeton for almost 20 years.  Other students included Ralph Fox, Paul Smith, Felix Browder, Norman Steenrod, John Tukey, Shaun Wylie and also people like C. Ehresmann who was a joint student with E. Cartan.  Many of these people also served on the editorial board of Annals of Math.  But there are only a few positions available.  What happened to the people that Lefschetz did not support?

Hodge reminisced about his interactions with Lefschetz \cite[p. 30]{HodgeLefschetz} around a theorem he proved (likely \cite{Hodge}) using ideas of Lefschetz from \cite{LefschetzCorr}. 
\begin{quote}
	I wrote a short paper....When this appeared Lefschetz ``made up his mind in a flash''-I was wrong and I must withdraw the paper. A correspondence lasting for several months took place, and during this period Lefschetz was travelling round Europe visiting mathematicians in various countries to whom he voiced his criticisms of my paper. Eventually we reached a state of armed neutrality, and he extended an invitation to me to visit Princeton. When I arrived there I was immediately instructed to conduct a seminar on my paper (owing to his characteristic the seminar actually lasted for six sessions), at the end of which he stood up and publicly retracted all his criticisms, and then, despite his handicap, wrote to all his European correspondents admitting that I was right and he was wrong. 
\end{quote}
At the time when this paper was written, Hodge was a young mathematician: he gained his MA degree at Cambridge in 1930.  However, the competitive aspects of mathematics at Cambridge have been discussed in \cite{Warwick}, and likely Hodge dealt well with Lefschetz's insistence.  In the 1930s, algebraic geometry was still in thrall to the Italian school, while Lefschetz's arguments were famously impenetrable; one wonders how the actual exchange went.

Undoubtedly such behavior and the concomitant editorial decisions have consequences.  One sample consequence: Mac Lane reminisces in his autobiography about {\em his} effect on a mathematician \cite[p. 155]{MacLaneAuto}
\begin{quote}
	He was remarkably cheerful as he reported that the critical reviews of his papers I had written in the Journal of Symbolic Logic had been used by French officials to block his promotion. This left me at a loss. 
\end{quote}
It is telling that the mathematician remains unnamed.

Hodge continues with discussion of Lefschetz's influence at Annals.
\begin{quote}
	He employed equally drastic methods in his capacity as editor of the Annals of Mathematics over a period of 25 years. No leniency was shown towards any paper submitted to the journal which was not up to his standards, and anyone who disagreed with his judgment had to work very hard to make him change his mind. But once he had decided that man was worth helping there was no limit to the aid he would give him. By these methods he made the Annals one of the top mathematical journals in the world, and he and his colleagues made Princeton a world centre of mathematics. In the course of this vigorous programme he made very few enemies indeed: one felt that there was no personal animosity in his bark, and no self-seeking: he just wanted to serve mathematics as best he could. 
\end{quote}

The treatment of papers whose authors Lefschetz did support could be rather different. Tucker recounts in \cite{TuckerInterviews} the case of a note that Tucker hoped would appear in PNAS:
\begin{quote}
	I sent this paper off to Lefschetz asking him to submit it to the Proceedings of the National Academy of Sciences. I assumed that that had been done, but when I returned to Princeton in September 1941 from the year’s leave of absence I discovered that Lefschetz had not submitted the paper to the Proceedings of the National Academy of Sciences because he was having a fight at that time with the editor of the Proceedings. Instead, he had submitted the paper to the Annals of Mathematics of which he was editor. I was very upset because the paper did not have complete details in it. It was merely an outline, a projection of what I intended to do, and it seemed to me that that was not an appropriate paper to be published in the Annals of Mathematics. So I withdrew the paper, and the paper has never been published.
\end{quote}

Rota recalls \cite[p. 19]{Rota}:
\begin{quote}
	He liked to repeat, as an example of mathematical pedantry, the story of one of E. H. Moore's visits to Princeton, when Moore started a lecture by saying, ``Let $a$ be a point and let $b$ be a point.'' ``But why don't you just say, 'Let $a$ and $b$ be points!''' asked Lefschetz. ``Because $a$ may equal $b$,'' answered Moore. Lefschetz got up and left the lecture room.
\end{quote} 
Recalling Lefschetz's view of ``general problems'' while simultaneously eschewing ``pedantry'' seems a striking distinction from a modern point of view.

These opinions give some idea of how Lefschetz defined the ``mathematical mainstream'' for himself in the 1930s.  Gone were the papers on classification of groups of a particular order or type that had proliferated during the Wedderburn era, or the papers on ``point-set'' style topology.  The emphasis was fully directed to ``combinatorial topology''.  Gone was the space for ``expository'' articles exemplified by Wedderburn's exposition of the theory of fields \cite{Wedderburnfields}, replaced only by pure ``research articles.''

From this point of view, one could imagine that Eilenberg was initially dismissive of Eilenberg's complaints about his singular cohomology theory, but eventually Lefschetz became a great supporter of Eilenberg (and as we see, eventually Eilenberg was an editor of Annals of Math).\footnote{In a letter to R. Brauer from Feb 13, 1958, Weil writes \cite[p. 13]{WeilFacultyFile} admiringly of Lefschetz's editorial work, giving at the same time his impression of the American Journal of Math:
	\begin{quote}
		Of the existing journals in the U.S.A., the American Journal of Mathematics has the longest and most impressive tradition. For many years, largely because of Lefschetz's superb editorship, the Annals of Mathematics have held the first place, while the American Journal was gradually sinking into mediocrity.
	\end{quote}
	He continues admiringly of Transactions of the AMS under Tamarkin's editorial guidance.  He closes with the statement ``For a number of reasons, scientific and sociological, it is unavoidable that a good deal of comparatively inferior material should get into print.''  To me, this reads as a reference to the fact that in the creation of status hierarchies, the presence of low-status individuals is necessary to enforce/support the high-status individuals. \cite{Ridgeway}.}

While publication is essential to mathematical production, it is worth analyzing the question: exactly what is the task of a ``good journal''?  Certainly it is supposed to make judgments about ``quality'' of mathematics.  One point of view was spelled out by Weil in \cite[p. 13]{WeilFacultyFile}:
\begin{quote}
	It is not necessary, but it seems clearly desirable and convenient, that the greater part of the best mathematical work done in each country should be concentrated in one or more first-class journals rather than diluted with inferior material.
\end{quote}
It is reflective of the volume of mathematical publication of the day that such an opinion was even plausible, but Weil's use of the word ``best'' makes it sound as if the ``best'' mathematics is self-evident.  However, he makes a number of supporting comments that give us some sense of how one determines ``best'' mathematics: anyone who wants to judge good mathematics must be ``broadly mathematically cultured'' and one shouldn't try to involve too many people because ``consensus'' editing leads to a weakest link problem.\footnote{He writes: 
	\begin{quote}
		When several editors have to work together, it is almost inevitable that the weakest one will drag his colleagues down to the level of his own standards. But no one can maintain high standards unless he has had considerable experience in more than one branch of mathematics, and, above all, unless he enjoys a secure and unassailable position in the mathematical world.
\end{quote}} While ``universal'' mathematicians, in the vein of Hilbert were already a thing of the past due to increasing production and specialization, Weil's view is again reflective of the scope of the mathematical establishment; it seems even harder to square with the concomitant size of the establishment today.

Weil in \cite[p. 14]{WeilFacultyFile} explains his principle as follows:
\begin{quote}
	rejecting many papers which, by the usual standards, were perfectly publishable. We took the view that, by doing so, we did little harm to the authors, provided it was done quickly; for such papers would get into print anyway. This applied particularly to the case of well-established authors; we thought that promising young beginners should receive a more considerate treatment and should at times
	be measured by less exacting standards. We could not hope to be exempt from errors; whether or not quick decisions imply a higher percentage of errors than a slower procedure is a debatable point. Had we been dictators over the whole field of mathematical publications, the weight of our responsibility would have been crushing; since, however, opportunities for publication of even moderately deserving papers remained plentiful even after rejection from the Journal, we
	found that we could carry it lightly.
	
	Unrealistically high standards would have stifled the Journal altogether; this had to be avoided. On the other hand, our policy made it possible to reduce to a bare minimum the backlog which is the curse of most journals. We found that we had no difficulty in publishing quickly whatever we thought worth publishing at all.
\end{quote}
Implicitly here, Weil and his coeditors get to decide who the ``promising young beginners'' and ``well-established authors'' are.  The community was small, so perhaps there was even something approaching a consensus, but once again, ``promising'' is a mark of status.

Lefschetz was a strong figure, and students like Tucker followed in his footsteps establishing Princeton as a center for topology in the 1940s and 1950s, buoyed by the presence of luminaries at the nearby IAS. Cartan's students were undoubtedly successful, but due to what personal efforts of Cartan.  There is a reciprocal mechanism acting between ``exciting mainstream mathematics'' and ``faculty hiring.''  Status gets dynamically transferred between these two systems and also journals. 
\aravind{Lots to expand here!}

\subsubsection*{Bourbaki and category theory}
By the end of the 1940s, category theory was in the air, but the question of its reception seems to me more complicated.  At this point one sees a certain fracture between American and European directions.  Whether people {\em used} categories themselves seems to be a question of individual taste, and to work on category theory was another matter altogether.  

The first crop of Eilenberg's students (he started taking students in the late 1940s and early 1950s), a group including Alex Heller, David Buchsbaum and Dan Kan were undoubtedly familiar with category theoretic notions and all of these people wrote papers that were directly category-theoretic related, though Kan's story is perhaps a bit more complicated than the others.\endnote{From \cite[p. 1043]{KanNotices}: ``In the
	Spring of 1954, Samuel Eilenberg came from Columbia University on
	a visit to the Hebrew University in Jerusalem. (At the time, Eilenberg
	was already one of the leading figures in algebraic toplogy, due to his
	work with Saunders Mac Lane and the influential Foundations of Algebraic Topology which he had just written with Norman Steenrod.)
	Kan knocked Eilenberg’s hotel room door, and explained his simplicial
	description of homotopy groups. Eilenberg asked him if he could prove
	the homotopy addition theorem, and Kan returned a week later with
	a proof. Eilenberg told Kan that he had a thesis there, engineered
	an ad hoc arrangement giving Kan the status of graduate student at
	the Hebrew University, and in the summer of 1954 Kan submitted his
	thesis. He formally received his PhD in 1955.''}  Buchsbaum's thesis was devoted to ``exact categories'', a predecessor to the later notion of abelian category \cite{Buchsbaumthesis}.  In the MathSciNet review of \cite{Buchsbaumexact} by H. Cartan, which expands on what appeared in Cartan--Eilenberg, we learn that Grothendieck had formulated a notion of abelian category by this point.  
	
	The situation involving Bourbaki and category theory has been treated elsewhere in much greater detail.  It is treated in great detail in \cite[Chapter 8]{Corry}, and also in great detail in Kr\"omer's paper \cite{KromerBourb}.  We mention just one episode that has some bearing on the later discussion.  In February 1955, Grothendieck discussed formulating sheaf cohomology in terms of categories of modules in \cite[Feb 26, 1955]{GrothendieckSerre}.  Serre responds \cite[Mar 12, 1955]{GrothendieckSerre} that Cartan was aware that sheaf cohomology of some flavor could be formulated in terms of the Cartan--Eilenberg notion of derived functor, and had asked Buchsbaum to develop this, but was unaware of what happened.  This correspondence makes it clear that Serre was familiar with and conversant with the theory of categories, but he scrupulously avoided using this terminology in his written work.

\subsection{The local to global paradigm in sheaf theory}
\label{ss:sheaves}
\addtoendnotes{\vskip .2em\noindent {\bf The local to global..} \vskip .2em}
The fourth cultural touchstone I'd like to discuss is Leray's theory of sheaves, which was invented by Leray around 1946 after his stay in a prisoner's camp during World War II.  This stay has been recounted many places, e.g., \cite{MillerLeray} and there are a number of sources describing aspects of the history of the theory of sheaves, e.g., \cite{HouzelSheaves} and \cite{GraySheaves}; I will draw on these works, but I'd like to start a bit earlier.  While the algebraic-topological influences on the theory of sheaves are very well-represented in the literature, Cartan's reworking of Leray's ideas will be very important for our narrative; Chorlay \cite{Chorlay} discusses the history of sheaf theory from this point of view,\footnote{Chorlay also has a nice discussion of the construction in the subtitle: ``local-to-global'' in \cite{ChorlayLocalGlobal}.} but I'd like to bring the discussion back through the lens of Cartan--Weil, viewing sheaves as a complementary perspective to that brought by the theory of the fiber bundles.  
  
Leray was a student at the ENS, overlapping with Weil, Cartan, Dieudonné and Chevalley, and received his degree in 1933.  Leray participated in some of the initial Bourbaki meetings, but Leray moved away from Bourbaki \cite{AndlerLeray}; the complicated relationship between Leray and Bourbaki has recently been reanalyzed \cite{Ricotier}, and plays a not inconsiderable role in our discussion.  A clear disjunction between Leray's mathematical tastes and those of other Bourbaki members is an evident source of tension, but the ways in which these tensions changed and eventually bubbled over during subsequent years is important.  We claim that there is some reflection of these tensions in the reception of sheaf theory.  Indeed, the theory of sheaves as initially developed by Leray was not immediately absorbed by the postwar French mathematical mainstream.  Undoubtedly Leray's status played a role in the willingness of others to grapple with his ideas in the first place.  

Leray announced his theory of sheaves in a collection of announcements published in the {\em Comptes rendus de l'Académie des Sciences} (CRAS), appearing as: \cite{Leray1,Leray2,Leray3, Leray4}.  After discussing some aspects of Leray's biography and early mathematical work that seem relevant to his presentation of ideas, I aim to show that there were two competing evaluations of Leray's initial CRAS notes.  On the one hand, there is a collection of individuals, one might call it the Eilenberg--Weil axis, who questioned the novelty of Leray's ideas.  On the other hand, there is an opposing Cartan--Hopf axis that was convinced that Leray's theory had much to offer.  I aim to unearth some sources for this tension and simultaneously to explain that Cartan's support, in particular his focused reworking of Leray's theory and dissemination by way of his seminars, cleared the way to wider applicability of these ideas.  

The change in perspective offered by Cartan on the nature of sheaves was {\em integral} to the mainstream acceptance of sheaf theory, especially outside France.  In the context of the broader narrative, I aim to show how Cartan's broadening of the sheaf concept, wresting the idea from Leray in a sense, is essential in bringing the theory of fiber bundles to bear on algebraic geometry.  In Section~\ref{ss:fiberbundlesandprojectivemodules}, I will return to this point and describe its influence on the reception of the projective module concept.

\subsubsection*{Leray--Schauder theory becomes applied topology?}
Leray's thesis related to applied analysis, specifically mathematical problems related to fluid flow; certainly looking back with modern eyes, this seems quite far from the topics discussed above, but let us not refer to Leray as an applied mathematician in an unqualified way.  In the Spring of 1933, Hans Lewy introduced Leray to Julius Schauder in a restaurant on the rue Soufflot \cite{LerayonSchauder}.  Leray continues:
\begin{quote}
	I immediately said to Julius Schauder: ``I have read your paper on the relationship between existence and uniqueness of solutions of a nonlinear equation.  I know now that existence is independent of uniqueness.  I admire your topological methods.  In my opinion they ought to be useful for establishing an existence theorem independent of the whole question of uniqueness and assuming only some a priori estimate.''
\end{quote}
In practice, the tools used were analogs of fixed-point like theorems of Brouwer in the setting of (infinite dimensional) function spaces.  Leray and Schauder quickly applied their techniques to establish existence result for certain non-linear partial differential equations; Leray refers to applications of the Dirichlet problem for two-dimensional convex regions as key. Of course, the Dirichlet problem itself has a long pedigree, going back to Riemann's ``problematic'' variational approach that was put on sound footing by Hilbert.  Evidently the Leray--Schauder method was existential and {\em non-constructive}, a point to which I will return momentarily, but undoubtedly these results established ``first contact'' between Leray and the tools of algebraic topology.
 
Leray's life has been described as curiously parallel to Weil, but in some ways opposite \cite{EkelandLeray}.  Both were born in 1906 and died in 1998, both were Normaliens.  Weil avoided WWII escaping to the United States, while Leray fought until he was captured.  Ekeland argues Leray is the ``applied'' yin to Weil's ``pure'' yang, writing \cite{EkelandLeray}:
\begin{quote}
	Leray viewed mathematics as a tool for modelling, and drew his inspiration from problems in mechanics and physics, such as fluid dynamics and wave propagation.  He was fond of explaining how the road from mathematics to applications is two way, and how a purely mathematical theorem (concerning, for instance, the existence and uniqueness of solutions of systems of partial differential equations) might have profound physical implications.
\end{quote}
This is reflected to some extent in Leray's own self-identification.

One gets a sense of Leray's tastes by looking at his commentary on some of his correspondence with Schauder from the period.  Schauder attended the Moscow Topology colloquium of 1935 and presented a talk ``Some applications of the topology of functional spaces'' \cite[p. 38]{ANSMoscow} touching on his joint work with Leray.  Schauder sent Leray a letter discussing his views on the conference.  Leray views Schauder's pronouncements as giving ``evidence of very sound judgment and great breadth of mind, analysing with penetration the capabilities and limits of everyone, whether they were admirably great or somewhat narrow'' \cite{LerayonSchauder}. He highlights Schauder's comment ``I am like you a man of applications.''\footnote{This correspondence was also discussed by Eckes \cite{EckesLeraySchauder} who adds some color to Leray's description by mentioning that Schauder's letter emphasizes the misunderstandings, in particular, of John von Neumann and Hans Freudenthal.  Eckes writes \cite{EckesLeraySchauder} \begin{quote}La lettre de Schauder s’achève sur un résumé de sa conférence, soulignant les mécompréhensions dont elle a fait l’objet de la part de John von Neumann et de Hans Freudenthal en particulier.\end{quote}}

For all this talk of a ``mindset of applications,'' one sees a slightly different side of Leray in the Bourbaki meeting of January 14, 1935 \cite{Delsarte}.  The extract we are about to give turns around the so-called Cauchy--Lipschitz theorem, an existence theorem for ordinary differential equations.\footnote{Weil trouve les hypothèses demandées, anti-naturelles ; Leray critique vivement la démonstration donnée par Goursat et attire l’attention sur l’intérêt qu’il y aurait à donner avant tous les théorèmes de ce type, théorèmes maintenant classiques, des théorèmes généraux, de caractère topologique, qui sont des théorèmes d’existence purs et non des théorèmes de calcul, mais qui permettent de prévoir quand il est possible d’énoncer un théorème de calcul.  Delsarte termine en critiquant la désinvolture avec laquelle chacun est arrivé à la présente séance. Il faut préparer ce qui a été demandé et ne pas se fier à ses facultés d’improvisation, si brillantes qu’on les estime.}
\begin{quote}
	Weil finds the stated hypotheses unnatural; Leray strongly criticizes Goursat's proof and draws attention to the interest that there would be in giving before all theorems of this type, now classical, general theorems, of a topological character, which are pure existence theorems and not theorems of calculus, but which make it possible to predict when it is possible to state a theorem of calculus. 
	
	Delsarte ends by criticizing the nonchalance with which everyone has arrived at the present session. It is necessary to prepare what has been requested and not to rely on one's faculties of improvisation, however brilliant they are considered to be.
\end{quote}
As Peter Lax writes in his foreword to Volume II of Leray's collected works \cite{LerayCollectedII}: ``physicists sometimes deride such existential pursuits by mathematicians.''  Certainly Leray's work seems applied by contrast to that of Weil.  Nevertheless, looking back with modern eyes, the idea that one would want to formulate ``pure existence theorems'' seems quite ``pure'' and makes the view that Leray was an ``applied mathematician'' seem over-simplified.

Leray's work in algebraic topology has been described as ``obscure'', and I'd like to analyze that characterization now.  To set the stage for this discussion, we need to review some of the development of Leray's point of view on algebraic topology, a key source of which is the Leray--Schauder theory (I will use this as a convenient summary for the joint work with Schauder above).  In the prisoner's camp, Leray had famously given a course on algebraic topology ``in captivity,'' eschewing his more ``applied'' interests for fear that they could be used by the German war machine.  Having limited interaction with development of algebraic topology at the time, save a few articles that had been procured for him by Hopf, Leray had developed his own perspectives, and consequently his own notational and terminological idiosyncracies.  

The results of this course were published in outline form as CRAS notes \cite{LerayCRAT1,LerayCRAT2,LerayCRAT3,LerayCRAT4} during the war. Eilenberg reviewed the CRAS notes for Math Reviews, and the conclusion sounds somewhat dismissive:\footnote{In the paper \cite{Cechbicompact}, Čech gives various equivalent characterizations of bicompact spaces, one of which is equivalent to what we now call a quasi-compact space.}
\begin{quote}
	The outline is very concise with no proofs and not very complete definitions but it seems likely that the cohomology theory developed will coincide with that of Čech (at least in the case of bicompact spaces). In the remaining two notes the author uses his homology theory to outline the topological principles that lead to existence theorems...
\end{quote}
Čech's approach to homology was published initially in 1933 \cite{CechI}, developing ideas of Alexandroff, and Leray's work postdates this by some time, so it seems worth exploring Eilenberg's allegation somewhat. 

\begin{entry}[Čech and homology]
	\label{par:cechhomology}
	In 1932 at the Zurich ICM, Čech delivered a lecture introducing higher homotopy groups.  He already knew that higher homotopy groups were abelian, and famously was ``persuaded by Hopf'' to withdraw his report; now only a paragraph appears \cite{CechICM}.\footnote{G.W. Whitehead claims that inklings of higher homotopy groups were also studied by Dehn \cite[p. 2]{Whitehead50years}.}  Alexandroff reportedly remarked: ``But my dear Čech, how can they be anything but the homology groups?'', ironic since the Hopf map $S^3 \to S^2$ had been defined by 1930 and was known to be homotopically non-trivial \cite{Whitehead50years}.
	
	 As mentioned earlier, a veritable zoo of approaches to homology theory developed in the wake of Poincar\'e's initial approach.  Čech's approach to homology developed Alexandroff's idea of combinatorialization of a space using, in modern terminology, limits of nerves of finite covers of spaces.  This approach to homology differed from other approaches in several ways.  For example, it differed from the approach of Vietoris which relied on a different combinatorialization that would perhaps nowadays be described as simplicial approximations to a space.   The variant approaches also differed in the class of spaces to which they were applicable.  On the one hand, this methodological variety led to differing perspectives on what ``key'' properties of homology should be.  From this point of view, the explicit methodological variations were later ``tamed'', or perhaps organized, by the Eilenberg--Steenrod axiomatic approach.  
	 
	 The Čech approach was distinguished by what are now called {\em continuity} properties: the resulting homology theory ``commuted'' with formation of limits of spaces.  Moreover, Čech's approach was understood to be well-suited to ``local'' analysis.  Indeed, by 1934 Čech had used his approach to homology to define ``local'' Betti numbers \cite{Cechlocal}.  The appearance of this paper in Annals of Mathematics can, in view of our discussion of publication practices in this journal, be viewed as tacit endorsement of this approach by Lefschetz by this time.  Lefschetz writes about Vietoris homology: ``It is noteworthy for its convenience in many applications, but as we shall show, it is in fact reducible to the Čech theory'' \cite[Chapter 7, Point 23]{Lefschetz}.  Thus, by the mid 1930s, the Čech approach to homology seems to be well-known. 
\end{entry} 

\subsubsection*{Leray after the war}
Leray's motivation for revisiting the foundations of homology were evidently grounded in his work with Schauder: he wanted an approach that was general enough to apply to Banach spaces.  In view of the vast variety of approaches to defining homology at the time, the search for yet another variant with slightly different domain of applicability seems unsurprising.  Eilenberg's assertion that Leray's theory ``coincides'' with the Čech theory thus seems, retrospectively, somewhat reductive.  

The fleshed out versions of the CRAS notes were submitted to Henri Villat at the Journal de Mathematiques Pures et Appliques in January of 1944.  Villat begins his preface to Leray's papers by explaining that Leray's treatment is the first ``didactic treatment'' since the Alexandroff--Hopf text.\footnote{Villat writes \cite{LerayCAT1}: \begin{quote} ``L’exposé rédigé par Mr J. Leray \guillemetleft Sur la forme des espaces topologiques  et sur les points fixes des représentations\guillemetright, est la première partie d’un Cours  de Topologie algébrique, appelé à faire quelque bruit dans le monde mathématique. Le sujet est neuf et de grande actualité. Mais en dehors d’un livre  de MM. Alexandroff et Hopf, il n’existe encore aucun traité didactique sur ces sortes de questions.'' \end{quote}}  He closes with:\footnote{Il est presque superflu d’insister sur l’opportunité d’une telle publication: on sait le renom de l’auteur, dont les nouvelles méthodes concernant les équations différentielles ou aux dérivées partielles, et les équations intégro—différentielles, ont en immédiatement un retentissement mondial. Ajoutons  que M. J. Leray, professeur à la Sorbonne, a écrit le présent travail en captivité (il est encore prisonnier, détenu à l’Oflag XVII). Le travail actuel a reçu de M. Hopf, professeur à l’Université de Zurich (savant d’une compétence notoire sur le sujet,) une adhésion enthousiaste...} 	 
\begin{quote}
	It is almost superfluous to insist on the timeliness of such a publication: we know the renown of the author, whose new methods concerning differential equations or partial derivatives, and integro-differential equations, immediately had worldwide impact. Let us add that Mr. J. Leray, professor at the Sorbonne, wrote the present work in captivity (he is still a prisoner, detained at Oflag XVII). The current work received enthusiastic support from Mr. Hopf, professor at the University of Zurich (a scholar of noted competence on the subject)...
\end{quote}
The papers appeared in sequence in 1945 as \cite{LerayCAT1,LerayCAT2,LerayCAT3}.  Eilenberg reviewed these more complete treatments as well writing:
\begin{quote}
	The starting remark is that the product of two abstract complexes is again an abstract complex, while the product of two simplicial complexes is not (without further construction) a simplicial complex. In order to take full advantage of this convenient property of abstract complexes the author by-passes the usual ``covering-nerve'' scheme and constructs a cohomology theory of topological spaces more directly connected with abstract complexes. The resulting cohomology theory seems to be isomorphic (at least for compact Hausdorff spaces) with the groups obtained following the Čech scheme.
\end{quote}
This review contains slightly more detail about the methods, but Eilenberg appears to confirm his earlier conclusion.  

Houzel characterizes Leray's approach \cite[\S 2]{HouzelSheaves} slightly differently writing: ``il cherchait \`a se d\'ebarrasser des hypoth\'eses inutiles et \`a associer aux espaces topologiques des invariants alg\'ebriques sans passer par des constructions interm\'ediaires''; this sounds to me more like a {\em reworking} of Čech theory, but this deserves some discussion.  Čech's initial investigation studied homology.  While the Alexandroff--Čech approach to homology is mentioned above, Čech also worked on product structures of Alexander--Kolmogorff.  Leray speaks of the ``homology ring'' referencing the work of Alexander--Kolmogoroff and the latter work of Čech that we now call the cohomology ring.  

Leray sought to avoid the complicated ``combinatorialization'' used by Čech involving the nerve of a space and the associated auxiliary constructions.  These papers of Leray never mention the word sheaf, speaking instead of abstract and concrete complexes.  Leray explicitly mentions work of Alexander \cite{AlexanderGratingI} extending the latter author's earlier work on cohomology rings (the connectivity ring in those papers) to more general settings using a notion called a ``grating''.  Gray suggests \cite{GraySheaves} that Leray's notion of concrete complex and subsequently ``couverture'' is a first step towards his invention of the sheaf notion: a concrete complex is an algebraic structure (an abstract chain complex) together with a rule or ``support function'' for attaching this algebraic structure to closed subsets of a suitable topological space.  

Leray mentions a good deal about de Rham's approach to cohomology citing what we now call the product in the exterior algebra.  We have already mentioned Weil's first meeting with Eilenberg where he brings up de Rham theory as a putative challenge to Eilenberg's axiomatic approach to homology.  Miller characterizes Leray's approach using ``concrete complexes'' and couvertures as an axiomatization of the de Rham theory of forms.  Specifically, the axiomatics of couvertures were built in order to encode abstractly the Poincar\'e lemma, that, viewed sufficiently ``locally'', closed forms are exact. 

Yet another source of this attention to ``locality'' in Leray's work stems from the Leray--Schauder theory and eventual applications to fixed point theorems: in order to compute the degree, one would like situations where it can be computed in terms of suitable local contributions \cite[p. 284]{Mawhin}.

In June or July 1945, just after Leray had returned from the prisoner's camp, Weil remarks on meeting him in Paris that:\footnote{\cite[p. 526-7]{WeilII}	``...il m'avait parl\'e de sa ``cohomologie a coefficients variables''; c'est ainsi du moins que le souvenir m'en \'etait rest\'e, et cette id\'ee, bien que vague dans mon esprit, m'avait frapp\'e.''} 
\begin{quote}
	...he spoke to me about his ``cohomology with variable coefficients''; at least that's the memory that remains with me, and the idea, although vague in my mind, struck me.
\end{quote}
 
In 1946, Leray published two further CRAS announcements where the notion of a sheaf finally does appear \cite{Leray1,Leray2}.  Miller characterizes the reception of these announcements as follows \cite[\S 3]{MillerLeray}:
\begin{quote}
	This work must have seemed incredibly obscure at the time—a sheaf, a new concept certainly containing a lot of information, was immediately fed into an unstudied cohomology theory, and then used to define an invariant of a map on which there appeared without any motivation a highly complex structure, all expressed in a terse and contrarian style.
\end{quote} 
There is much to say about the reception of these works, which I don't believe Miller's description completely captures. While they were certainly terse (could a CRAS announcement be anything but compressed?) much of our discussion will turn on questions of novelty of the new notions.

Eilenberg was the reviewer of these CRAS notes for Math. Reviews again.  His review can be viewed as providing the word ``bundle'' as an English translation of ``faisceau'', and he states that these bundles of groups can be viewed as coefficient system for a homology theory. Unlike previous works that were interested in cohomology of a space, Leray was interested in what he called the ``homology ring of a representation.''  Leray was explicitly interested in Steenrod's approach to ``local coefficients''.  Furthermore, he aimed to associate cohomological information to a continuous map, rather than to a topological space.  Houzel goes further \cite{HouzelSheaves}, suggesting, by way of the fact that Leray references work of Steenrod on the homology of fiber bundles, that he had in mind computations of homology of fiber spaces of various sorts. This interpretation is indeed supported by some of Leray's later work it's important to remember that Picard and Lefschetz had analyzed various classes of fibrations establishing results that know go under the name ``Picard--Lefschetz theory'', and Leray was undoubtedly aware of these results.

Eilenberg was undoubtedly familiar with Steenrod's theory of local coefficients.  His review mentions in particular, that given a continuous map of topological spaces, the cohomology of the fibers with coefficients in the bundle of groups again yields a bundle of groups, and closes with comment that ``the second article enters in more detail into the structure of this new group and states without proofs a number of applications.''  This description demonstrates awareness that the situation Leray considers is more general than fiber spaces in topology, though that would eventually be an application. In retrospect, this is the germ of the Leray spectral sequence, and Eilenberg mentions no consequences of Leray's analysis.


\subsubsection*{Sheaves filtered through Cartan and Weil}
Before the war broke out, there was a certain amount of jockeying for mathematical positions.  Tensions between Weil and Leray were heightened during these various campaigns.  One clear view into these tensions stems from a letter from Weil to Jean Coulomb in 1938 around the promotion of Bourbaki member Mandelbrojt that has been discussed in detail in \cite{Ricotier}.\footnote{Jean Leray s’est retourné vers Émile Picard et Gaston Julia, a fait une campagne acharnée (il a eu le toupet d’écrire à [Henri] Cartan qu’il avait passé 15 jours à Paris à « étudier l’état d’esprit de l’Institut »), et, en faisant valoir (il a eu l’impudence de me le déclarer) l’argument xénophobe et aussi, je crois, l’argument antisémite, il a réussi à récolter 21 voix contre 26 à M[andelbrojt]. Voilà donc Julia et Villat, qui l’ont porté aux nues à l’Institut, et qui ne peuvent plus s’en dédire à la Sorbonne}
\begin{quote}
Jean Leray aligned with Émile Picard and Gaston Julia, campaigned fiercely (he had the nerve to write to [Henri] Cartan that he had spent 15 days in Paris “studying the spirit of the Institute”), and, by putting forward (he had the impudence to tell me) the xenophobic argument and also, I believe, the anti-Semitic argument, he managed to garner 21 votes against 26 for M[andelbrojt]. So here are Julia and Villat, who praised him to the skies at the Institute, and who can no longer go back on their words at the Sorbonne.
\end{quote}
The view of increasing tension between Leray and Bourbaki around this affair goes further.  Julia ran a seminar around the time, and various Bourbaki members had contributed talks.  Julia unceremoniously decides to remove all Bourbaki members from his seminar, and Leray's role in this removal is described by Weil \cite[4.3]{Ricotier} thus: \footnote{[Leray] s’est arrangé pour faire comprendre à Julia qu’il avait eu grand tort de laisser à Bourbaki la conduite de son séminaire, que le séminaire marchait mal, qu’il n’en sortait pas de Travaux, et qu’avec sa collaboration à lui Leray tout marcherait bien autrement. Ce qui explique qu’à la fin des deux derniers exposés de l’année (Chevalley et Pisot, exposés excellents d’ailleurs) Julia a engueulé les conférenciers avec un parti pris manifeste, et a proprement limogé notre équipe, en annonçant que l’an prochain le séminaire aurait lieu « dans des conditions différentes ».}
\begin{quote}
	[Leray] managed to make Julia understand that he had been very wrong to let Bourbaki lead his seminar, that the seminar was going badly, that no Travaux were coming out of it, and that with his collaboration, Leray, everything would work out quite differently. Which explains why at the end of the last two presentations of the year (Chevalley and Pisot, excellent presentations by the way) Julia yelled at the speakers with a clear bias, and properly dismissed our team, announcing that next year the seminar would take place ``under different conditions''.
\end{quote}
These accusations are strong, but nevertheless we saw above that in 1945 Weil was still talking to Leray, and evidently interested in his mathematical opinions.  

In 1945, Weil moved to Sao Paolo, Brazil and had assumed a temporary role as an instructor.  This period takes up but a few pages in his autobiography--effectively an epilogue.  About his mathematical interactions he writes \cite[p. 192]{WeilAuto} the circumspect sounding:
\begin{quote}
	Besides, despite Zariski's presence in Sao Paulo in 1945, and Dieudonn\'e's in 1946 and 1947, it was impossible not to wish for a more stimulating scientific milieu. My Parisian friends thought it possible to arrange for my appointment to the College de France when Lebesgue's retirement left a chair vacant, but this plan did not materialize. 
\end{quote} 
The Cartan--Weil correspondence from this period paints a rather different picture and involves Leray on several fronts.  The appointment at the College de France to which Weil refers, in 1946, was eventually assumed by Leray.\footnote{Leray's published recollections of this, which I take from \cite[\S 5]{MeyerLeray}, are rather bland: \begin{quote}C’est [ce travail de topologie alg\'ebrique] qui m’a fait entrer \`a mon retour de captivit\'e au Coll\`ege de France...Il y eut un drame \`a cause de l’attitude d’Andr\'e Weil pendant la guerre: elle ne fut pas admise par cette maison, qui a hautement le sens du devoir national...\end{quote} Much more is written about Weil's view of this episode, and we will discuss it momentarily.}  The correspondence itself refers to this episode as ``L'affair du college''.  One cannot help but keep all these scuffles in mind as one confronts the tone of the letters exchanged between Cartan and Weil and, consequently, the mathematics they discuss. 

Weil's comment about wishing for ``a more stimulating scientific milieu'' can also be unpacked a bit.  Certainly the institute in Sao Paulo had a library, but which periodicals did they receive?  How did Weil get to hear about current mathematics?  On the one hand, the Cartan--Weil correspondence itself answers this question: there is rather a lot of it saved from this period, with Weil's epistolary efforts presenting something of a bridge to the outside mathematical world, supplemented by the visits he mentioned above.  Indeed, it appears Weil receives a cache of interesting articles from visiting friends, attached in letters he receives and through French periodicals delivered to the embassy.\endnote{Let us take a moment to understand that Weil is processing all of this information from Brazil, with apparently limited access to mathematical journals.  One has to imagine that he is using the resources to which he has access, Cartan's letters to effectively recreate these ideas.  The frustrations of this method are probably clear and summarized in Weil's letter from the 14th of March \cite[p. 204]{CartanWeil}\footnote{Il n'y a aucun espoir pour moi de voir les memoires de Leray avant mon retour en France, a moins que tu ne me les envoies (sur la question des periodiques francais, comme sur tout le reste, la pauvre Mme Mineur, sur laquelle nous avions fonde de si grands-espoirs, s'est montree d'une nullite total).  Pour le Grand Oeuvre, je crois que ca ne vaut guere la peine et les frais d'un envoi (il sera bien temps en octobre), mais tu pourrais bien acheter les numeros des C.R. ou tu dis qu'ont paru ses notes sur les coefficients variables, et m'envoyer les notes in question dans ta prochaine lettre.}
	\begin{quote}
		There is no hope for me of seeing Leray's memoirs before my return to France, unless you send them to me (on the question of French periodicals, as with everything else, poor Madame Mineur, on whom we had pinned such high hope, turned out to be useless).  For the Great Work, I believe that it is hardly worth the trouble and expense of sending it (it will be good time in October), but you could well buy the issues of the C.R. where you say that his notes on variable coefficients appear, and send me the notes in question in your next letter.
\end{quote}
Gabrielle Mineur had been cultural attaché at the French Embassy in Brazil since January 1945, and had also worked at the CNRS \cite[p. 542]{CartanWeil}.  The correspondence also contains an extract from a letter from Weil to Yves Rocard \cite[p. 55]{CartanWeil} that is even more direct about Mme. Mineur:\footnote{Je ne comprends rien a ce qui lui arrive depuis qu'elle est au Bresil: elle autrefois si active, et si intelligemment active, ne repond plus aux lettres de qui que ce soit, a echoue dans tout ce que nous attendions d'elle, et en est au point de ne meme plus nous avertir de sa presence a San Paulo quand elle y vient.  Je presume que, d'une part l'atmosphere de pourriture de l'ambassade aura deteint sur elle et que d'autre part l'ete de Rio lui aura ote tous ses pouvoirs de reaction.}\begin{quote}
	I don't understand anything that has happened to her since she was in Brazil: she was once so active, and so intelligently active, no longer responds to anyone's letters, has failed in everything we expected of her, and is to the point of no longer even warning us of her presence in San Paulo when she comes there.  I presume that, on the one hand, the rotten atmosphere of the embassy will have rubbed off on her and that on the other hand the summer of Rio will have taken away all her powers of reaction.
\end{quote}}
Thus, the ``scientific milieu'' in which Weil found himself sounds rather spartan, and given the context in which Weil left France, it strikes one almost as mathematical exile.\footnote{He expresses a sentiment in this direction in a letter to Cartan \cite[p. 93]{CartanWeil} from November 1944: \begin{quote}...pour moi, je n'ai pas tr\'es brillamment r\'eussi; en ce pays o\'u d'ailleurs je n'ai jamais eu l'intention de me fixer.  Non sans h\'esitation, j'ai accept\'e une chaire au Br\'esil, ou je dois rester deux ans; mais je compte venir passer en France les vacances de 1945-46, d'octobre \`a mars.\end{quote}} 

Against this backdrop, in January 1947, Weil wrote a letter to Cartan where he sketches a proof of the de Rham theorem using local--to--global ideas: he asserts the existence of suitable ``simple'' coverings; on each covering patch, one can appeal to the Poincar\'e lemma, and then appeal to combinatorics to conclude.  While Weil works with closed coverings, some variant of this is one of the standard ``modern'' proofs of the result.

Weil's letter to Cartan is reproduced in \cite[pp. 45-47]{WeilII}, but a different version appears in \cite[p. 139--142]{CartanWeil}.  The version in Weil's collected works makes no mention of Leray, leading Miller \cite{MillerLeray} to write: ``It is curious and surprising that it was Andr\'e Weil rather than Leray himself who found the modern proof of the de Rham theorem, since this proof was a vindication of the local methods espoused by Leray. While Weil does not acknowledge Leray's influence...'' However, the version of this letter in \cite[p. 142]{CartanWeil} includes two further paragraphs, which do mention Leray \footnote{Weil writes: \begin{quote}
		En particulier, si les groupes de cohomologie, \`a domaines de coefficients variables, de l'illustre Leray, sont bien ce que j'ai cru comprendre quand il m'en a dit un mot \`a Paris, ils doivent s'encadrer tout naturellement dans ton expose\'e de la cohomologie, et dans la principe de la me\'ethode ci-dessus; L. a-t-il commence\'e \`a publier sur ce [barr\'e:sujet?] sujet?  Je n'ai pas encore vu son grand opus de Liouville, mais d'apr\'es l'analyse d'Eilenberg dans Math. R. il ne s'y trouve rien sur ce sujet, et m\^eme il ne s'y trouve rien de bien nouveau.  Son invention des coefficients variables serait-elle orien\'ee vers une th\'eorie genre Morse?  Il conviendrait de savoir le plus t\^ot possible s'il y a quelque chose \`a en tirer pour nous.
\end{quote}  Recall that the {\em Journal de Mathématiques Pures et Appliquées} was founded by Liouville.}
\begin{quote}
	In particular, if the cohomology groups with local coefficients of the illustrious Leray, are indeed what I thought I understood when he said a word to me about them in Paris, they should fit quite naturally into your exposition of cohomology, and into the principle of the method above; has L. begun to publish on this subject? I have not yet seen his great opus in Liousville's journal, but according to Eilenberg's analysis in Math. R. there is not much new to be found in this treatment. Would his invention of variable coefficients be oriented towards a Morse-like theory?  It would be nice to know as soon as possible if there is something to be gained from it for us.
\end{quote}
This translation does not suggest that Weil viewed Leray's comments to him as ``inspiration''.  I feel like this reading is supported by the last line of the quote above that seems to orient Weil, and by association Cartan, in opposition to Leray.  It seems natural to read Weil's mention of the ``illustrious'' Leray as mocking (indeed, he uses this construction numerous times in letters to Cartan) or, perhaps, even jealous.\footnote{We have already discussed some of the tensions between Weil and Leray above, but one can add to this a further rupture caused by the war. In the spirit of Ekeland's ``opposing forces'' we can speak of Leray the ``illustrious'' war hero vs. Weil the shamed\aravind{Wrong word?}.  In support of this point of view, in addition to Leray's comments above, I mention here the discussion in \cite[p. 553]{CartanWeil} of ``L'affair College'' which describes Weil's aborted candidacy for a position at the College de France, contra Leray and stymieing Weil's {\em redemptive} return to France.}  Weil seems to equivocate: he appears simultaneously skeptical of the interest of Leray's mathematical ideas, echoing Eilenberg's comments, and open to the idea that there are things from which he and Cartan could draw.  Weil was aware of the Alexandroff--Čech approach to homology and his opinion seems to reflect Eilenberg's opinion that at an appropriate distance ``local methods'' are already contained in this work.  Borel marks this letter in \cite{LerayNotices} as the source of Cartan's interest in the theory of sheaves, and below we will analyze subsequent correspondence between Weil and Leray, justifying our assertion of Weil as the ``skeptic'' to Cartan as Leray's supporter.  

Looking back at Leray's initial works on sheaves and homology, one sees an approach that contains an amalgam of concepts.  The new notion of sheaf appears simultaneously with new computational tools like the still intimidating notion of spectral sequence.  As we will detail below, the reception of these new ideas was not uniform: some members of the establishment viewed the ideas somewhat suspiciously, perhaps an elaborate dressing-up of ideas that are inherently simpler, others view the techniques hopefully as the claimed applications begin to appear.  

From this point of view, it is no surprise that Leray's approach has been described as obscure.  Note: Leray's work is not obscure in an anachronistic way; Leray still uses ``axiomatic'' methodology and presentation, but instead his notation, terminology, or perhaps something more ineffable, ``style'' perhaps, are not reflective of the mainstream presentations of ideas of algebraic topology, even in France.  As a consequence, it seems completely unclear at this stage whether Leray's ideas would assume a place of importance in mathematics.  What changed people's minds?

The early January 1947 letter from Weil to Cartan marks the beginning of an explosion of correspondence about homology and sheaves that we now have access to.  We can track the evolution of the notion of sheaf by its reflection in these letters.  Later in January \cite[p. 143]{CartanWeil}, Weil writes to Cartan about a putative Bourbaki algebraic topology book: ``the more I reflect on it, the more I think that fiber spaces, differential geometry, etc., cannot be separated from algebraic topology and should constitute only one book...''\footnote{``Plus je r\'efl\'echis et plus je pense que les espaces fibr\'es, g\'eom\'etrie diff\'erentielle, etc. ne peuvent \^etre s\'epar\'es de la topologie alg\'ebrique et doivent constituer un seul livre avec celle-ci...''}.  Weil goes on to suggest that this treatment should be based on open coverings instead of triangulations, in the manner of this previous letter.  
In February he writes with further discussion of de Rham theory, as well as questions about the cohomology theory for groups discussed by Eilenberg--Mac Lane.  Weil sends another letter on February 2 \cite[p. 145]{CartanWeil} clarifying some points of his preceding treatment and making mention of the fact that his proof depends on some lemmas, the first of which (A) is the fact that given a compact differentiable manifold $V$ there exists a covering of $V$ by open sets $A_i$ such that every non-empty intersection is homeomorphic to a ball.  Nowadays we call this a good cover, and this seems to be one of the first places where this notion appears.  The second allows Weil to identify the cohomlogy of $V$ in terms of the nerve of this covering in the sense of Alexandroff-Čech.

Cartan responds shortly thereafter on February 5, saying that he has a thousand things to say.  He questions a key step in Weil's sketched proof of the de Rham theorem, which Weil remarked was a bit delicate in his first letter: cover a given space by subsets all of whose intersections are homeomorphic to balls.  The existence of a ``good cover'' of a differentiable manifold is now frequently attributed to Weil, but Cartan brings this up to support Leray: ``this follows from one of the essential theorems of Leray''.  Cartan then asks: ``do you know this theorem, or are you using something else?'' and writes\footnote{Quant \`a ta d\'emonstration elle-m\^eme, j'ai aperçu sa parent\'e (sans vouloir ta vexer!) avec une d\'emonstration que Leray fait plusieurs fois, sous une autre forme d'ailleurs d\'emonstration sur laquelle j'ai longuemont travaill\'e l'an dernier quand je cherchais \`a debrouiller le m\'emoire de Leray (qui, entre parenth\`eses, pr\'esente plus d'int\'er\^et que tu ne sembles le croire).  J'ai song\'e a t'addressser un rapport d\'etaill\'e sur ce que j'avais tir\'e l'an dernier de ce m\'emoire, mais je crois que c'est maintenant en partie inutile car {\em je viens de trouver bien mieux}, si c'est exact! Je n'ai pas encore toutes les d\'emonstrations en noir sur blanc, mais je n'attends pas davantage pour t'\'ecrire pour ne pas abuser de ta patience.}
\begin{quote}
	As for your proof itself, I saw its relation (without wishing to offend you) to a proof that Leray gives several times, in another form, moreover, a proof which I studied at length last year when I was trying to unravel Leray's memoir (which, incidentally, is more interesting than you seem to believe).  I thought about sending you a detailed report on what I learned last year from this memoir, but I think it is now useless because I have just found much better, if it is correct! I don't yet have all the proofs in black and white, but I didn't want to wait any longer to write to you so as not to try your patience.
\end{quote}
The letter continues with discussion of some purely algebraic notions subordinate to generalizations of the de Rham theorem: differential graded algebras, tensor products of such, eventually returning to Cartan's reconsideration of Leray's notion of {\em couverture}, which he calls a {\em carapace}.  Among these examples is included the ``delicate point'' in Weil's approach to the de Rham theorem.  The letter concludes with announcements of results unifying Leray's presentation, the de Rham theorem and various other approaches to cohomological computations.  

The back and forth continues: Weil further discusses the points on which Cartan asked for clarification.  Cartan draws back some of his claims, and Weil responds with further questions.  On February 18th, Weil continues to be skeptical of Leray's approach linking his opinions to a possible Bourbakiste treatment of algebraic topology \cite[p. 178]{CartanWeil}:\footnote{...mais, pour une th\'eorie bourbachique fond\'ee sur les \guillemetleft carapaces \guillemetright, il y a l'objection de la complication et de l'artificialit\'e--a laquelle je ne sais si la m\'ethode Leray \'echappe (tu dis qu'elle fournit aussi des \guillemetleft carapaces \guillemetright): j'en doute apr\`es Math. Rev.; d'ailleurs, toute theorie qui introduit des complexes g\'en\'eraux (cellulaires) et non pas seulement simpliciaux me para\^it condamn\'ee d'avance, pour un expos\'e bourbachique, \`a cause des horribles matrices d'incidence.}
\begin{quote}
	...but, for a Bourbakiste theory based on ``carapaces'', there is the objection of complication and artificiality--I do not know if Leray's method escapes from this (you say that it also provides ``carapaces''): I doubt it after Math. Rev.; moreover, any theory which introduces general (cellular) complexes and not just simplicial ones seems to me doomed in advance, for a Bourbakiste presentation, because of the horrible incidence matrices.
\end{quote}

Still later in February, Weil's skepticism of the novelty in Leray's approach continues.  In a letter date February 24 \cite[pp. 181-186]{CartanWeil}, which contains a large number of postscripts, he brings up the discussion of Alexander's theory of gratings (which we discussed earlier; see \cite{AlexanderGratingI,AlexanderGratingII}).  The third postscriptum reads:\footnote{D'apr\`es ce que tu me dit la \guillemetleft m\'ethode Leray \guillemetright, la m\'ethode Alexander (gratings) n'est au fond pas autre chose que l'application de cette m\'ethode (avec quelques annees d'anteriorit\'e!) a l'ensemble des recouvrements binaire (je veux dire, form\'es chacun de deux ensembles seulement) de l'espace \'etudi\'e; c'est donc, dans son principe, une {\em simplification} de la \guillemetleft m\'ethode Leray \guillemetright; d'o\`u r\'esulte que si celle-ci (comme tu affirmes) donne des \guillemetleft carapaces \guillemetright, celle que donne la method Alexander doit certainement \^etre plus simples.}
\begin{quote}
	According to what you tell me about``Leray's method'', Alexander's method (gratings) is basically nothing more than the application of this method (with a few years precedence!) to the set of two-fold coverings (I mean, each consists of only two sets) of the space studied; it is therefore, in principle, a simplification of ``Leray's method''; from which it follows that if the latter (as you assert) gives ``carapaces'', the one given by Alexander's method must certainly be simpler.
\end{quote}
This postscriptum ends with the Weil's vivid characterization of the Leray's approach:\footnote{En somme, le monstre Leray serait un rejeton b\^atard issu de l'union incestueuse du cube et du simplexe (ou encore: d'Alexander et e Cech)!!!}
\begin{quote}
	In short, the Leray monster would be the bastard child formed of the incestuous union of the cube and the simplex (or even of Alexander and Cech)!!!
\end{quote}

The sixth \cite[p. 186]{CartanWeil} returns to Leray, reflecting on the conversation Weil had with him in 1945 around the notion of local coefficients:\footnote{Si c'est bien l\`a ce que j'ai cru comprendre en 1945, c'est l\`a une id\'ee tr\`es f\'econde, mais d;une telle g\'en\'eralite et d'une telle port\'ee qu'on ne peut gu\`ere esp\'erer la mettre en {\oe}uvre, pour l'instant, que dans des cas particuliers...{\em Si} c'est l'id\'ee de Leray, et {\em si} il l'a d\'ej\`a mise en {\oe}uvre, il serait de premi\`ere importance pour nous de savoir o\`u il en est; sinon il faudrait explorer nous-m\^emes le terrain, bien que ça semble difficile.}
\begin{quote}
	If this is indeed what I thought I understood in 1945, it is a very fruitful idea, but of such generality and such scope that we can only hope to implement it, at the moment, in particular cases...If it is Leray's idea, and if he has already implemented it, it would be of primary importance for us to know where he stands; otherwise we would have to explore the terrain ourselves, although that seems difficult.
\end{quote}
Weil goes on to speculate once again about using ideas of local-coefficients to move toward Morse theory, and includes the telling footnote:\footnote{Mine inepuisable de sujets de these, pour eleves (astucieux) de Bourbaki.} ``An inexhaustible mine of thesis topics for (clever) Bourbaki students.

On March 4th, Cartan responds with a long letter.  He discusses a number of things before writing, ``and now we come to homology'' responding to a number of Weil's criticisms \cite[p. 191]{CartanWeil}, in part around how Bourbaki should approach homology, but also about Leray.  First, I mention for context that Cartan points out that the Eilenberg--Steenrod axiomatization of homology only establishes certain kinds of uniqueness statements, but of a manifestly different sort than what one would like.  Cartan references ``complexes'' presumably of the combinatorial sort, which is at odds with Weil's point of view about what a proper ``Bourbaki'' account of homology might look like.  Cartan continues:\footnote{D'ailleurs tu sembles un peu avoir perdu de vue l'importance de ce th\'eor. d'unicit\'e, t'\'etant surtout attach\'e \`a la question d'existence, que tu croyais devoir r\'esoudre par la n\'egative en g\'en\'eral.  Je pense que maintenant tu es convaincu de cette existence, pusique tu as retrouv\'e Leray chez Alexander! Mais, si je juge Alexander par ce que tu m'en dis, la diff\'erence avec Leray (ou plut\^ot, ce que je t'ai dit de Leray) est bien minime; et cela ne veut pas la peine de s'exciter pour savoir si on proscrira enti\`erement les complexes simpliciaux pour ne plus tol\'erer que les sous-complexes d'un cube (car c'est cela que tu r\'eclames avec tes recouvrements binaires d'Alexander).  Dire que dans ta lettre ant\'erieure (de 18), tu voulais bannir tous les complexes autres que simpliciaux, pour antibourbachisme!  Tu passes un peu d'un extreme et \`a l'autre.  Mon opinion est qu'on a {\em besoin} des complexes simpliciaux {\em et de leurs produits}; avec cela, on peut tout faire...C'est e\'videmment beaucoup plus que qu'il n'en faut.}
\begin{quote}
	Besides, you seem to have lost sight a little of the importance of this uniqueness theor., having been especially attached to the existence question, which you believed you had to resolve by the negative in general.  I think you are now convinced of this existence result, since you found Leray's ideas already in Alexander's work! But, if I judge Alexander by what you tell me, the difference with Leray (or rather, what I told you about Leray) is very minimal; and there is no need to get excited about whether we will completely ban simplicial complexes and only tolerate the subcomplexes of a cube (because that is what you are demanding with your binary Alexander coverings) .  To think that in your previous letter (of the 18th), you wanted to ban all complexes other than simplicial, for anti-Bourbachism!  You go from one extreme to the other.  My opinion is that we {\em need} simplicial complexes {\em and their products}; with this, you can do anything...It's obviously much more than is necessary.
\end{quote}
Cartan continues with an implicit admonishment of taking the content of Math Reviews of long memoires too seriously:\footnote{Tout de meme, il est regrettable que tu ne puisses pas consulter le gros memoire de Leray, car les Math.Reviews ne suffisent pas: en general, les comptes rendus des memoires importants y sont si mal faits que le lecteur ne peut m\^eme pas se douter qu'il passe \`a c\^ot\'e de quelques chose d'important.}  ``All the same, it is regrettable that you cannot consult Leray's big memoir, because the Math.Reviews are not enough: in general, the reports of important memoirs are so poorly done that the reader cannot even understand that he is missing out on something important.''

Cartan continues to give further explanations for the importance of Leray's work. Could it be that he knows he cannot change Weil's mind on his own? After mentioning some results of Gysin, Cartan continues:\footnote{Leray pr\'etend qu'il a int\'egr\'e cela dans sa th\'eorie des coefficients variables, il a d\'ej\`a paru des Leray, l'\'et\'e dernier, 3 ou 4 notes tr\`es compliqu\'ees et difficilement pigeables; il faudrait plusieurs mois pour d\'ebrouiller cela.  Il obtient des resultats que je ne suis pas a meme d'apprecier, mais qui ont vivement int\'eress\'e Hopf l'\'et\'e dernier (Hopf ayant auparavant s\'ech\'e vainement sur les memes questions).} ``Leray claims that he integrated this into his theory of variable coefficients; Leray already published 3 or 4 very complicated and difficult to understand notes last summer; it would take several months to sort this out.  He obtains results that I am not in a position to appreciate, but which greatly interested Hopf last summer (Hopf having previously tried and failed with the same questions).''  

Cartan continues on March 11th \cite[p. 198]{CartanWeil}, having looked closely at Alexander's work on ``gratings'' to ostensibly defend Leray's originality, but ceding some ground on the precedence of Alexander's approach:\footnote{Mon opinion actuelle sur les \guillemetleft gratings\guillemetright: c'est un {\em proc\'ed\'e} particulier pour d\'efinir un anneau de cohomologie, tout comme celui de Leray (avec \'evidemment plusieurs ann\'ees d'ant\'eriorit\'e), mais presque tout le m\'emoire est consacr\'e \`a un expos\'e alg\'ebrique du proc\'ed\'e, et il en tire peu de chose, pas m\^eme la cascade (et pourtant, c'est de cette mani\`ere qu'elle se pr\'esente tr\`es intuitivement).  Leray lui, en tire des tas de choses, ne serait-ce que son th\'eor\`eme sur la d\'etermination de l'anneau d'homologie par celui d'un recouvrement \`a \guillemetleft supports\guillemetright simples.  D'autre part, pas plus que Leray, Alexander ne semble vouloir se donner la peine de confronter ses d\'efinitions avec d'autres d\'efinitions existantes, telles que celle de Cech.  Enfin, il se trouve que les \guillemetleft gratings\guillemetright, comme les trucs de Leray, donnent des {\em carapaces}, mais il y en a d'autres comme celles de De Rham avec les formes diff\'erentielles.  Donc, si Alexander (ou Leray, ou une synth\`ese des deux) peut et doit servir \`a prouver l'existence de carapaces (sup., bien entendu), le th\'eor\`eme d'{\em unicit\'e} de l'anneau d'hom. des carapaces est une pi\`ece essentielle, qu'on ne trouve pas chez Alexander ni ailleurs, pas plus que chez Eilenberg; et il donne imm\'ediatement, de la mani\`ere la plus simple, le th\'eor\`eme de la cascade.}
\begin{quote}
My current opinion on ``gratings'': it is a particular process for defining a cohomology ring, just like that of Leray (obviously several years prior), but almost the entire work is devoted to an algebraic presentation of the process, and he gets little out of it, not even the cascade (and yet, this is how it presents itself very intuitively).  Leray draws lots of things from it, if only his theorem on the determination of the homology ring by that of a covering with simple ``supports''.  On the other hand, no more than Leray, Alexander does not seem to want to take the trouble to compare his definitions with other existing definitions, such as that of Čech.  Finally, it turns out that the ``gratings'', like Leray's tricks, give ``carapaces'', but there are others like those of De Rham with differential forms.  So, if Alexander (or Leray, or a synthesis of the two) can and must be used to prove the existence of carapaces (sup., of course), the uniqueness theorem for the homology ring of a ``carapace'' is an essential piece, which one does not find in Alexander or elsewhere, nor in Eilenberg; and he immediately gives, in the simplest way, the cascade theorem.
\end{quote}
Here the term ``cascade'' is a reference to what we now call a ``long exact sequence''.  It is difficult to decipher exactly which ``long exact sequence in cohomology'' is being considered, but given the discussion of local coefficients, it is natural to speculate that they are speaking of the long exact sequence in cohomology attached to something like a short exact sequence.  Cartan follows this up soon thereafter with a letter on the 12th of March \cite[p. 200]{CartanWeil} where he describes another benefit of the theory of carapaces it\footnote{...entra\^ine le th\'eor\`eme vainement cherch\'e par Eilenberg et aussi la d\'etermination de l'anneau d'homologie d'un espace-produit.  Je suis persuad\'e que la m\'ethode Alexander--Leray est la bonne, aux d\'etails de pr\'esentation pr\`es.  Ce sont ces d\'etails que d\'ecideront de sa \guillemetleft simplicit\'e \guillemetright ou de sa \guillemetleft complication \guillemetright apparentes.}
\begin{quote}
	...results in the theorem Eilenberg sought in vain and also in the determination of the homology ring of a product space.  I am convinced that the Alexander-Leray method is the right one, up to some details in the presentation.  It is these details that will decide its apparent ``simplicity'' or ``complication''.
\end{quote}
Cartan also draws a distinction between the ideas of Čech, which ``seduced him'' and the method of Alexander--Leray, to which we will return shortly.

The subsequent discussion turns to other topics, but there are still a few remarks to be made.  On the 6th of April \cite[p. 215]{CartanWeil}, Weil writes:\footnote{O\`u en es-tu de l'homologie?  Il me parait essentiel que tu r\'ediges le plus t\^ot possible, tr\`es br\`evement mais compl\`etement, ta th\'eorie carapatique, avec, non seulement, le th\'eor\`eme d'unicit\'e mais le th\'eor\`eme d'existence; sans quoi nous n'avancerons pas.}\begin{quote}
	Where are you on homology?  It seems essential to me that you write down as soon as possible, very briefly but completely, your carapatic theory, with, not only the uniqueness theorem, but also the existence theorem; without which we won't press forward.
\end{quote}
Cartan's response to further questions in this letter is also striking because it gives a first hint of the influence of Eilenberg's categorical point of view on his thinking and presentation \cite[p. 222]{CartanWeil}.  \footnote{Je r\'eponds maintenant, au moins en partie, aux questions de ta lettre homologique du 6 avril,--sans attendre une r\'edaction complete de ma part!...Mais ce supplement au th\'eor\`eme d'unicit\'e (\`a savoir que les repr\'esentations en questions sont \guillemetleft naturelles\guillemetright comme dit Eilenberg) est facile...}
\begin{quote}
	I now respond, at least in part, to the questions in your homological letter of April 6,--without waiting for a complete redaction on my part! ...But this supplement to the uniqueness theorem (namely that the representations in question are ``natural'' as Eilenberg says) is easy..
\end{quote}
Cartan continues, now confident that his point of view allows a comparison of the Čech and Leray approaches to homology.

Finally, on the 2nd of May, Weil writes back \cite[p. 225]{CartanWeil}:\footnote{Remercie Koszul pour les notes de Leray...Les d\'efinitions de ces notes sont exactement ce que j'avais reconstitu\'e de mon cot\'e; les r\`esultats qu'il y donne ne semblent pas aller tr\`es loin pour l'instant, et ce qu'il dit des applications qu'il a en vue est un peu d\'ecevant aussi; mais l'id\'ee est \`a retenir.}
\begin{quote}
Thank Koszul for Leray's notes...the definitions in these notes are exactly what I had reconstructed on my own; the results he gives do not seem to go very far for the moment, and what he says about the applications he has in mind is also a little disappointing; but the idea is worth remembering.
\end{quote}
In a letter on the 9th of May \cite[p. 227]{CartanWeil}, Weil continues:\footnote{Je me suis mis \`a l'\'etude des notes Leray; c'est fort int\'eressant, d'ailleurs pas difficile, et je t'en recommande l'\'etude. Ce qui est essentiellement nouveau, ce sont les homomorphisms canoniques d\'efinis dans la 2e note (N.B. j'attends pour bient\^ot les suivantes)}
\begin{quote}
	I began to study the Leray notes; It's very interesting, and besides not difficult, and I recommend you study it. What is essentially new are the canonical homomorphisms defined in the 2nd note (N.B. I am waiting for the following ones soon).
\end{quote}
I think I laughed out loud the first time I read Weil telling Cartan to study Leray, but I can only imagine how Cartan must have responded to this nudge.  He continues in a postscript:\footnote{Bien entendu, le fait que Leray utilise, dans ces notes, sa technique des \guillemetleft couvertures\guillemetright (sic) n'est pas un argument en faveur de celles-ci; il n'y a pas de difficult\'e a transposer tout ça dans n'importe quel syst\`eme (Cech, Alexander, \guillemetleft carapaces\guillemetright) qu'on pourra d\'ecider d'adopter.  Il ne me semble pas, pour l'instant, qu'il y ait un rapport quelconque entre Leray et Gysin (sauf que Leray affirme qu'il retrouve les resultats de G. relatifs aux fibr\'es a fibre sph\'erique). }
\begin{quote} Of course, the fact that Leray uses, in these notes, his technique of ``couvertures'' (sic) is not an argument in favor of these; there is no difficulty in transporting all of this into any system (Čech, Alexander, ``carapaces'') that we may decide to adopt.  It does not seem to me, for the moment, that there is any connection between Leray and Gysin (except that Leray claims that he deduced G.'s results relating to sphere bundles).
\end{quote}
The letters continue, but the parties seem to have drawn their sides.  I was struck at just how hard Weil works, over an extended correspondence to assert that Leray's point of view is really not new.

In the meantime, the reception of these ideas in the United was somewhat different, and this is discussed in detail in \cite{MillerLeray}.  George Whitehead stated that he found Leray's writing obscure and William Massey said: ``“In the late 1940’s and early 1950’s all of us were studying Leray’s papers to try to understand how he got the marvelous results he claimed.  To be perfectly frank, I never got to 1st base in this enterprise, it was very frustrating. Leray was a horrible expositor.”; Evidently this is at odds with Weil's assessment of the results.  Miller goes on, quoting an interview with Leray where he remarks on the distress he felt by the reception of his work in the United States, drawing particular attention to the Math. Reviews which he felt did not accurately represent their content. 

What are we to take away from all of this?  In the face of Weil's criticisms, Cartan continued to investigate Leray's ideas, was a proponent for them, and decided to internalize them.  What would have happened to the theory of sheaves without the support of Cartan and Weil's challenges?  In retrospect, the differences between all of these approaches have disappeared.  Cartan's theory of ``carapaces'' and Leray's original theory of sheaves have been subsumed in another more ``modern'' theory.  The variant different approaches have been unified with Leray's ideas of ``spectral sequence'' providing a useful tool for comparing all the different approaches, but none of this seems evident to me at this stage.  Indeed, in order for the theory to develop, it seems necessary for the tools to find an audience.  

\subsubsection*{Sheaves in the Cartan seminar}
In the meantime, at the instigation of Cartan, J.-L. Koszul was investigating Leray's work.  Cartan already remarks on Koszul's works in the letter of April 20th \cite[p. 223]{CartanWeil}, where he indicates that Koszul has recovered ``all of Leray's results and those of some others'' in the case of ``representations'' of a compact Lie group to a homogeneous space.  These results eventually appeared as \cite{Koszul1}, which isolates the algebraic structures that Leray considered--in modern terminology aspects of the theory of filtered differential rings--, and \cite{Koszul2} which announces the claimed applications to homogeneous spaces. 

In the summer of 1947, the CNRS hosted an international colloquium on algebraic topology in Paris.  Leray gave a talk at this conference, and Cartan comments directly on this talk \cite[p. 244-5]{CartanWeil}:\footnote{Leray ne m'a rien appris d'essentiellement nouveau, sauf peut-\^etre le suppl\'ement qu'il a fait pour un cours du colloque et qui se rapporte aux espaces dans lesquels op\`ere un groupe d'automorphismes.  Il faudra que je tire au clair les quelques indications sybyllines qu'il donne; ces questions se rattachent au m\'emoire d'Eilenberg des Transactions, que j'ai pu lire rapidement avant de quitte Strasbourg.}\begin{quote}Leray taught me nothing essentially new, except perhaps the addition he made for his talk at the conference and which relates to spaces on which a group operates by automorphisms.  I will have to clarify the few cryptic indications he gives; these questions relate to Eilenberg's memoir of the Transactions, which I was able to read quickly before leaving Strasbourg. \end{quote} Here, Leray makes reference to Eilenberg's paper \cite{EilenbergActionsI} analyzing homology of spaces on which a group acts; he uses Steenrod's theory of local coefficients and makes connections to his earlier work on group cohomology.  In his colloquium talk, Leray applied his theory to the analysis of finite (Galois) coverings.  Cartan's ``clarifications'' amounted to an extension of Leray's results to infinite coverings and were announced in two CRAS notes \cite{Cartan48I,Cartan48II}.  

In 1947-48, Leray gave a course on his approach to cohomology via sheaves at the College de France; the notes from this course were eventually published in a polished form in 1950 \cite{LerayMain}.  Cartan visited Harvard in 1948 and gave a course on algebraic topology.  The notes of this course were written up by George Springer and Henry Pollak \cite{CartanHarvard}, and it is unclear what influence Cartan had on the final treatment (see \cite[Footnote 4]{MillerLeray}).  The aim of the course was to end up at Weil's proof of the de Rham theorem.  According to \cite{MillerLeray}, the resulting treatment makes no mention of Leray or Weil.  Moreover, Chern's MathSciNet review of this course makes no mention whatsoever of Leray's work on ``couvertures'', though it does make mention of Alexander's theory of gratings.  In view of the amount of discussion of the relative importance of the ideas of Alexander and Leray in the Cartan--Weil correspondence, one wonders how Cartan reacted to this lack of attribution. 

Shortly after the correspondence above, in the 1948-49 academic year, Cartan devoted a portion of the S\'eminaire Cartan to analyze sheaves and Leray's approach to homology.  The write-ups from these expose\'es were later withdrawn as Cartan's viewpoint changed, and only the portion of the seminar devoted to other topics in algebraic topology is still available.  We draw a few comments from the exposes that do remain.  Leray's name is explicitly mentioned in the expose on Čech--Alexander cohomology \cite[Expose 6]{SemCartan194849}.  In the expos\'e on local coefficients \cite[Expose 10 \S 11]{SemCartan194849}, one finds a definition of ``local system of groups'', highlighted as a new definition; this is the modern notion of {\em presheaf}.  Here ``local coefficient systems'' are used as coefficients for Čech cohomology.




We will continue this development below, but while we keep referring to Leray's theory of sheaves, it is telling that Leray published essentially nothing about sheaves or homology after 1950 (save an erratum to his various works discussed above).  The comments above about Leray's frustrations with the reception of his work, the rapid internalization of the theory by Cartan and the subsequent developments of Leray's tools by other authors, e.g., spectral sequences and Koszul as we described above, make it clear that the theory was very quickly no longer Leray's.  One is led to ask the question: was this for reasons of taste or frustration?  Instead, the theory of sheaves was subject to significant and repeated re-invention in the 1950s.

\subsubsection*{Interlude: on the mathematical mainstream}
Leray's work was very much a topic of discussion among a certain group of mathematicians in the 1940s and 1950s: one might call this group mainstream French mathematicians interested in algebraic topology and related subjects.  This treatment of Leray's work on sheaf theory, I hope, demonstrates that there is a real sense in which Leray never entered this community, as evidenced by his mathematical interests moving in other directions.  Leray was in some ways an outsider to this community to begin with, and it seems unquestionable that issues of taste played a role in this outcome.  We hear from various people that Leray's work is ``difficult'', ``obscure'' or that Leray is a ``poor expositor''.  But one could just as well argue that there was a simple disjunction between Leray's intuitive style and the Bourbaki-style that was sweeping across the landscape of algebra and algebraic topology where Leray's ideas found their audience.  Mawhin describes Leray's writing in the following way: ``With his caustical mind, which reminds Voltaire or Henri Poincar\'e, with his love of French language, whose severe elegance follows from a constant care
for concision...''.  Meyer echoes this description \cite[\S 5]{MeyerLeray}:\footnote{Pour Jean Leray, la beaut\'e et la force des math\'ematiques suffisaient. Il \'ecrivait donc soigneusement au tableau ce qu’il voulait transmettre \`a ses auditeurs. Son style \'etait sobre, sec, sans concession et les auditeurs de ses cours \'etaient parfois perdus. }
\begin{quote}
	For Jean Leray, the beauty and strength of mathematics were enough. He therefore carefully wrote on the board what he wanted to convey to his listeners. His style was sober, dry, uncompromising and the listeners of his courses were sometimes lost.
\end{quote}

And thus we come back to Leray's independent-mindedness \cite{Mawhin2}:
\begin{quote}
	The essential character of my publications is however their diversity: the problems which attracted me required techniques still unused in the specialty where they were classified; their solution required or suggested improvements of those techniques, such that, to develop them, I had to change my area.
\end{quote}
While Leray won a number of mathematical honors, he gave only a contributed paper in the 1950 ICM, and was never a plenary speaker at any of the ICMs that took place during his lifetime \cite{Mawhin2}. 

I have already mentioned Tucker's review \cite{TuckerSFBR} where an awareness of the existence of a mainstream in mathematics, or at least an opinion that a mathematical mainstream exists is professed.  The discussion of mainstream mathematics, especially around sheaf theory surrounds Bott's recollections of his second visit to the IAS \cite[p. 532]{Duren2} from 1955-1957:
\begin{quote}
	if my previous stay was dominated by a feeling of awe at the brilliance of the older generation, during this second one, my dominant feeling was one of awe, alas, mixed with envy, at the brilliance of a generation younger or contemporary to my own! 
	
	And no wonder. During that time, and largely at Princeton, I met Serre, 
	Thorn, Hirzebruch, Atiyah, Singer, Milnor, Moore, Borel, Kostant, Harish- 
	Chandra, James, Adams,..., 
\end{quote}
He speaks of the ``great topologist, Serre'' and how Serre tried to teach him sheaf theory.  By 1950, sheaf theory was very much part of the culture of ``mainstream mathematics'' about which Bott also recorded some interesting things.

Baum wrote, in volume 3 of the collected works \cite{BottCollected3}:
\begin{quote}
	What did I learn from Bott? I learned--perhaps I should say relearned--from Bott things that I already knew: that there is great beauty and wonder and fun in mathematics; that the best mathematics is not necessarily the most difficult; that there is a mainstream to our subject that talented mathematicians like Bott have an instinctive feel for; that ideas should be presented clearly and without nonsense so that other mathematicians can understand and appreciate them.
\end{quote}  Allyn Jackson discusses this in her interview with Bott \cite{JacksonBott}.  Of Baum's comments she writes: ``How does one come to understand what that mainstream is? You're born with it? You learn it? You pick it up from the environment?''

For Bott, finding the mainstream was always a matter of following his tastes.  ``And sometimes my taste led me in directions that weren't fashionable but that luckily turned out to be fashionable later on! But these things are dangerous, the fashions change, and it's hard to tell in retrospect whether you were in the mainstream. I was just very affected by the early development of sheaf theory, and especially the combination of analysis with topology that then ensued. Suddenly complex variable theory fitted in with topology and even certain aspects of number theory. So I think that at that time it was very easy to discern this development as a main road in mathematics.''  While we haven't gotten there yet, this discussion of ``mainstream sheaf theory'' seems entirely divorced from the theory that Leray was describing.

Bott continues: ``But I've seen the mainstream change considerably over my lifetime. For instance, if I think of Princeton before sheaf theory, the emphasis was very different. When I first came there, much of topology in those early years had to do with very abstract questions of pathological spaces, comparing fifteen different cohomology theories, and such. This was what I would have said at first was the mainstream. Then topology moved more to what I felt was the real world: the study of compact manifolds and their invariants. Lower-dimensional topology was not emphasized then, but
in the 1990s it came to the fore again. So there is really a tremendous difference in perspective over the years.''

And finally, we come to the question which touches on the inevitability of this kind of mathematics: ``But isn't there a core of mathematics that is vital and lively, independent of fashion, and there are other fields that are more outlying; and one needs a sense of what is central and what perhaps is not so central?''  To which Bott responds: ``I don't know to what extent I believe this.  I think, for instance, that Bourbaki had that feeling, and I was always a little skeptical of Bourbaki.  The subject is just too big. It doesn't just have one main road. There are too many unsuspected branches. So although I was in a sense very much influenced by Bourbaki, I don't really subscribe to the belief that there is just one way of looking at things...If you had one really good main highway, it would be dangerous, because then everybody would be marching along it!''  

It is worth contrasting this opinion to that put forth by Serre in an interview from 2004 \cite{SerreNotices}. Serre was asked directly: ``Do you feel that there are core or mainstream areas in mathematics? Are some topics more important than others?'' He responds: ``A delicate question. Clearly, there are branches of mathematics which are less important: those where people just play around with a few axioms and their logical dependencies. But it is not possible to be dogmatic about this. Sometimes a neglected area becomes interesting and develops new connections with other branches of mathematics.''  Serre goes on to mention the Riemann hypothesis and the Langlands program as ``obvious'' cases of questions that are central to our understanding of the mathematical world.  


\subsubsection*{Cartan and analytic sheaves}
The War also marked a turn for Cartan and we discuss the implicit question from the Cartan--Weil correspondence: why was Cartan so receptive to Weil's approach to the de Rham theorem (effectively immediately lecturing on it at Harvard) but simultaneously less receptive to Weil's opinions about sheaf theory? We have already mentioned Weil's pushing Cartan to integrate ``bundle-theoretic thinking'' into his approach (see \S\ref{paragraph:WeilFB}).  

In 1944, Cartan published the second part of his work on the Cousin problem \cite{CartanIdeals} and, we recall some aspects of Chorlay's analysis  \cite[\S 4]{Chorlay} of this period.  Chorlay refers to this as the beginning of Cartan's ``structural'' transition.  Cartan remarks on this paper to Weil in a letter of December 10, 1944 \cite[p. 96]{CartanWeil} explaining that it answers questions posed to him by Weil in December 1940. This paper embarks on a project of a ``global'' study of ideals of holomorphic functions.  Cartan embeds the multiplicative Cousin problem into a problem of {\em patching} invertible matrices with holomorphic coefficients:\footnote{pour passer de donn\'ees locales \`a une existence globale, on proc\`ede par assemblages successifs de morceaux. Chaque stade d’assemblage consiste en ce que nous appellerons une op\'eration \'el\'ementaire.}
\begin{quote}
to pass from local data to a global existence, parts are assembled in turn. Each assembling step consists in what we will call an elementary operation.
\end{quote}
Chorlay argues that the algebraic structures used in Cartan's theory are a visible example of the structuralist ideal of Bourbaki (viz. the theory of ideals and modules being imported wholesale into holomorphic function theory).  Chorlay also highlights Cartan's use of ``coherent systems'' where algebraic data are seen varying with the points on a space.

In the 1944 letter to Weil discussed above, Cartan also mentions that his interests have turned a bit to algebraic topology, especially the notion of homology, and that he's curious if the conclusions he's drawn are concordant with Eilenberg's views. In 1945, he publishes \cite{CartanModern}, which Steenrod reviews for Math. Reviews asserting: ``The interest lies chiefly in the methods used. Once the general properties of the Lefschetz groups are established, all subsequent results are deduced from these properties alone, thus avoiding any additional use of complexes or assumptions of triangulability. These basic properties are almost exactly the ones used by Eilenberg and Steenrod in their axiomatic development of homology theory'' closing with reference to \cite{EilenbergSteenrodPNAS}.  Cartan's published work in these years seems quite far from algebraic topology, for example, he publishes numerous papers on potential theory, but in the Cartan--Weil correspondence we see Cartan grappling at the same time with sheaves.  

\subsubsection*{Conclusion}
The years 1950-1952 show the theory of sheaves rapidly diverging from path initiated in Leray's approach. Cartan's synthesis of the bundle theoretic approach to the Cousin problems espoused by Weil, his understanding of Leray's theory of sheaves, and developments of algebraic topology play out in the seminar in 1949-1952.  The aborted 1948-49 treatment of sheaf theory is revisited in the 1950-51 version of the S\'eminaire Cartan providing a first approximation to Cartan's ``definitive'' treatment of sheaf theory.  

In this ``second generation'' approach (\cite[p. 3]{Chorlay}) Leray's definitions have effectively been excised.  In their place, one finds a new more geometric definition of sheaf due to Lazard.  This new definition allows one to retain the idea of sheaves as coefficients, which we saw was key to Weil's initial fascination, but expands the purview of sheaf theory to allow sheaves as complementary geometric objects to fiber bundles.  The key point I'd like to make here, echoing, and I believe extending some aspects of Chorlay's analysis, is that this conjunction of ideas: sheaves as coefficients, fiber bundles, and sheaves as way-station in formulating ``global'' problems in complex analysis is {\em unexpected}.

Simultaneously, tremendous use was made of Leray's spectral sequences in the hands of Koszul, Borel and Serre; these results have been well-documented elsewhere so I will not say much about them here.  At the end of this period, we see that Leray's influence has been transcended as the notion of sheaf is subsumed in a paradigm evolving in the Cartan seminar.  

In the background, we see academic squabbles, perhaps anodyne to the French reader, but which to my modern American sensibilities, read like tales of intrigue from courtly salons.  Rather than viewing this as mere pettiness, I feel that this tension led to a particular style of treatment of Leray's work because of the way it plays out in the Cartan--Weil discussion.  Furthermore, the political drama in the background, brings to the fore questions such as whom should we credit for what discoveries and how much should we credit them!  What is said seems to matter just as much as who says it!  One wonders what Leray's role in algebraic topology would have been if the politics had been different: would the reception of his ideas would have played out differently?
  
\newpage
\section{Projective modules proper} 
\label{s:projective} 
\addtoendnotes{\vskip 1.5em \noindent \begin{large}{\bf Projective modules}\end{large} \vskip .2em}
With the incantation: ``a module $P$ is projective if given any homomorphism $f: P \to A''$ and any epimorphism $g: A \to A''$, there exists a homomorphism $h: P \to A$ with $gh = f$'' the ``notion of projective module'' was invented by Cartan and Eilenberg in their influential treatise on homological algebra \cite[p. 6]{CE}.\footnote{Nowadays, I think it is standard practice to isolate a definition in its own environment, especially if one is isolating a new concept or one that plays a central role in the text.  This kind of separation did already occur in the 1950s: it appears in \cite{EilenbergSteenrod}, but it does not appear in Cartan--Eilenberg.}  While this text was published in 1956, it was written earlier. The preface to \cite{CE} indicates that the text was submitted to Princeton University Press by September 1953, but Cartan states \cite{CartanNotices} that he has no idea why it only appeared in 1956.\footnote{Mac Lane opens his review of \cite{MacLaneCEReview} with the amusing:
	\begin{quote}
		At last this vigorous and influential book is at hand. It took nearly three years from completed manuscript to bound book; Princeton is penalized 15 yards for holding. 
	\end{quote}
	He continues with
	\begin{quote}
		In spite of the delay in its publication, widespread acquaintance with the manuscript and with the ideas of this book has already played an important role in the development of this lively subject.
	\end{quote}
}  

Weil writes to Cartan \cite[p. 329]{CartanWeil} indicating that the 150 page manuscript he received, evidently a preliminary version of \cite{CE} could be published by Hermann.  Thus, a rather detailed treatment already existed by 1951 and apparently enough to start considering where the resulting monograph should be published.  A preliminary Bourbaki report on ``homological algebra'' \cite{BourbakiHA178}, notably written in English, also contains a treatment of aspects of homological algebra, once more before the official publication of Cartan and Eilenberg's text, but explicitly mentioning projective and injective modules. 

Preliminary treatments of the ideas in the Cartan--Eilenberg text also appeared in the S\'eminaire Cartan and there is another Bourbaki redaction of an approach to homological algebra \cite{BourbakiHA130}, this time written in French and making no mention of projective modules.  Instead, this text contains a preliminary notion of a complex of projective modules: a graded (left) module over a ring equipped with a graded projection operator, makes no mention of Cartan--Eilenberg (so one imagines that it predates the inception of that text), but seems not to come equipped with a precise date of preparation.

At the very least, making precise the date of invention of the notion of projective module seems much more complicated,\endnote{Mac Lane includes \cite[p. 623]{MacLaneCEReview} a fascinating personal opinion about editorial choices made by Cartan-Eilenberg as regards the material to be included in their book:
	\begin{quote}
		The reviewer might also add his strictly personal opinion that the authors have not kept sufficiently in mind the distinction between a research paper and a book : a good research paper presents a promising new idea when it is hot—and when nobody knows for sure that it will turn out to be really useful; a good research book presents ideas (still warm) after their utility has been established in the hands of several workers. This book contains too large a proportion of shiny new ideas which have nothing to recommend them but their heat and promise...The reviewer is not claiming that...these...notions will not later have significant...uses: some of them will, but until that time comes their presence clutters up the book. Another danger of shiny new notions is that sometimes the shine proves illusory. 
	\end{quote}
	Where do projective modules fall on the spectrum suggested by this description?}, but we will see that by 1956 the ideas contained in the Cartan and Eilenberg treatise were well accepted in a wide circle around these figures.  Having set the stage with a discussion of the mathematical ``cultures'' before 1951 into which projective modules were born, we now focus on the interactions of these various cultures and their co-evolution in the early 1950s.  Key among these will be further attempts to bridge a perceived ``gap'' between sheaf theory and Cartan--Eilenberg's homological algebra.

As discussed above, Cartan and Eilenberg defined projective modules in \cite[p. 6]{CE}.  In fact, the notion appears before this in the preface \cite[p. vii]{CE}.  I no longer remember when I first read the definition of projectivity, but on looking at the definition again recently I had the experience of a familiar word becoming unfamiliar.  

Cartan and Eilenberg aim to introduce/motivate their work: ``During the last decade the methods of algebraic topology have invaded extensively the domain of pure algebra, and initiated a number of internal revolutions. The purpose of this book is to present a unified account of these developments and to lay the foundations of a full-fledged theory'' \cite[p. v]{CE}.  This phrasing strengthens Eilenberg's description from \cite{EilenbergBull}: `` The method of study is also purely algebraic but is the replica of an algebraic process which has been widely used in topology, thus the words ``topological methods'' could be replaced by ``algebraic methods suggested by algebraic topology.''''  The ``invasion'' was on three fronts including cohomology theories for groups, Lie algebras and associative algebras.\endnote{Rather than go into 
	detail about these notions, we mention here only some references.  In the spirit of the work of Cartan--Eilenberg, arguably the ``earliest'' of these novel cohomology theories is that for Lie algebras as it was studied by Elie Cartan \cite{CartanELie}: he considered a complex of left-invariant differential forms on a compact connected Lie group in order to determine its Betti numbers (here Poincar\'e polynomial) of such a space.  Such computations were very much of interest in the 1920s and 1930s.  Cartan's result was also a precursor to de Rham's theorem allowing a description of Betti numbers of manifold in terms of differential forms. As we have seen de Rham's theorem was a motivating force in all of the developments of axiomatic cohomology theories because it did not fit into the framework of initial axiomatizations.  Cartan's theory underwent a first modernization in \cite{ChevalleyEilenberg} and later in \cite{KoszulLie}; the latter also containing a good amount of information regarding Koszul's understanding of Leray's theory of spectral sequences.
	
	The (co)homology theory for groups mentioned here was the cohomology theory of discrete groups invented by Eilenberg--Mac Lane \cite{EilenbergMacLaneCollected}, which we discussed as inspiration for the invention of functoriality.  This theory had a resurgence in the early 1950s as it was coupled to arithmetic ideas arising in Galois theory and class field theory.
	
	The (co)homology theory for associative algebras was that of Hochschild \cite{Hochschild1,Hochschild2}, building on earlier results of his thesis, which eventually was published as \cite{HochschildSS}.  This theory is rather special, as Hochschild observes in \cite{Hochschild1}: low-dimensional cases were apparently studied in unpublished results of J.H.C. Whitehead, and cohomology is {\em degenerate} in the sense that higher cohomology is determined by what happens in degree $1$.  Hochschild mentions a cohomological characterization of separability of an algebra.
	
	In modern terms, all of these (co)homology theories are obtained by taking (co)homology of explicit complexes, and the ``unification'' that Cartan and Eilenberg speak of amounts to treating them all as derived functors.}  

The authors stand ``in part'' upon K\"unneth's works on homology of a product space \cite{Kunneth1,Kunneth2}.  As discussed earlier, K\"unneth's description of the homology of a product space was given in Poincar\'e's style, separating Betti numbers and ``torsion numbers''.  Cartan and Eilenberg describe their goal as ``stating these results in a group-invariant form'', which is vague enough, but granted their introductory remark above, can be viewed as code for presenting this result within their unified treatment.  As a ``first step'' in this direction the goal is to study ``the homology groups of the tensor product of two (algebraic) complexes.''  The aim is to describe the homology of this tensor product in terms of homology of the contributing complexes, but they quickly observe that the problem involves in addition ``a second product called their torsion product. The torsion product is a new operation derived from the tensor product.''

I'd like to begin by attempting to reconstruct how this must have looked in 1953.  One can imagine two vantage points: those who were participating in the Cartan seminar or in the orbit of Eilenberg and those that were not.  In the course of defining the torsion product, Cartan and Eilenberg use ``free resolutions'' before immediately commenting that freeness is ``unnecessarily restrictive''.  What are we to make of this ``unnecessary restriction''?  I will argue that, in fact, the notion of projective module is a technical contrivance, introduced for aesthetic reasons, as opposed to any of the standard reasons mathematicians use to justify making definitions.  The idea that projective modules could transcend their ``technical'' role is mediated by an observation of Serre, which seemingly only could have been made within the culture of fiber bundles, sheaves, homology etc. that we have described.  

With this context in place, projective modules appear prominently at the beginning of the book, followed quickly by a related notion of injective module.  No mention is made of sources of inspiration for the definition, and projective modules appear as an essentially technical tool and, as I will argue, not even {\em technically} relevant if one looks carefully.  Bourbaki makes one comment about motivation for the notion of projective modules in \cite[p. 145-6]{BourbakiHistoire}.  They assert:
\begin{quote}
	several of the fundamental notions of Homological Algebra (such as those of projective module and that of the functor Tor) were born as a result of the close study of the behaviour of modules over a Dedekind ring relative to the tensor product, a study undertaken by H. Cartan in 1948.
\end{quote}
Bourbaki concludes (I use the French, since the English translation seems to me to have a slightly different feel:
\begin{quote}
	Inversement, on pouvait prévoir que les nouvelles classes de modules introduites de façon naturelle par l'Algèbre homologique...jetteraient une lumière nouvelle sur l'Algèbre commutative. Il se trouve que ce sont surtout les modules projectifs...qui se sont révélés utiles...
\end{quote}
Indeed, by the late 1950s projective modules took on a life of their own, enough to be retrospectively mentioned in Bourbaki's own historiography, becoming a rich theory in its own right.  How did this happen?  

\subsection{Bourbaki's algebra in context}
\label{ss:bourbakialgebra}
\addtoendnotes{\vskip .2em \noindent {\bf Bourbaki's algebra..} \vskip .2em}
Almost immediately in the introduction to Cartan--Eilenberg's text one is confronted by the tensor product for modules and the statement ``It is important to consider the behavior of this construction in relation to the usual concepts of algebra: homomorphisms, submodules, quotient modules, etc.'' \cite[p. vi]{CE}.  I think nothing of such statements now, but how were they received in 1953?  

I will take the above statement at face value and attempt to disentagle where one might have learned about the relevant notions.  After all, the notion of tensor product of abelian groups was defined by Whitney by 1938 \cite{Whitneytensor}, but the more general notion of tensor product of modules over rings was considered only later in Bourbaki's algebra treatise, reprised as \cite{BourbakiTensor}.  When did it become ``important'' to study, in modern parlance, functorial properties and who might have done so?  More widely, how was Bourbaki's approach to algebra received?\footnote{I think it is safe to assert that Bourbaki's algebra (and Lie theory) texts have had the most long-standing influence on modern mathematics; so analyzing the reception of these texts among mathematicians of the era seems useful especially as regards questions of ``heritage'' as discussed in Section~\ref{ss:fiberbundles}.}  Having gone through those exercise, I will then return to the introduction of Cartan--Eilenberg.  I aim to show that the text was written in a way that would have made it easily accessible to only a small group around these authors. 

\subsubsection*{The evolution of tensor products}
We'd like to trace the development of tensor products of modules and its reception in the wider mathematical world.  Bourbaki says only a sentence about the development of this notion tracing it back to Kronecker: they mention Kronecker's work on matrices (what is now called the Kronecker product) \cite[p. 88]{BourbakiHistoire}; in fact Bourbaki characterize Kronecker's work as the introduction of tensor products, in an fashion that was not ``intrinsic'' and without feeling the need to give it a name.  

Certainly for the uses we have in mind, this description seems wholly inadequate.  As an example of the kind of thing that requires comment: in order to formulate something like the ``universal property'' of tensor products, we need to have at hand the notion of functoriality.  This kind of notion seems to have first been considered by Eilenberg and Steenrod, so something like a ``universal'' characterization of tensor products, which appears in \cite{BourbakiTensor} and on which Artin remarks in his 1953 review, is unlikely to have been formulated earlier, i.e., {\em even by Bourbaki}.\endnote{Regarding the evolution of Bourbaki's algebra treatise, we refer the reader to the discussion in \cite[pp. 247-249]{Beaulieu}.  The initial foray into algebra is described as ``timid'' as Bourbaki ``did not want to rewrite van der Waerden's'' text.  It was only the first chapter on ``Algebraic structures'' that was published in 1942.  The chapters on linear and multilinear algebra were conceived later and published after 1945, though many decisions on the treatment had been made by 1942-43.}  Indeed, one of the first public treatments of universal constructions appears in the work of Samuel \cite{SamuelUniversal}, which was submitted in August 1947.  Thus, the treatment in the referenced version (which does include the universal property) likely was a later modification, though the presentation does appear to have stabilized by 1948, as Cartan and Serre describe Bourbaki's theory of tensor products in the \cite[Expose 11]{SemCartan194849}.  

Weibel suggests in \cite{WeibelHHA} that the notion of tensor product appears as well in \cite{ArtinNesbittThrall}.  MathSciNet lists a publication date of 1944, agreeing with the copyright date of the text.  The fact that the book mentions results of Artin--Whaples indicates the book was published after 1943, but other sources, e.g., the chronology in the Artin--Hasse correspondence list other dates of publication \cite[p. 44]{FreiRoquette}.  The chapters of this book in which we will be interested, originated in a course given by Artin in 1941 while he was visiting the University of Michigan (from Indiana University, Bloomington) for a semester.  During this period, it is unclear what, if any, interaction Artin had with topologists like Eilenberg or Wilder; though based on his publications his research focus seemed rather different.  

The Artin--Nesbitt--Thrall monograph speaks of vector spaces over a ring (see \cite[I.3]{ArtinNesbittThrall}); this is in the spirit of the terminology of groups with operators, rather than that of module-theoretic terminology.  This terminology presents a first disjunction with that used by Bourbaki and serves to highlight the variety of different names for this kind of notion.  The notion of Kronecker product of ``spaces'' is then defined using elements which are given by formal juxtapositions of symbols \cite[VI.2]{ArtinNesbittThrall}; the object itself is denoted by $\times$ rather than by means of some other symbol, though notation may have been circumscribed by typographical limitations imposed by the printers.  Some multilinear aspects of the Kronecker product spaces are described.  Unfortunately, the Artin--Nesbitt--Thrall monograph gives no references, so I am unable to determine with certainty whether the authors were familiar with Whitney's work or the work of Bourbaki mentioned above.\endnote{A metaphor that I'd like to propose here is that these two developments of tensor products are examples of something like ``convergent evolution'' in biology.  We can think of two species of tensor products $\tensor_{ANT}$ and $\tensor_B$ due to Artin--Nesbitt--Thrall and Bourbaki.  These two species exhibit similar morphological features, say interpreted in terms of the kinds of results that are established about them, but the similarity of those features is likely not due to any direct link between them: neither evolved from the other because of the psychological separation imposed by WWII and the lack of available channels of communication.  Bourbaki's theory ends up considerably more developed (it is successful), whereas $\tensor_{ANT}$ is its own evolutionary endpoint in some sense (even the notation disappears).}  However, Artin's stay at Michigan was temporary and dissemination of texts during the war was difficult.  The eventual goal is to construct Kronecker products of rings.   Note that Bourbaki, like Artin--Nesbitt--Thrall, mentions tensor products of algebras, and credits this notion to B. Pierce \cite[p. 150]{BourbakiHistoire}.

In any event, Artin's experience with Kronecker products places him in a unique position to appreciate the significance of the tensor product construction, which he later does rather explicitly when he highlights the universal property in his 1953 review of Bourbaki's algebra text \cite[pp. 476-7]{ArtinBourbakireview}.  I am therefore led to conclude that this discussion of Kronecker products of ``spaces'' is a parallel development, not overlapping with tensor products as used later in topology, and consequently playing no role in how the notion was used later; only retrospectively could it be viewed as a precursor in any sense.  

Postwar, Princeton (and the IAS) had become the center for ``mainstream'' mathematics \cite{Aspray} and Artin moved there in the Fall of 1946.  Undoubtedly, the sphere of Artin's interaction and influence grew here, and his students in the period up to his review included John Tate and Serge Lang.  Rota characterizes his personality here thus \cite[p. 14]{Rota}:
\begin{quote}
	A great many mathematicians in Princeton, too awed or too weak to form opinions of their own, came to rely on Emil Artin's pronouncements like hermeneuts on the mutterings of the Sybil at Delphi. He would sit at teatime in one of the old leather chairs (``his'' chair) in Fine Hall's common room and deliver his opinions with the abrupt definitiveness of Wittgenstein's or Karl Kraus's aphorisms.
\end{quote} 

Artin had particularly strong views on algebra, and from Rota's point of view Artin's views were consonant with those of Gauss and Dirichlet \cite[p. 14]{Rota}: ``a piece of mathematics was the more highly thought of, the closer it came to Germanic number theory.''  In this regard, the Cartan--Eilenberg text appealed to Artin's sensibilities, as aspects of the cohomology of groups relevant to the Artin--Tate formulation of class field theory are developed there \cite[Chapter XII]{CE}.  Artin's review of Bourbaki's algebra text is also revealing in several respects \cite{ArtinBourbakireview}.  

Artin begins with a striking general observation of the nature of exposition of mathematical results.
\begin{quote}
	We all believe that mathematics is an art. The author of a book, the lecturer in a classroom tries to convey the structural beauty of mathematics to his readers, to his listeners. In this attempt he must always fail.  Mathematics is logical to be sure; each conclusion is drawn from previously derived statements. Yet the whole of it, the real piece of art, is not linear; worse than that its perception should be instantaneous.  We all have experienced on some rare occasions the feeling of elation in realizing that we have enabled our listeners to see at a moment's glance the whole architecture and all its ramifications. How can this be achieved? Clinging stubbornly to the logical sequence inhibits the visualization of the whole, and yet this logical structure must predominate or chaos would result. Bourbaki is quite aware of this dilemma. The fact that his work is subdivided into books, the fact that exercises are given which utilize more advanced parts of the theory show this awareness. However I feel that in some instances the subdivision into books is not enough.
\end{quote} 
As evidence of the ``not enough'' of the final line, Artin describes the Bourbaki treatment of Galois theory, evidently one of the fields closest to Artin's heart, lamenting Bourbaki's avoidance of arithmetic ideas:
\begin{quote}
	...a heavy price has to be paid for the fact that one is not permitted to use number theory...almost no example over the rational field can be given, since it is next to impossible to show the irreducibility of a polynomial without some arithmetic.
\end{quote}
Once again, laid bare is the fact that Bourbaki has made a choice about what results are primary and what results are secondary, and that these choices are not universally agreed upon.  Artin nevertheless closes with his approval:
\begin{quote}
	In concluding I wish again to emphasize the complete success of the work. The presentation is abstract, mercilessly abstract. But the reader who can overcome the initial difficulties will be richly rewarded for his efforts by deeper insights and fuller understanding.
\end{quote}
Thus, in the end, Bourbaki is forgiven for his sins, and absolution granted.  Nevertheless, Artin's criticisms highlights something about the reading of Bourbaki that I would like to explore, recalling some of Zariski's criticism of the strict logical progression of Weil's {\em Foundations of algebraic geometry}.  Even though Bourbaki's texts were nominally textbooks, how did one learn mathematics from them?  

\subsubsection*{Outside Bourbaki's algebra}
The reception of the ideas of Bourbaki's algebra, even by those in Artin's circle in the early 1950s was not uniform, and the reception of Bourbaki texts outside the scope of algebra was even less uniform.  We probe these two directions separately.  We begin by discussing Chevalley's Fundamental concepts of algebra \cite{ChevalleyAlgebra}, released in June 1956, as a proxy for Bourbaki's algebra.

Arthur Mattuck received his Ph.D. from Princeton in 1954 under Artin's direction.  Mattuck seems poised at the cusp: by this point in time Bourbaki's ideas are relatively widely available and well-exposed, certainly in Artin's circle.  Mattuck reviewed Chevalley's text for the Bulletin; his review suggests an attraction to the intuitive nature of pre-Bourbaki treatments of algebra, but simultaneously aware of the usefulness of the new approach, which he  nevertheless characterizes as ``austere''.\footnote{For the sake of comparison, Rota describes Mattuck's lectures \cite[p. 16]{Rota} as ``an exercise in high motivation.''}

Beyond the reasons suggested above, we analyze the reception of Chevalley's text because it seems to be one of the first sources that mentions projective and injective modules, via the lifting property, but, in support of the theme we will pursue in Section~\ref{ss:axiomaticprojectivity} defers them to an exercise \cite[p. 130 Exercise 17]{ChevalleyAlgebra}.\footnote{According to Mattuck \cite{MattuckReview}:
	\begin{quote}
		There are a number of exercises after each chapter—well over a hundred in all. Almost none of then are routine, nor are there many in the nature of specific examples; most of them are not easy. Many provide significant extensions and applications of the theory, and are the equivalent of many additional pages of text. Such, for instance, would be the exercises on the derived groups of a group, projective and injective modules, projective limits, quadratic forms and Clifford algebras, and representations of Lie algebras.
\end{quote}}

Mattuck's review \cite{MattuckReview} begins with the following characterization of Chevalley's project; the austerity motif seems to appear repeatedly in description of Bourbaki's algebra:\footnote{The reference in the quote below is to the first paragraph of Weil's review of Chevalley's ``Introduction to the theory of algebraic functions of one variable'' for the Bulletin \cite{WeilChevReview}.}
\begin{quote}
	Readers...who agreed at the time with André Weil's dictum ``algebraic austerity can go no further'' may decide that a counterexample has been produced.
\end{quote}
Mattuck comments in detail on his view of the perceived level of abstraction in Chevalley's treatment and whether it is appropriate for the subject matter at hand, but aesthetic judgments appear in this evaluation.
\begin{quote}
	In considering the treatment the book presents of its subject, one must recognize that it is extremely abstract, and the level of abstraction each man likes in his mathematics is as personal a taste as the amount of perfume he likes on his wife. When linear transformations made their comeback over matrices, it was easy to point to the shortened proofs and to the gain in clarity resulting from the triumph of abstraction over algorithm. An intense sans-culottism has since made the subject perhaps a bit top-heavy in concepts (after all, we still do have a multiplication table). Each reader will have to decide for instance whether Chevalley's seven page intrinsic treatment of the Pfaffian here is the height of beauty and elegance, or of absurdity, and if the former, whether the associated aroma is that of ripe bananas or of freshly-roasted coffee.
\end{quote}
Evidently, I wish to view the reception of Chevalley's treatment as indicative of the reception of Bourbaki's ideas more generally, but it is perhaps difficult to the draw the line so that I am not simply extracting features of Chevalley's particular presentational idiosyncracies (writing as a group certainly served to rounded out such individual tendencies that leaned extreme in Bourbaki's presentations).\endnote{Since our analysis here turns on ideas of Chevalley's mathematical style, it is useful to keep in mind that Chevalley had himself written on ``mathematical styles'' in 1935 \cite{Chevalley1935-CHEVDS}.  Two points from this essay seem quite relevant to me.  First, Chevalley asserts that authorial style not withstanding, different eras have recognizable styles that are influenced by strong mathematical personalities:\footnote{Le style math{\'e}matique, tout comme le style litt\'eraire, ne va pas sans subir d'une \'epoque \`a l'autre d'importantes fluctuations.  Sans doute, chaque auteur poss\`ede-t-il un style proper; mais on peut apercevoir \`a chaque \'epoque une tendance g\'en\'erale assez bien reconnaissable.  Ce style subit, de temps \`a autre, sous l'influence de personnalit\'es math\'ematiques puissantes, des r\'evolutions qui infl\'echissent l'\'ecriture, et donc la pens\'ee, pour les p\'eriodes qui suivent.}\begin{quote}
		Mathematical style, like literary style, cannot but undergo significant fluctuations from one period to another. No doubt, each author has his own style; but one can perceive in each period a general tendency that is quite recognizable. This style undergoes, from time to time, under the influence of powerful mathematical personalities, revolutions that influence writing, and therefore thought, for the periods that follow.
	\end{quote}
Chevalley goes on to illustrate an example of the matheamtical style of an epoch by discussion Weierstrass and what has been called by other Bourbaki members ``epsilontics''.  As regards the rise of ``epsilontics'' he explains that this new style arose from criticisms of lack of rigor around arguing with infinitesimals.  But Chevalley's essay then turns to the existence of a new style, arising again from criticisms of the old style; I find his stated reasons fascinating \cite[pp. 379-80]{Chevalley1935-CHEVDS}: \footnote{``On peut donc dire que les définitions constructives de l'analyse, si elles ont les premières permis les raisonnements rigoureux, ont eu souvent l'effet de cacher profondément la nature de ce q'uelles  cherchaient à définir ou de confondre indûment des domaines mathématiques en réalité distincts les uns des autres. De là résultent les complications inutiles qui se rencontrent dans beaucoup de démonstrations classiques, du fait de l'emploi de méthodes n'ayant rien à voir avec le résultat escompté, on pourrait dire : de méthodes n'admettant pas le même groupe de transformations que le résultat.}
\begin{quote}
	We can therefore say that the constructive definitions of analysis, if they were the first to permit rigorous reasoning, have often had the effect of deeply hiding the nature of what they sought to define or of unduly confusing mathematical domains that were in reality distinct from one another. From there result the unnecessary complications that are encountered in many classical demonstrations, due to the use of methods having nothing to do with the expected result, one could say: of methods not admitting the same group of transformations as the result.
	\end{quote}
Thus, Chevalley is perfectly aware that particular choices of abstract definitions can obscure meaning.  In light of what we have just read, is it thus the case that this is just an early view of Chevalley that matures into the ``formalist'' perspective we describe Bourbaki as having?  I think that is not the case: Chevalley viewed the level of abstraction he used as illuminating.  Moreover, the later Bourbaki style was in fact representative of the style of an era in Chevalley's sense, undoubtedly under the influence of a number of strong mathematical personalities.}
Nevertheless, Mattuck's comments on Chevalley's tone seem quite relevant as they directly reflect the reception of Bourbaki's ideas among the ``older generation'' and we have heard their echo in the reception of other Bourbaki-aligned presentations.
\begin{quote}
	The voice that we hear resounding is that of an Old Testament prophet, but the mental attitude is more like a tenth-grade Latin teacher's, reeking with the old theory of formal disciplines. Thinking rigorously demands first of all a firm grasp of the concepts, otherwise the sort of proof-following which passes for thinking is only a very sophisticated version of computation-checking. It is the difference between a rat running physically through a maze and a man running his pencil through one of the Sunday supplement mazes: one has the over-all understanding, the other does not. It is downright unfair for an older generation which learned these ideas in an intuitive fashion in which they were well-adapted for thought to foist off on a younger one, in the name of rigorous thinking and without any further explanation, such a construction as this one...
\end{quote}
Nevertheless, Mattuck's comments are tempered by momentary cracks of admiration and assent: ``This basic and essential usefulness of the book should be kept in mind as overshadowing any critical remarks made below.''  Mattuck also strongly praises the ``unity'' of Chevalley's algebra.
\begin{quote}
	One cannot emphasize too strongly the beautiful unity of the book: the pruning has been severe, perhaps, but at least one will not forget the essentials that have been left.
\end{quote}
Simultaneously, Mattuck seems to highlight the tension between unity, in this case revealed to us by presentational choices, and intuition, for which there is no easy explanation.

Rather than leave Mattuck's analysis on its own, we can also consult another review of Chevalley's book, this time by Irving Kaplansky for Math Reviews, which observes that Chevalley's treatment of multilinear algebra is quite reflective of the Bourbaki presentation.
\begin{quote}
	A generation of algebraists grew up for whom ``modern algebra'' meant Van der Waerden's book, or possibly one of several similar later texts. Time has passed and (happily) mathematics has not stood still. In particular algebraic topology has exhibited an insatiable appetite for algebraic gadgets. In response, modern algebra has changed.
	
	What distinguishes the new modern algebra from the old? The latter emphasized groups, rings, and homomorphisms as the basic concepts. Modules, more or less sitting astride groups and rings, were prominent, though perhaps not sufficiently prominent. But at least two things, now clearly of central importance, were completely missing: the tensor product of modules, and the generalization of every object to a graded object.
	
	Chevalley's book is timely and it will be widely studied; the meaty exercises will invite a diligent reader to educate himself some more. Teachers may find it ``futile to disguise the austerity'' (last sentence of the preface). It goes without saying that large sections are similar to Bourbaki's ``Algèbre multilinéaire'' [Actualités Sci. Ind., no. 1044, Hermann, Paris, 1948; MR0026989].
\end{quote}
Once again, we see the austerity motif.  

Reading these criticisms first, the start of Chevalley's rather circumspect preface might come as a shock:
\begin{quote}
	Algebra is not only a part of mathematics; it also plays within mathematics the role which mathematics itself played for a long time with respect to physics. What does the algebraist have to offer to other mathematicians?  Occasionally, the solution of a specific problem; but mostly a language in which to express mathematical facts and a variety of patterns of reasoning, put in a standard form. Algebra is not an end in itself; it has to listen to outside demands issued from various parts of mathematics. This situation is of great benefit to algebra; for, a science, or a part of a science, which exists in view of its own problems only is always in danger of falling into a peaceful slumber and from there into a quiet death. But, in order to take full advantage of this state of affairs, the algebraist must have sensitive ears and the ability to derive profit from what he perceives is going on outside his own domain, Mathematics is changing constantly, and algebra must reflect these changes if it wants to stay alive. This explains the fact
	that algebra is one of the most rapidly changing parts of mathematics: it is sensitive not only to what happens inside its own boundaries, but also to the trends which originate in all other branches of mathematics.
\end{quote}
The final statement here is prescient, especially in light of the ``cohomological turn'' concurrent with the release of this text.  While the initial impetus was a development of algebraic notions in the service of topology, Bourbaki's treatment of algebra and the focus on constructions such as tensor product and their properties, allowed the flow of information to turn in the opposite direction: this breed of algebraists was prepared to allow cohomology to ``clarify'' algebra.  This began with Eilenberg and Mac Lane's cohomology of groups, buoyed by the then recent {\em mainstream} applications and developments of this cohomology theory to class field theory in the work of Artin--Tate, the cohomological formulation of the theory of separable algebras studied by Hochschild, and the incorporation of cohomological ideas in formulation of Lie theory.  Even if cohomological formulations did not prove new results, they could be used to {\em unify} and {\em streamline} or incorporate old results.  The algebraic language became the one people needed to speak.  In retrospect, the transition in styles from 1945 to 1956 seems practically revolutionary.  

That people would adopt this language and terminology was certainly not a foregone conclusion, as one can see by comparing with the reception of other Bourbaki texts.  Halmos reviewed Bourbaki's ``Integration'', which was itself published in 1952, for the Bulletin \cite{Halmos}.  He finds several aspects of the treatment artificial
\begin{quote}
	Owing, no doubt, to the authors' predilection for using as definiens what for most mathematicians is the definiendum, there are many spots at which the treatment appears artificial.
\end{quote}
He then goes on to detail several examples of this artificiality, explaining how these are forced by the choices the authors have made, for example writing, ``the effect of the definition is to help perpetuate the myth that the measurability of a function can only be defined by reference to a measure—a myth quite as unfounded as the (fortunately moribund) myth that the continuity of a function can only be defined by reference to a metric.''  Halmos, complaining about the unsuitability of Bourbaki's approach to probability theory closes with:
\begin{quote}
	My conclusion on the evidence so far at hand is that the authors
	have performed a tremendous tour de force; I am inclined to doubt
	whether their point of view will have a lasting influence.
\end{quote}
In the review proper, he carefully details a host of problems he sees with the choice of presentation of material.  Indeed, Halmos' predictions seem to have been borne out; one modern review of Bourbaki's Integration by R. Schilling \cite{Schilling} writes:
\begin{quote}
	The Bourbaki volumes on integration are nowadays mainly of historic value, which is partly due to the (dogmatic) misconception which limited the theory to a locally compact setting and the, especially for probabilists, cumbersome notion of measurability.
\end{quote}

This dual treatment (Halmos and Artin) of the Bourbaki project did not go unnoticed by Weil who remarks on it \cite[p. 340]{CartanWeil}:\footnote{Le papier de Artin sur notre Alg\`ebre nous d\'edommage largement du torchon de Halmos.  Ne conviendrait-il pas de lui envoyer d\'esormais nos r\'edactions d'alg\`ebre suffisamment avanc\'ees, afin d'avoir une critique qui nous serait utile?}
\begin{quote}
	Artin's paper on our Algebra largely compensates us for Halmos's rag.  Would it not be appropriate to send him from now on our sufficiently advanced algebra essays, in order to have a critique that would be useful to us?
\end{quote}
And following up on the reviews of Bourbaki's Integration, the notes in the Cartan--Weil correspondence \cite[p. 622]{CartanWeil} suggest that Bourbaki's ``Integration'' is still useful in modern $p$h-adic arithmetic settings.


\subsubsection*{Receiving the introduction of Cartan--Eilenberg}
Having analyzed the reception of Bourbaki's algebra texts, we return to analyze how the Cartan--Eilenberg monograph may have been received.  First, there were practical differences between the notation and terminology of this text and notation and terminology that were more ``widely'' in use.  For example, the notion of kernel of a homomorphism was widely used, but consulting papers in group theory from this era and before (e.g., those papers we will discuss in Section~\ref{ss:axiomaticprojectivity} as regards the evolution of the notion of ``injective'' module), one sees that when considering a single homomorphism it was standard to give the kernel a special name.  Uniform notation for kernels only becomes useful when one considers many homomorphisms simultaneously, which of course was the case in homology.  The notation $\operatorname{Ker}$ does not seem to be widely used even by 1952.  Even in \cite{EilenbergSteenrod} which did consider many homomorphisms simultaneously, one frequently sees $\operatorname{Kernel}$ spelled out in its entirety.  Cartan and Eilenberg seems to be the first place where this notation is systematically used.  Second, Cartan and Eilenberg systematically use diagrams and what might be called diagrammatic arguments.\endnote{It is worth saying a word about what a ``diagrammatic argument'' looks like.  A nice example of this kind of thing can be seen in the $5$-lemma, which appears as \cite[Proposition 1.1]{CE}.  The proof of this result as given in Cartan and Eilenberg is perfectly fine, but usually one performs this by means of a ``diagram chase''.  If you haven't seen this before, it makes sense to watch a video of its performance, of which you can find many.}  While diagrams themselves may have been psychologically useful, diagrammatic arguments have a mechanical, performative nature that in my experience is better communicated by demonstration than by reading.  Keeping all that in mind, I think it is fair to say that the nature of the introductory example must have seemed fairly inscrutable to those outside a small circle around Cartan and Eilenberg.  

Mac Lane claims that \v{C}ech ``first introduced (but did not name) the torsion product'' \cite[p. 172]{MacLanehomology} making reference to \cite[p. 34]{CechUnivCoeff}.  Indeed, does \v{C}ech consider the $n$-th homology group with coefficients in some abelian group $A$.  He aims to describe homology with coefficients with $A$ of some ``infinite complex'' in terms of previously defined invariants, e.g., Betti numbers.  His formulation of the universal coefficient theorem is contained in two results: \cite[Th\'eor\`eme II and III]{CechUnivCoeff}.  He identifies two contributions $\mathbf{B}_1^n(A)$ and $\mathbf{B}_2^n(A)$ to the $n$-th homology with coefficients in $A$.  The first result states that the first contribution is determined entirely by $A$ and the $n$-th Betti group, while the second result asserts that the second contribution is completely determined by $A$ and the torsion in the $n-1$st homology group.  These descriptions are all given in terms of equations involving generators and relations arising from the combinatorial description of the complex.  

As we described earlier, the tensor product of abelian groups had not yet been invented (this construction predates even \cite{Whitneytensor} who himself made no mention of anything like torsion product).  Only someone who is well-versed in the formal aspects of the tensor product would be able to extract the ``tensorial'' description of the torsion product from this presentation.  As a consequence, it seems better to say that \v{C}ech highlighted a  contribution to a homology group, which was retroactively realized to be a concrete instantiation of the functor described by Eilenberg and Mac Lane.  This is not a distinction without a difference: from \v{C}ech's point of view, his problem of describing the homology group in question was solved by the generators and relations picture: he effectively gave a formula for the determination of the groups from the combinatorial description of a geometry.  

Bourbaki's analysis of tensor products is reprised in the Cartan seminar in an expose attributed jointly to Cartan and Serre \cite[Expos\'e 11]{SemCartan194849}.  Bourbaki analyzes behavior of homomorphisms under tensor products and formulate various ``exactness'' statements.  The current version of Bourbaki considers this section ``Proprietes..relatives aux suites exacts'' but the notion of exact sequence was not widely used yet.\todo{Check!}  One point to which Cartan and Serre draw attention, which is not mentioned in Bourbaki is the behavior of tensor product with respect to homomorphisms that are {\em projectors} (while the latter are defined in Bourbaki's algebra treatise, no special mention is made of projectors).  At this stage, it makes sense to measure the failure of left exactness of tensor product, but no one appears to explicitly do so until Cartan and Eilenberg's book.  

Imagine, then, opening Cartan and Eilenberg and immediately being confronted with free resolutions, higher torsion products etc., all of which appeared to be {\em new} notions at the time.  The introduction continues to layer definition on top of definition.  After having been confronted with the notion of free resolution, Cartan and Eilenberg immediately add yet another layer of abstraction:
\begin{quote}
	The condition that...free in the definition...is unnecessarily restrictive. It suffices [to]...be projective...
\end{quote}
What evidence do they give for the pronouncement that this restriction is {\em unnecessary}?  Perusing Cartan--Eilenberg, we will argue in the next section that the answer is, apparently, none.  

After the introduction, one has the sense that a vast new landscape has been unveiled: infinite sequences of new homology groups.  What do these groups mean...what do the new groups measure?  What tether to older more concrete notions do we have?  At this stage, the key notion is unification of constructions and proofs.  

\subsection{On the axiomatic necessity of projectivity}
\label{ss:axiomaticprojectivity}
\addtoendnotes{\vskip .2em\noindent {\bf On the axiomatic..} \vskip .2em}
In the introduction to Cartan--Eilenberg, projective modules are introduced via the lifting property described above.  In Chevalley's \cite{ChevalleyAlgebra}, projective modules are mentioned, once again via the lifting property, but in this case they are relegated to an exercise.  The particular style of these mentions leads one to speculate that projective modules were not conceived initially as {\em central} actors, but relegated to a purely technical role with no ``practical'' mathematical value.  This section and the next will, undoubtedly, be the most purely mathematical sections of the text.

To make this case about the ``status'' of projective modules, I will employ several devices.  First, I will show that in early formulations of homological algebra, projective modules make no explicit entrance whatsoever: the results which the theory aims to establish can be established by consideration the considerably more limited and intuitively simpler notion of {\em free} modules.  Second, while free modules serve as examples of projective modules, I will argue that non-free projective modules appear as  ``curiosities'' early-on: non-free projective modules are never highlighted explicitly, only implicitly, if one extrapolates from the text.  Finally, I will suggest that bearing the above in mind, the introduction of projectivity is guided by axiomatic ``aesthetic'' concerns.

\subsubsection*{From freeness to projectivity}
Here we counterpoise Eilenberg's treatment of homology of groups for the Cartan seminar in 1951 \cite[Expose I-II]{SemCartan195051} with earlier discussions of group cohomology in the work of Eilenberg--Mac Lane and later discussion of Cartan--Eilenberg to track the evolution of the theory.  In particular, our discuss here centers on the contention from Cartan--Eilenberg's introduction that to be ``free...is unnecessarily restrictive'' \cite[p. viii]{CE}.

The first expose, dated November 13, 1950 sets out Eilenberg's aim: he wishes to give a systematic development of the homology and cohomology of discrete groups, following the axiomatic method.  To understand his goal further, we must investigate precisely what he means by ``axiomatic'' as compared to previous treatments of (co)homology of groups (we have already observed that, for example, Mayer had given axiomatics for homology), and thus we should review what previous treatments of the cohomology of groups looked like, at least schematically.  

As we have sketched it, up to this point, homology of some kind of structure has been presented in terms of one explicit complex.  In more detail, the original Eilenberg--Mac Lane treatment of (co)homology of groups \cite{EilenbergMacLaneCollected}, say even as summarized later in \cite{EilenbergBull} is given by (co)homology of what we now call the bar resolution.  The cohomology of associative algebras as developed by Hochschild was given as cohomology of the Hochschild complex \cite{Hochschild1}.  The cohomology of Lie algebras was given in \cite{ChevalleyEilenberg} also in terms of cohomology of an explicit cochain complex.  

I would like to frame this in terms of a dichotomy emphasized in the Cartan--Weil correspondence: the tension between ``existence and uniqueness''.  To this point, there is no question of {\em existence} of any of the cohomology theories mentioned above: they are all given by cohomology of some explicit complex.  We thus find ourselves in a situation entirely unlike the situation for (co)homology of manifolds where there were multiple putative methods for building up workable (co)homology groups (see Paragraph~\ref{par:manyhomologies}), and thus the {\em uniqueness} question does not even arise for each of the theories above.  
 
Nevertheless, the starting point of Eilenberg's approach is that there are {\em in principle} many approaches to group homology.  For the reader's convenience, we summarize his treatment in modern terminology.  Start with a discrete group $G$ and a (left) $\Z[G]$-module $A$ (he calls $A$ a $G$-group).  Group homology, is then any procedure that:
\begin{itemize}[noitemsep,topsep=1pt]
	\item[(1)] attaches to each such $A$ a sequence $H_q(A)$ of (abelian) homology groups, together with 
	\item[(2)] for each homomorphism $A \to A'$ of $\Z[G]$-modules an induced homomorphism $f_*:H_q(A) \to H_q(A')$, and 
	\item[(3)] given a short exact sequence of $\Z[G]$-modules $0 \to A' \to A \to A'' \to 0$, connecting homomorphism $\partial: H_q(A'')\to H_{q-1}(A')$.
\end{itemize}  The main claim that Eilenberg makes is the assertion that if such a group homology theory furthermore satisfies suitable naturality, and normalization axioms, then the homology procedure is unique.  Thus, Eilenberg has, with this presentation, shifted the discussion: having described homology axiomatically, and answering the uniqueness question, the existence question arises anew.  One wonders if a conservative reader would have supposed that a new problem has been artificially created.

For Eilenberg, this was evidently not an issue, and he quickly describes a procedure for producing many different potentially inequivalent homology theories for a group.  On the very first page of this note, he makes definitions of a {\em free module} and an {\em injective module}.  The ``new'' cohomology groups arise via different {\em free resolutions}.  By way of passing, the definition of injective used here is identical to that in Cartan--Eilenberg, i.e., the modern use, and this appears to be the first instance of this word in the literature. 


Having created all these potentially new versions of (co)homology by means of different resolutions, we are forced to consider new questions.  Given two such approaches, how can we even compare them?  It's worth pointing out that Eilenberg's axiomatization only makes sense when one considers {\em arbitrary $G$-modules}, as opposed to considering ``classical'' group cohomology using only ``trivial'' coefficients.  Thus, axiomatization only seems reasonable when we move to the ``variable coefficient'' setting.  This potential variety of new approaches, now unified by the axiomatic approach, is once again viewed as a positive feature of the theory.

Cartan--Eilenberg indicate in the introduction to \cite{CE} a further point that they refer to as unifying, in the context of the three theories alluded to in the discussion above: 
\begin{quote}
	In addition an interplay takes place among the three specializations; each enriches the other two.  The unified theory also enjoys a broader sweep. 
\end{quote}
Here they allude to applications of the relationship between cohomology of commutative rings and Hilbert's theory of syzygies.  

What Cartan and Eilenberg {\em fail} to mention is a psychological trade-off.  Instead of having to consider a {\em single} complex to compute (co)homology in a given setting, one is now confronted with a multitude of different complexes.  Of course, the uniqueness theorem tells us that the multitude of new computational approaches to (co)homology yield the same answers, but there is a cost: the unification comes at the expense of development of considerable (at the time) overhead: they spend the bulk of Chapters III-V developing tools for comparing different complexes computing (co)homology.

In \cite[Expose I-II]{SemCartan195051}, we also do not yet see the full sweep of different resolutions, nor even the terminology of a resolution (Eilenberg only speaks of acyclic complexes).  What is furthermore conspicuous, by comparison with later treatments is a lack of symmetry: injective resolutions are {\em not considered} and the notion of {\em projectivity} makes no appearance.  

Returning to the fact that Eilenberg computes cohomology using {\em only} free resolutions, at least at the beginning freeness was definitively {\em not} viewed as ``unnecessarily restrictive'': it was sufficient for establishing uniqueness of the relevant (co)homology theories.  Furthermore,  while Eilenberg mentions that every module can be embedded in an injective module, he does not provide a construction (in contrast to the fact, used almost immediately in the Expose, that every module admits a surjection from a free module).

Given the way that results from Cartan--Eilenberg trickled out (150 pages by May 1951 as Weil mentions in his letter), one wonders precisely when projectivity was introduced in this period.  We will analyze the parallel story of injectivity momentarily.  Having just argued that free resolutions are sufficient to establish the uniqueness theorems, the {\em only} motivation for introducing projective modules and projective resolutions into the story is axiomatic minimality: they provide the answer to the question what is the minimal structure necessary to write the proof of uniqueness. Slightly more broadly, if one is interested in simply {\em defining} satellites of a given functor, and then computing the resulting (co)homologies, then free resolutions are entirely sufficient.  

\subsubsection*{The red-herring of rational reconstruction: invertible modules}
The standard modern treatments of projective modules quickly turn to the Dedekind--Weber theory of invertible ideals for examples.  Since the work of Dedekind--Weber was clearly a motivating force in approaches to modern algebra, and as such a theory with which one expects Bourbaki-aligned mathematicians to be familiar, it is tempting to frame the theory of invertible ideals as a source of {\em motivation} for the theory of projective modules.  Thinking this true, I asked Serre for confirmation, only to be rebuffed with the response \cite{Serreprivate}:
\begin{quote}
	Another thing : you seem to believe that the correspondence ``invertible line bundle'' $\leftrightarrow$ ``ideal classes'' (under suitable conditions on the ring) was standard before the 50s. I don't think so: it was clarified in the period 1952-1955.
\end{quote}
I would like to explore this here further. 

First, I take up the fact that not a single ``non-trivial'' example of a projective module is given throughout the Cartan--Eilenberg text.  Here, we need to qualify the characterization ``non-trivial'' and in doing so we will return to the discussion of projective modules in Bourbaki's history. Evidently free modules are projective as follows immediately from the lifting property.  However, does the Cartan--Eilenberg text give examples of non-free projective modules?  One is tempted to say that the answer is ``no'' as there is no {\em explicit} mention made at any point in the book of such an example.  There is, however, {\em implicit} discussion of non-free projective modules in the treatment of modules over hereditary rings from \cite[I.5]{CE}.  

By definition a (not necessarily commutative) ring $R$ is called {\em left hereditary} if every left ideal is a projective module; this ring-theoretic notion appears to have also been defined by Cartan and Eilenberg.  Cartan and Eilenberg make reference to some work of Kaplansky \cite{KaplanskyDedekind} which studies modules over a Dedekind domain or a valuation ring, ``modernizing'' and extending earlier work of Steinitz \cite{SteinitzDedekind}.\footnote{From the standpoint of classical algebra, projectivity is a notion that has a somewhat subtle relationship with finite generation.  Indeed, the rational numbers $\Q$ is a uniquely divisible abelian group and thus provides an example of an (infinitely generated) torsion-free $\Z$-module that fails to be free.  Thus, while many authors had isolated {\em torsion-free} groups as interesting, since many algebraists of the day, e.g., Kaplansky, were concerned with not-necessarily finitely generated abelian groups, given the tenuous relationship between projectivity and finite generation, it strikes me as unlikely that an algebraically-inspired mathematician of the age would have isolated projectivity as a useful concept.}  Cartan and Eilenberg state \cite[Theorem I.5.3]{CE} that if $R$ is a left-hereditary ring, then every submodule of a free module is isomorphic to a finite direct sum of modules, each of which is isomorphic to a left ideal of $R$.

Recall that the structure theorem for finitely generated modules over a principal ideal domain implies that finitely generated torsion-free modules are free (and thus projective); this was known in some form, e.g., via the theory of Smith normal form from the 1860s or, as we described it in Noether's conception, the ``elementary divisor theory''.  It took quite some time to consider more general rings.  However, the theory of finitely-generated modules over more Dedekind domains was considered in detail by Steinitz \cite{SteinitzDedekind}.  From a matrix-theoretic point of view, Steinitz analyzed finitely generated torsion-free modules over Dedekind domains (this is the formulation given by Kaplansky in \cite{KaplanskyDedekind}) observing that a torsion free module is a finite sum of invertible ideals.\footnote{Continuing the discussion of the preceding footnote: Kaplansky is interested in {\em not necessarily finitely generated} modules over Dedekind domains.  Once again, the isolation of projectivity seems somewhat unnatural from the point of view of general modules.}  Steinitz's result was considered classic, and also made its way into Bourbaki's algebra treatment (and so should be taken as well-known to Bourbaki members).  

Now, the notion of invertible ideal was studied for a long time in arithmetic situations.  However, it appears to only have been realized {\em after} making the definition of projectivity that invertible ideals gave rise to non-free projective modules.  Indeed, granted the analysis of torsion-free modules in Kaplansky's 1952 paper, it would appear that projectivity was still not widely used by this point.  That being said, it was widely known that the ideal class group of a Dedekind domain could be non-trivial, i.e., there were, in general, invertible fractional ideals that were not principal.  In that case, Cartan and Eilenberg write ``it will be seen later'' that hereditary integral domains (now in the commutative case) are precisely the Dedekind domains.

Further evidence for propagation in this direction arises from the notes of Kodaira and Spencer \cite{KodairaSpencer,KodairaSpencerII} discussing line bundles on complex varieties and divisor class groups on algebraic varieties in terms of sheaf cohomology.  These results can be viewed as published versions of some of the ideas exposed in Weil's letters to Cartan on the same topics. Thus, only {\em after} one formulates the general definition of projectivity, is there a search for examples, and it seems unlikely that considerations of pure algebra would have inspired the definition.

Moreover, this background use of projectivity, and the fact that the only non-trivial examples of such even implicitly given in the Cartan--Eilenberg text are invertible ideals contributes once again to the feeling that projective modules are simply technical devices.  Moreover, to a ring-theorist of the day, most of the non-trivial examples of projective modules given were simply torsion-free.  One can easily imagine a skeptic wondering whether this new notion was just some esoteric notion amongst the zoo of special types of modules of the day.

\subsubsection*{Injective modules, a counterpoint}
As a complement  to our discussion of projective modules, we now survey some episodes in the development of the notion of injective modules, taking a tour through some aspects of finite group and module theory in the 1930s and 1940s.  Weibel's history of homological algebra discusses injectivity \cite[\S 3.2 p. 817]{WeibelHHA} and makes the flat assertion that ``injective $R$-modules were introduced and studied by R. Baer''.  Moreover, Weibel suggests a linear, cumulative narrative surrounding the development of this notion, fitting nicely into a point of view about mathematics suggested by H. Hankel \cite[p. 15]{Moritz}:
\begin{quote}
	In Mathematics alone each generation builds a new story to the old structure.
\end{quote}
I refer the reader to \cite[\S 3 p. 263]{CroweMisconceptions} for further criticism about this point of view, but here I aim to argue that the mathematical record suggests a much more complicated narrative.

The name {\em injective module} appears, to the best of my knowledge, for the first time in Eilenberg's treatment of the cohomology of groups for the Cartan seminar \cite[Expose I]{SemCartan195051}.   While Eilenberg phrases the property in the context of groups with operators, in module-theoretic language he gives precisely the characterization of injectivity in terms of the lifting property, i.e., if $I$ is a left $R$-module, then $I$ is injective if and only if given any $R$-module map $A \to I$ and an injective map $A \hookrightarrow B$, there is an induced map $B \to I$ extending the map $A \to I$; I will refer to this as the ``lifting'' characterization of injectivity.  

The lifting characterization of injectivity has some immediate consequences already in the case of abelian groups, i.e., $\Z$-modules.  Indeed, if $I$ is any injective abelian group and $a \in I$ is an element, then we can consider the homomorphism $\Z \to I$ sending $1$ to $a$.  The multiplication by $n$ map $\Z \to \Z$ is injective, and the lifting property then guarantees that we can find an element $a' \in I$ such that $na' = a$.  In other words, the multiplication by $n$ map is surjective.  To say this slightly differently, we have deduced that injective abelian groups are {\em divisible} abelian groups.  A similar argument can be made for injective modules over an integral domain, but maintaining the distinction between injectivity and divisibility will be important in what follows.

To develop the story, we need to turn back the clock a bit.  In making the ``modernist'' transition away from matrices to abstract structures, one apparent source of problems in the 1930s group and module theory surrounds transitioning from results about finitely generated modules to various not-necessarily finitely generated situations.  For example, the structure theorem for finite abelian groups or finitely generated abelian groups was well-known in the 19th century, say in terms of Smith normal forms for matrices or the elementary divisor theory. On the other hand, Pr\"ufer introduced his eponymous rings in the 1920s and 1930s \cite{PruferRing}; now-a-days we think of such rings as not necessarily finitely generated analogs of Dedekind domains.  Even in 1943, it was open question \cite{Helmer} whether or to what extent one could salvage the elementary divisor theory (a.k.a. structure theory for finitely generated modules) in this context.\footnote{It's also worth pointing out that, in this context, the distinction between finitely generated and finitely presented modules was still not standard.  For example, in 1949, Kaplansky  \cite{KaplanskyElemDiv} feels the need to write:
\begin{quote}
	the ability to reduce finite matrices does not carry with it results on all finitely generated modules, but only those whose ``relations'' are also finitely generated.
	\end{quote}}
We will not try to analyze why these questions were posed, but merely observe that they were.

Another distinction that seems to reappear in the literature of the time was between torsion and non-torsion elements, especially in abelian groups.  As a sample ``extension'' of the structure theorem for finite abelian groups (thus finitely generated) to torsion abelian groups (not necessarily finitely generated) is the primary decomposition.  Baer mentions primary decomposition several times in papers that we will reference shortly; Kaplansky refers even in 1954 to the original source of the primary decomposition as being ``lost to antiquity'' \cite[\S 20 p. 73]{KaplanskyIAG}.


In a ring-theoretic direction, we have already mentioned Steinitz's version of the structure theorem for finitely generated modules over Dedekind rings.  While we now consider the theory of abelian groups as part of the theory of modules over commutative rings, the papers in the 1930s indicate that considerable attention was paid to the analysis of abelian groups proper, distinct from more ``abstract'' or general module-theoretic analyses, and the latter were typically carried out only for special classes of rings (which have themselves been named, Dedekind rings, Pr\"ufer rings, later Krull domains.)

With that context in mind, in 1936, Baer published \cite{BaerTorsion}, analyzing so-called {\em mixed} abelian groups, i.e., those containing elements both of finite order and infinite order.  If $A$ is an abelian group, Baer considers the subgroup $F(A)$ of elements of finite order and, in modern terminology the corresponding quotient $A/F(A)$ (he calls this the classgroup).  The first result, appearing as result \cite[(1;1)]{BaerTorsion} states that if $S$ is a subgroup of $A$ that is stable by multiplication by $n$, then $S$ is a direct summand of $A$.  Baer observes in a footnote that ``The statement (1 ;1) is well known. But we prove it here, since, to the authors knowledge, it has not been published before.''  In the same paper, Baer introduces the notion of $p$-complete group and complete group (which we would now call $p$-divisible and divisible respectively).  Note that he definitively does {\em not} use the terminology ``injective''.  

In late 1934, L. Zippin submitted a paper entitled ``Countable torsion groups'' to Annals of Math \cite{Zippintorsion}; which appeared in early 1935.  This paper defines a notion of ``root subgroup'' of a countable torsion subgroup \cite[Definition 2]{Zippintorsion}, which coincides with the notion of a $p$-torsion $p$-divisible group; he makes no mention of the notion of ``root subgroup'' outside the torsion case. Theorem 3.3 of this paper establishes that root subgroups are automatically direct summands. Kaplansky highlights the fact \cite[p. 74]{KaplanskyIAG} that Zippin's arguments make no use of the countability hypothesis. Zippin's paper appears to be completely independent of Baer's treatment, and is likewise not mentioned by Baer.  It also appears not to be mentioned in most subsequent treatments we will discuss.  

The case of torsion divisible groups appears elsewhere, in particular in analyses of duality \`a la Pontryagin.  Pontryagin initially studied in \cite{PontryaginCharacters} the group of continuous characters of a topological abelian group; here a continuous character is a continuous homomorphism to the circle group, for which he writes $\mathbf{K}$.  Among other things \cite[Lemma 1]{PontryaginCharacters} states that any character of a subgroup of a group $G$ can be extended to a character of the whole group.  Near the end of the proof, he notes in passing something related to divisibility ``in the group $\mathbf{K}$ division is always possible, though perhaps not unique.''  Analysis of divisibility in this context appears repeatedly.  In \cite[\S 2]{Whitneydivisible0}, Whitney mentions the notion of an {\em infinitely divisible} group.  In Section 3 of this paper, he remarks that the group of rational numbers modulo $1$ is an example of an infinitely divisible group and no group with finitely many generators can be.  In the Remarks on p. 38, Whitney makes the statement: ``If $A'$ is a subgroup of $A$ with division (i.e., $ma$ in $A'$, $m \neq 0$, implies $a$ in $A'$), and $f$ is a homomorphism of $A'$ into $G$, we may always extend $f$ over $A$.'' As Whitney remarks, this results extends the observation of Pontryagin stated above.  


Whitney revisits these ideas in \cite{Whitneydivisible} (submitted November 1943, published in 1944), which introduces the notion of ``completely divisible group.''  His putative goal here is to analyze universal coefficient theorems for cohomology, using ideas related to Pontryagin duality.  Here, \cite[Lemma 1]{Whitneydivisible} expands on the Remarks from the preceding paper and contains various equivalent characterizations of complete divisibility, among which include the ``lifting'' characterization of injectivity that Eilenberg later isolated.   Once again, this paper does not use the word ``injectivity'' to modify abelian groups, seems not to be mentioned in subsequent work, even in later more expository treatments such as \cite{KaplanskyIAG}.  Moreover, Whitney's title suggests that this work is intended to be the first part of some treatment, the subsequent portions of which never appeared.  I view these independent treatments as support for Baer's idea that results such as these were ``well-known''.\todo{Comment on inevitability}



In \cite{Baer}, which appeared in 1940, Baer extends his previous results on abelian groups to abelian groups with a ring $R$ of operators, i.e., $R$-modules.  In so doing, he introduces a modification of the notion of $p$-completeness for an abelian group: for an ideal $I \subset R$, he introduces a notion of $I$-completeness of an $R$-module.  An $R$-module $M$ is called $I$-complete if for every $m \in M$ and every element $i \in I$, we can find an element $m' \in M$ such that $im' = m$.  Baer then analyzes modules that are $I$-complete for all ideals, calling these {\em complete} modules.   He also focuses attention on the case of modules that are $I$-complete for $I$ ranging among the principal ideals of $R$. With this terminology, \cite[Corollary 2]{Baer} contains a converse to the observation that a complete abelian group is a direct summand of every group in which it is a subgroup, and \cite[Theorem 3]{Baer} states that every $R$-module is a submodule of divisible $R$-module.  Once again, injectivity in the guise introduced by Eilenberg is nowhere to be found, nor is any discussion around extension of homomorphisms.


In a note submitted to the PNAS in May 1948 \cite{MacLaneDuality}, Mac Lane formulated the lifting property characterizing injectivity in another analysis of ``duality'' phenomena, also motivated in part by Pontryagin duality.  He recapitulated these ideas in a talk at the Western Sectional meeting of the AMS in November 1948, and submitted a more complete treatment of these ideas appears in December 1949, appearing finally in 1950 as \cite{MacLaneDualityFull}.  According to \cite[Chapter 26]{MacLaneAuto}, Mac Lane was officially at Harvard through the end of the 1947-48 academic year, but had spent considerable time in Europe that year. still in Europe at this time, and only returned to Chicago in the Fall of 1948.  Mac Lane's paper makes no mention of Baer's work, or the work of Whitney, to which it seems most similar, even though Whitney was a colleague at least when the PNAS note was published.  Mac Lane does make a number of other references.  Among other things, Mac Lane gives a formulation preliminary to the eventual universal property of free groups (existence of a lifting along a surjection).  Furthermore, he remarks on the importance of ``duality'' in his approach to cohomology and homology making reference to joint work, still to appear, with Eilenberg \cite{EMOperatorsII}.  The discussion of \cite{MacLaneDualityFull} makes the view that this work is a precursor to Eilenberg's treatment in \cite{SemCartan195051} seem apparent.  Indeed, \cite[Theorems 1.2 and 1.2']{MacLaneDualityFull} explicitly formulate the computation of extensions of abelian groups in terms that look strikingly like the beginning of a free or injective resolution.

Kaplansky introduced the terminology of divisible modules in \cite[\S 2]{KaplanskyDedekind}, which was submitted to Transaction of the AMS in July 1950.  In \cite[Theorem 6]{KaplanskyDedekind}, Kaplansky shows two things.  First, if $R$ is a Dedekind domain and $D$ is a divisible $R$-module, then given any injective ring homomorphism $S \to M$ and a homomorphism $S \to D$, there is an extension $M \to D$.  Second, if $D$ is a divisible submodule of some $R$-module $M$, then $D$ is a direct summand.  Of course, the first point is precisely the condition of injectivity introduced by Eilenberg, to which Kaplansky makes no mention.  Kaplansky then adds a footnote stating that what he is doing is showing that Baer's more complicated notion of completeness is equivalent to the simpler notion of divisibility.  Kaplansky is thus aware of Baer's work.  However, Kaplansky seems unaware of Eilenberg's treatment of injectivity.

Eilenberg's treatment from \cite[Expose I-III]{SemCartan195051} makes no reference of Kaplansky's work or Baer's work, but in fact makes no reference to any other papers either.  In his proof of the uniqueness of group cohomology, he mentions only\footnote{``La démonstration est basée sur une construction (non donnée ici) qui fait correspondre à tout $\Pi$-groupe $A$ un $\Pi$-groupe $\Pi$-injectif $Q$ contenant $A$.}
\begin{quote}
	The proof is based on a construction (not given here) that associates with every $\Pi$-group $A$, a $\Pi$-injective $\Pi$-group $Q$ containing $A$.
\end{quote}
Eilenberg's formulation does not allude at all to the work of Baer--Kaplansky.  

Mac Lane returned to Chicago in Fall 1948, and Kaplansky had just been hired the preceding year.  Kaplansky was Mac Lane's first student, so one imagines that they would have interacted.  Nevertheless, Kaplansky's paper \cite{KaplanskyDedekind} also makes no reference to Mac Lane's PNAS note \cite{MacLaneDuality}, the subsequent treatment \cite{MacLaneDualityFull}, nor to the Whitney paper, so one wonders whether Kaplansky was in fact aware of these ideas at this time.  Nevertheless, \cite[Theorem 6]{KaplanskyDedekind}, establishes the equivalence, for modules over Dedekind domains, of the notion of divisibility and the lifting characterization of injectivity, once again without ever mentioning injectivity.

Kaplansky's book \cite{KaplanskyIAG}, published in 1954 does mention Mac Lane's ideas of duality in the closing ``guide to the literature'' \cite[p. 97]{KaplanskyIAG}. Here Kaplansky writes:
\begin{quote}
	Another approach is to...concentrate instead on an axiomatic study of duality itself...A typical instance of this duality is afforded by divisible groups and free groups.  A group is divisible if and only if it is a direct summand whenever it is a subgroup: a group is free if and only if it gives rise to a direct summand whenever it is a homomorphic image.
\end{quote}
Once again, the lifting characterization of injectivity without any mention of injectivity...

This all leads one to the idea that there were roughly two parallel evolutions, not quite simultaneous, of the eventual notion of injectivity: one that was ``algebraically'' oriented, and one that was ``topologically'' oriented.  A more complete description would suggest various evolutionary dead ends, like the theory of Zippin or Whitney that appear not have been further developed.  Another piece of support of this parallel development theory is the publication in 1953 of the work of Eckmann and Schopf \cite{EckmannSchopf}.  Here, the authors explicitly mention Baer's paper \cite{Baer} in a footnote clarifying:\footnote{``..es ist dort nicht von Injektivit\"at, sondern einer damit \"aquivalenten Eigenschaft (4.2 (c) unserer Note) die Rede.'' }
\begin{quote}
	...it does not speak of injectivity, but of an equivalent property...
\end{quote}
Eckmann and Schopf make no mention of Whitney, Mac Lane or Kaplansky, though the content of Section 4.2 of their paper clearly overlaps with these treatments.  


\subsubsection*{Axiomatic aesthetic}
The treatment of individual modules in Cartan--Eilenberg's homological algebra was remarked upon by Hochschild in his MathSciNet review.  Indeed, from the point of view of homological algebra, focus was to be placed on {\em relations} between modules rather than modules themselves.  This diminishing of the role of individual modules, in conjunction with the arguments above leads to my conclusions that, at the outset, projective modules were not to be viewed as interesting objects in their own right.  Their definition arose only in axiomatizing the property of free modules that was used in the construction of resolutions that were used to compute (co)homology.  Those resolutions themselves were thought of as auxiliary black boxes with the actual cohomology groups being of primary importance.   The relevant lifting property also was more naturally suggested by category-theoretic considerations rather than simply considerations of abstract algebra and thus reflects a ``modernist'' or perhaps ``minimalist'' axiomatic sensibility.\endnote{One speculation, which I think plays rather well with the discussion from the introduction is that this axiomatic aesthetic, coupled with the place Bourbaki held in the mathematical community of the day, was reason enough for people to start studying any mathematical notion. Thus, people may have started studying projective modules in their own right with no further impetus (examples of this abound in the literature).  A vote in this spirit comes from commentary of H. Behnke on Cartan's description of Bourbaki from \cite{CartanonBourbakiArbeits}.  He writes:\footnote{Ich möchte zu zwei Fragen sprechen, die sich dem in die heutige Entwicklung der Mathematik nicht eingeweihten Zuhörer sicherlich aufdrängen. 1. Wie weit beeinflußt das Werk Bourbakis heute schon das mathematische Leben in Deutschland?}
	\begin{quote}
	I would like to address two questions that will certainly occur to listeners who are not privy to the current development of mathematics. 1. To what extent does Bourbaki's work already influence mathematical life in Germany? 
	\end{quote}
After some general remarks, Behnke writes a testimonial to the influence of Bourbaki in Germany:\footnote{Nun zum Punkt 1: Bourbakis derzeitiger Einfluß in Deutschland. Dieser Einfluß ist viel größer, als er vom Herrn Vortragenden aus Bescheidenheit dargestellt wurde. Dieser Einfluß begann sich vor gut 5 Jahren zunächst in der Topologie bemerkbar zu machen. In der Topologie haben wir seit 1935 das glänzende Standardwerk von P. Alexandroff und Heinz Hopf.  Wir hatten uns gerade daran gewöhnt, auf die dort eingeführten Begriffe uns zu stützen, da kam Bourbaki mit seiner neuen Topologie, mit seinem anderen Begriff der Kompaktheit, mit seinen Filtern und seinen uniformen Strukturen. Die jüngste Generation begann immer mehr, sich auf diesen Aufbau der Topologie zu stützen, und jetzt ist es schon beinahe selbstverständlich geworden, daß man sich in der Topologie der Begriffe und Terminologien von Bourbaki bedient. Die Vorlesungen in allen Disziplinen - wenigstens bei uns in Münster - werden schon davon beeinflußt. Es ist z. B. schon deutlich im 2. Teil der Anfängervorlesung zur Infinitesimalrechnung zu spüren. Aber je höher die Vorlesungen und je jünger die Dozenten sind, um so mehr macht sich dieser Einfluß bemerkbar. Besonders deutlich ist er auch bei unseren letzten Dissertationen. Gerade jetzt sind bei uns Bourbakis uniforme Strukturen mit großem Erfolg in der Forschung benutzt worden.
	
	Bei den Studenten wird natürlich gar nicht für Bourbaki geworben. Und doch werden von unseren wirklich guten Studenten die Bücher von Bourbaki geradezu ``gefressen''. Es gehört bei ihnen zum guten Ton, Bourbaki zu lesen und sich seiner Terminologie rücksichtslos zu bedienen.}
\begin{quote}
	Now to point 1: Bourbaki's current influence in Germany. This influence is much greater than the lecturer modestly presented. This influence first began to make itself felt in topology a good five years ago. In topology we have had the brilliant standard work by P. Alexandroff and Heinz Hopf since 1935. We had just gotten used to relying on the concepts introduced there when Bourbaki came along with his new topology, with his different concept of compactness, with his filters and his uniform structures. The youngest generation began to rely more and more on this structure of topology, and now it has almost become a matter of course that Bourbaki's concepts and terminology are used in topology. The lectures in all disciplines - at least here in Münster - are already influenced by this. For example, it is already clearly noticeable in the second part of the introductory lecture on infinitesimal calculus. But the higher the level of the lectures and the younger the lecturers, the more noticeable this influence becomes. It is particularly evident in our latest dissertations. Bourbaki's uniform structures have been used with great success in research.
	
	Of course, Bourbaki is not promoted to students.  And yet our really good students practically ``devour'' Bourbaki's books. It is considered good form for them to read Bourbaki and to use his terminology without hesitation.
\end{quote}
Behnke continues about the participation of students in the Bourbaki seminar, but Behnke was an early collaborator of Cartan, and Behnke's students included Hirzebruch, Grauer, Remmert, and Stein, who would all further develop Cartan's complex analytic ideas and the theory of coherent analytic sheaves.  Certainly at Princeton and the Institute for Advanced Study, granted the presence of Artin, Borel and Weil by the late 1950s, that the influence of Bourbaki here among students was comparable, even if more complicated.} Coupled with the lack of examples, how did projective modules leave this background role? 

Having clouded the waters as regards the linear narrative of injectivity, let us indicate some further complicating factors echoing Hochschild's view that homological algebra brought a shift away from study of individual modules and toward relationships between modules.  The Baer--Kaplansky notion of divisibility could be viewed as a focus on divisible modules in their own right.  Indeed, there is a fascinating ``structure'' theory for divisible abelian groups\footnote{A divisible abelian group is a direct sum of copies of the rational integers and Pr\"ufer $p$-groups (a.k.a., quasi-cyclic groups); \cite[Theorem 4]{KaplanskyIAG}.}, which Kaplansky is explicitly developing in his paper \cite{KaplanskyDedekind}.  The transition from divisibility to injectivity that we have described above seems to coincide with a transition away from establishing structural results for divisible modules.  Instead, the many mentions of injective modules that begin to appear only seem to care about cohomological triviality of injective modules, rather than any explicit description of such modules and there appear to be few attempts in the literature to study divisibility in its own right.\footnote{But not none!  \aravind{Add references?}}

Given this development of the theory of injective modules, I speculate that there has to be a reason that projective modules developed differently: instead of fading away as an obscure technical tool, projective modules instead become a focal point of mathematical research.  One can see shadows of a development of this form in modular representation theory of finite groups today where the projective modules are precisely the {\em uninteresting} modules: one forms the stable module category by essentially working ``up to projectives'' and thereby systematically ignoring their contributions.  Thus, the change in the theory of projective modules arose from some additional, external ``vitalizing'' influence.

\subsection{Fiber bundles and projective modules}
\label{ss:fiberbundlesandprojectivemodules}
\addtoendnotes{\noindent {\bf Fiber bundles and..} \vskip .2em}
The 1951-52 version of the S\'eminaire Cartan was devoted to the theory of analytic functions of several variables \cite{SemCartan195152}.  One point that we have not addressed earlier is how widely available the polycopies of the Cartan seminar were.  Wide publication of some of these seminars had to wait until 1955 or so, as the publication dates indicate.  The last three exposes of this seminar are particularly important for our discussion as Cartan's presentations develop aspects of the theory of coherent sheaves on Stein spaces, in particular formulating results that eventually became ``Theorems A and B''.  These results of Cartan allowed explicit links between the Cousin problems, and sheaf cohomology, and the final Expose of this year is given by Serre who explains how to use Cartan's results in applications, especially to Cousin-type problems.  

This edition of the Seminar gives a first public development of the vision we followed in the Cartan--Weil correspondence around fiber bundles in complex analytic settings (see Paragraph~\ref{paragraph:WeilFB}).  Rather than using fiber bundles directly, however, Cartan has used cohomological formulations involving his theory of coherent sheaves (as discussed in Section~\ref{ss:sheaves}).  That these results from the Cartan seminar were more widely available is attested by other public treatments of this work, e.g., Cartan's expose at a Conference in Brussels \cite{CartanBruxelles} and Serre's subsequent treatment \cite{SerreBruxelles}.  

Serre's initial expose in the Cartan seminar was given on the 16th of June, but the write-up was revised in November of 1952.  Already at this stage, a number of the proofs, especially around applications of Cartan's theorems to theory of analytic fiber bundles are modeled on Weil's treatment in his 1952 course on the algebraic theory of fiber bundles.  This talk already contains results like the characterization of Stein spaces in terms of vanishing of higher cohomology of ideal sheaves. 

In May 1953, Serre gave a talk in S\'eminaire Bourbaki entitled ``\'Espaces fibr\'es algebriques'' \cite{SerreFSBBKI}.  In this talk, Serre discusses Weil's theory of fiber bundles in algebraic geometry, largely using the language of Weil's foundations of algebraic geometry \cite{WeilFoundations} and building on the work of \cite{WeilFBCourse} and the talk upon which that work was based.  In Section 3, he observes that the ``most important'' example of algebraic fiber bundles are the algebraic vector bundles.  In particular, guided by the discussion in the Cartan seminar, one sees Serre's interests moving away from Cartan's analytic treatment into corresponding algebraic situations.

These developments, especially the theory of analytic sheaves, were gaining broad acceptance at the same time, especially around Princeton.  Indeed, Kodaira and Spencer were using Cartan's sheaf theory for complex manifolds more generally, this time in conjunction with analytic techniques.  In 1953, they wrote \cite{KodairaSheaf}, \cite{KodairaSpencer} and \cite{KodairaSpencerII} which adopt the techniques and terminology of Cartan's version of sheaf theory.    

In this section, I aim to explore how projective modules made the transition from ``technical'' actors to objects of ``central'' mathematical focus.  To illustrate this transition, I give a ``rational reconstruction'' of the bridge between projective modules and fiber bundles, and then compare this to the actual presentation, which appears in Jean-Pierre Serre's influential treatise {\em Faisceaux alg\'ebriques coherents} \cite{SerreFAC}, submitted to Annals of Mathematics in October 1954, and appearing in March 1955.  Since Serre appears as a key player at this stage, we highlight something of his mathematical status by inclusion of some discussion around the 1954 Fields Medals, especially as it relates to the place of sheaf theory and fiber bundles in the mathematical mainstream.

\subsubsection*{Projective modules as afterthought}
The link between projective modules and fiber bundles is spelled out in \cite{SerreFAC}.  I will momentarily engage with Serre's proof of this fact, but I want to begin instead with another rational reconstruction of the link.  The link between projective modules and algebraic vector bundles requires first observing that geometric notion of algebraic vector bundles can be described in purely ring/module-theoretic terms; this requires admitting an algebro-geometric point of view; the corresponding module-theoretic notion is that of a ``locally free'' module.  

\begin{entry}[Interlude: cultural features of proof]
	\label{par:projectivelocallyfree}
	I would like to give a proof that finitely generated projective modules over a commutative unital ring $R$ are precisely locally free modules.  This proof is, I believe, a standard proof now, and I think from a modern point of view it might even appear ``organic''; I would like to indicate the cultural assumptions that are required to make that assertion.  I do this both to situate the ideas to a modern reader, and also to emphasize that the idea this proof is more elementary depends itself on various changes in perspective and mathematical values. 
	
	There are two implications to establish.  For the forward implication, we need to know that finitely generated projective modules are locally free.  In other words, given a projective module $P$, we want to show that for any prime ideal ${\mathfrak p} \subset R$, the localized module $P_{{\mathfrak p}}$ is itself free.  This amounts to establishing that finitely generated projective modules over a local ring are free.  
	
	Nowadays, the conceptual scaffolding of this part of the argument is built with an appeal to ``Nakayama's Lemma''.\endnote{According to Matsumura \cite[p. 8]{Matsumura}, the result which I have referred to as ``Nakayama's lemma'' was attributed by Nakayama himself to W. Krull and G. Azumaya.  The proof given there uses a ``determinant trick'', which seems exactly the kind of argument that Bourbaki would have avoided.  Nakayama established the result in \cite[I]{Nakayama}, where it is stated as a generalization of a result of Azumaya from \cite[Theorem 5]{Azumaya}, and no mention is made of Krull's results.}  To be precise, let $S$ be a local ring;  write ${\mathfrak m}$ for the maximal ideal of this local ring and set $\kappa = S/{\mathfrak m}$, usually called the residue field.  If $P$ is a finitely generated projective $S$-module, then we get a (finitely generated) $\kappa$-module $P \tensor_{S} \kappa$, which is a finite-dimensional vector space.  In this context, ``Nakayama's lemma'' is the assertion that a basis for this finite-dimensional vector space can be lifted to provide a basis for $P$.  This class of results was published in the early 1950s.
	
	For the reverse implication, suppose we admit that someone is familiar with the properties of tensor products of modules and localization in commutative algebra in the following form:
	\begin{itemize}[noitemsep,topsep=1pt]
		\item[(1)] over any commutative ring $R$ the assertion that $M$ is projective is equivalent to the condition that the functor $\hom_R(M,-)$ preserves surjections, 
		\item[(2)] an $R$-module homomorphism $N_1 \to N_2$ is surjective if and only if it remains so after localization at all prime ideals ${\mathfrak p}$ of $R$, and 
		\item[(3)] the construction $\hom_R(-,-)$ behaves functorially under localization.
	\end{itemize}  
	The first statement is a reformulation of the ``lifting characterization'' of projective modules stated at the beginning of Section~\ref{s:projective}, but depends on reformulating this property in terms of a functor, i.e., admitting a category-theoretic point of view.  The second statement perhaps becomes ``natural'' only when one adopts a geometric point of view on rings and modules: elements of a module are ``generalized functions'' on an associated geometric space (the prime or maximal ideal spectrum); this algebro-geometric point of view arguably only gained currency due to Weil and Serre's work.  Finally, admitting the last statement also requires some category-theoretic familiarity: for example, one could imagine thinking about this by viewing localization in terms of tensor products and appealing to universal properties. 
	
	Now, suppose $M$ is a locally free module and we want to check that $M$ is projective.  Appealing to (1) we will show that $\hom(M,-)$ preserves surjections.  Localizing at a prime ideal ${\mathfrak p}$, properties of tensor products imply that 
	\[
	\hom_R(M,N) \tensor_R R_{{\mathfrak p}} \cong \hom_{R_{\mathfrak p}}(M_{\mathfrak p},N_{\mathfrak p}).
	\]   
	In that case, since $M_{\mathfrak p}$ is free by assumption, we conclude that $\hom_{R_{\mathfrak p}}(M_{\mathfrak p},-)$ preserves surjections and we conclude.  This rational reconstruction could have been constructed by someone in the Cartan--Eilenberg circle, but doing so requires, I think different primary notions.  It is plausible that this proof could have been constructed by someone in the Cartan--Eilenberg circle, but doing so requires, I think, different primal notions.  Granted what we have to accept, in what sense is this argument more elementary?   
\end{entry}

So how does this presentation compare to the original discussion?  The establishment of the two implications in the statements of \ref{par:projectivelocallyfree} appear separately, and we discuss each such implication separately.

\begin{entry}[Projective implies locally free]
The assertion that projective modules over local rings is contained in \cite[Theorem VIII.6.1']{CE}.  They give a ``cohomological'' criterion for freeness of a module involving the $\operatorname{Tor}$-functor.  More precisely, Cartan--Eilenberg prove that if $M$ is a finitely generated module over a local ring $R$ with maximal ideal ${\mathfrak m}$ and residue field $\kappa := R/{\mathfrak m}$, then whenever $\operatorname{Tor}_1^R(M,\kappa) = 0$, $M$ is necessarily free.  Of course, this vanishing holds when $M$ is a projective $R$-module and thus one concludes that projective modules over local rings are free.\footnote{I highlight the fact that Cartan--Eilenberg do not assume $M$ is projective in the statement of this result, perhaps because for them, local rings need not be commutative.  Thus, on a first reading, projective modules are simply examples of modules that satisfy the hypotheses of the theorem. In fact, from a psychological point of view, it seems to me that the focus on vanishing of $\operatorname{Tor}$ draws one's attention away from particular modules because it implicitly depends on the ring $R$.  Nevertheless, once the theorem is phrased in this way, it immediately raises the question: are there other modules that satisfy this vanishing condition?  As it turns out, the answer to this question is ``no'': Serre established later \cite[Proposition 1]{SerreSplitting} that the vanishing condition in \cite[Theorem VIII.6.1']{CE} is equivalent to $M$ being projective.}


The proof of \cite[Theorem VIII.6.1']{CE} is ``immediate'' from two results in \cite[VIII.5]{CE}.  These auxiliary results are also phrased in terms of ``lifting generators'' of $M \tensor_R \kappa$ to $M$, but the discussion of \cite[VIII.5]{CE} uses additional notation and terminology (that of ``proper, faithful and allowable'') that a dutiful reader would have to process.  No reference is made to ``Nakayama's lemma'', perhaps because such results were not packaged as standard yet.  Moreover, these auxiliary results make no explicit mention of projective modules.  Only one implication of the claim in Paragraph~\ref{par:projectivelocallyfree} appears in \cite{CE}.  To my knowledge, the first published mention of ``Nakayama's lemma'' in this context appears in the proof of \cite[Proposition 1]{SerreSplitting}.
\end{entry}   

What about the reverse implication?  Serre begins \cite{SerreFAC} with the statement:\footnote{On sait que les m\'ethodes cohomologiques, et particulierement la th\'eorie des faisceaux, jouent un role croissant, non seulement en th\'eorie des fonctions de plusieurs variables complexes (cf. [5]), mais aussi en g\'eom\'etrie alg\'ebrique classique (qu'il me suffise de citer les travaux r\'ecents de Kodaira-Spencer sur le th\'eoreme de Riemann-Roch). Le caractere alg\'ebrique de ces m\'ethodes laissait penser qu'il \'etait possible de les appliquer \'egalement a la g\'eom\'etrie alg\'ebrique abstraite; le but du pr\'esent m\'emoire est de montrer que tel est bien le cas.}
\begin{quote}
	We know that cohomological methods, in particular sheaf theory, play an increasing role not only in the theory of several complex variables ([5]), but also in classical algebraic geometry (let me recall the recent works of Kodaira-Spencer on the Riemann-Roch theorem). The algebraic character of these methods suggested that it is possible to apply them also to abstract algebraic geometry; the aim of this paper is to demonstrate that this is indeed the case.
\end{quote}
The reference to $[5]$ here is to Cartan's theory of coherent sheaves on analytic spaces.  Certainly holomorphic vector bundles were known examples of coherent {\em analytic} sheaves.  

Serre observes early in \cite[\S 41 Proposition 5]{SerreFAC} that algebraic vector bundles also give examples of coherent algebraic sheaves.  Here, Serre makes explicit reference to Weil's work on algebraic fiber bundles and, in particular, the Chicago lecture notes on algebraic fiber spaces \cite{WeilFBCourse}.

\begin{entry}[Locally free implies projective]
Projective modules enter only in \cite[\S 50]{SerreFAC}.  Proposition 4 in this section contains the assertion that projective modules are {\em precisely} the locally free modules.  Let us compare Serre's proof of this result to that in Paragraph~\ref{par:projectivelocallyfree}. We have already discussed the projective implies locally free direction of the statement.  The reverse implication again relies on results in Cartan--Eilenberg: in particular, the reader is required to know about notions of projective dimension of modules and a characterization of projective modules as modules of projective dimension $0$, for which the reader is referred to \cite[VI.2]{CE}.  To establish the relevant dimension estimate, Serre relies on an exercise in Cartan--Eilenberg: \cite[Chapter VII Exercise 11]{CE} which bounds projective dimension in terms ``local'' projective dimension.  These results appear more than 150 pages into the Cartan--Eilenberg text.  

In private communication, Serre remarked that once he saw that projective modules were locally free, he saw the converse ``essentially immediately''; but the proof given makes the dependence of the phrase ``essentially immediately'' seem to depend on detailed familiarity with the ideas of Cartan--Eilenbeg.  The proof Serre gives makes free use of the results in Cartan--Eilenberg, which I emphasize remained unpublished from submission through publication of Serre's paper.  Anyone not part of the circle consisting of people with detailed familiarity of the Cartan--Eilenberg text could not plausibly have read and understood this proof, and it seems a stretch to describe the proof at that time as elementary as relied on what can only be described as cutting-edge mathematical ideas.  Moreover, I imagine processing this proof without some kind of roadmap to the Cartan--Eilenberg text would be difficult even for someone to whom that text {\em was} available.  As such, the cultural context of Serre's proof is paramount: it is clearly written for people with particular familiarity, and emphatically not a familiarity that most mathematicians of the time could have had.\endnote{This instance is by no means isolated: one could analyze many other proofs in \cite{SerreFAC} similarly.  Let me give two examples to illustrate this.
\begin{itemize}[noitemsep,topsep=1pt]
	\item Serre's use of sheaves mentions Leray's published treatments, but granted what we have said about the reception of Leray's treatment of sheaves (see Section~\ref{ss:sheaves}) it seems hard to believe that those not familiar with Cartan's revision would be able to process this discussion.  Once again, Cartan's revised view of sheaves was only available in polycopies of exposes from the Cartan seminar whose distribution we also discussed.
	\item Serre makes free use of Cartan and Eilenberg's homological algebra in his treatment of sheaf cohomology, in particular, he makes implicit use of the spectral sequence of a double complex; once again, treatments of these ideas are perhaps most easily found in the Cartan seminar.
	\end{itemize}
It seems impossible for me to imagine the audience of Serre's text as anything but a small group of mathematicians.}
\end{entry}
  
Immediately thereafter, Serre remarks that projective modules are equivalent to algebraic vector bundles.  In fact, \cite{SerreFAC} is replete with references to Weil's theory of fiber spaces in algebraic geometry, here in the form of the Chicago lecture notes, so we conclude that Serre has studied this paper carefully.\footnote{Further, in January 1954, Weil remarks to Cartan \cite[p. 348]{CartanWeil}: ``J'ai reçu une lettre de Serre il y a quelque temps; je songeais vaguement \`a y r\'epondre un jour, mais je n'arrive plus \`a remettre la main dessus, et d'ailleurs je n'avais rien d'int\'eressant \`a repondre.  Il faisait des objections \`a quelques passages des notes sur les vari\'etes fibr\'ees...''}  Moreover, \cite[\S 50 Corollaire]{SerreFAC} contains all aspects of the dictionary: projective modules, locally free sheaves of modules and algebraic vector bundles.  

All of a sudden, it seems everything has changed. Bott characterized the ideas of sheaf theory as part of the ``mathematical mainstream'' (see Section~\ref{ss:sheaves}), and we will see additional evidence for this belief shortly.  By means of the identification Serre has provided, projective modules are now part of this conversation also: they have been contextualized in a world mixing homological algebra, algebraic fiber bundles, and sheaf theory.  Among the algebraic fiber bundles, we have already seen the algebraic vector bundles described as ``the most important examples''.  The full force of the Cartan--Weil vision can thus be harnessed to motivate the study of projective modules as rich {\em in its own right}, at least to those who were exposed to this vision.   It will take a bit of time for this change in perspective to percolate, but the ball has been set rolling.


\subsubsection*{From analogies to concrete problems}
Famously, Weil addressed the guiding role of analogy in mathematics \cite[p. 408]{WeilII}:\footnote{Rien n'est plus f\'econd, tous les mathematiciens le savent, que ces obscures analogies, ces troubles reflets d'une th\'eorie \`a une autre, ces furtives caresses, ces brouilleries  inexplicables; rien aussi ne donne plus de plaisir au chercheur. Un jour vient o\`u l'illusion se dissipe; le pressentiment se change en certitude; les th\'eories jumelles r\'ev\`elent leur source commune avant de dispara\^itre; comme l'enseigne la Gita on atteint \`a la connaissance et \`a l'indiff\'erence en m\^eme temps.}
\begin{quote}
	As all mathematicians know, nothing is more fruitful than these obscure analogies, these troubled reflections from one theory to another, these furtive caresses, these inexplicable confusions; also, nothing gives more pleasure to the seeker. A day comes when the illusion dissipates; presentiment changes into certainty; twin theories reveal their common source before disappearing; as the Gita teaches we achieve knowledge and indifference at the same time.
\end{quote}

Projective modules have been identified as algebro-geometric analogs of vector bundles in topology.  Given the primacy of fiber bundles in Weil's worldview, which has evidently been passed to this generation of Cartan's students, it is natural to wonder how this analogy will be used: at which stage in the evolutionary arc of this analogy can the identification of algebraic vector bundles with projective modules be located?  Is the analogy still to be viewed as obscure, or have we reached a stage of knowledge and simultaneous indifference?   

Two questions are posed at this stage, and to address the preceding question we will track their reception in the mainstream mathematical conscience.   In \cite[\S 4(3)]{SerreFSBBKI}, Serre considers the following situation: if $V$ is a non-singular complex variety, and $G$ is an algebraic group, we can consider the comparison map from isomorphism classes of algebraic fiber bundles with structure group $G$ to the corresponding set of isomorphism classes of analytic fiber bundles.  In this case, Serre writes: it appears very likely that the comparison map is a bijection.  

Serre himself would go on to analyze this question in \cite{SerreGAGA} with the additional hypotheses that $V$ be a {\em complete} algebraic variety.  It is unclear what evidence existed for this conjecture at this stage.  The techniques of the day were sufficient to establish that this question of Serre had a negative answer if $V$ is a any non-singular complex affine curve of genus $g > 0$: the divisor class group of a non-singular affine curve was well-known to be non-trivial, but Cartan's theorems show that all analytic line bundles on a non-singular complex affine curve are trivial for cohomological reasons (Cartan's formulation of the so-called Oka principle).  

In \cite{SerreFAC}, Serre also makes the comment \cite[p. 243]{SerreFAC}:
\begin{quote}
	Signalons que, lorsque $V = K'$ (auquel cas $A = K[X_1,\ldots, X_r]$), on ignore s'il existe des $A$-modules projectifs de type fini qui ne soient pas libres, ou, ce qui revient au m~me, s'il existe des espaces fibr\'es alg\'ebriques \`a fibr\'es vectorielles, de base $K^r$, et non triviaux.
\end{quote}
Serre raises this question, no doubt motivated by his training as a topologist, by the analogy between algebraic and topological vector bundles.  The belief that all projective modules over polynomial rings over a field are trivial has sometimes been recast in the literature as the ``Serre conjecture'', but Serre has very adamantly asserted that he did not have an opinion about the existence or non-existence of non-free projective modules in this situation: it was only a problem \cite[Footnote p. 1]{LamSerreProblem}.  
 
There are two cases where one can show that projective modules over polynomial rings are trivial, with the known technology at the time.  The first is the case where $r = 1$, since in this case, the result follows immediately from th structure theorem for finitely generated modules over a principal ideal domain.  The second case requires a bit more speculation at this stage.  It seems likely known to Krull that in a unique factorization domain, all height $1$ prime ideals are principal.  That the class group of a unique factorization domain is trivial then follows immediately.  Thus, granted the dictionary that had become available between the divisor class group, line bundles and rank $1$ projective modules, we conclude that rank $1$ projective modules over a unique factorization domain are necessarily trivial.  

\subsubsection*{Serre and the 1954 Fields medal}
One fundamental difference in reception between the  two questions just posed is that, in the meantime, Serre and Kodaira win the 1954 Fields medals.  Ignoring the question of the cachet of the Fields Medal at this point in history, what this award does do is solidify the mathematics of the winners as {\em mainstream} mathematics.  Attention is thus focused on the style of work that I've been discussing: topology, sheaf theory, complex analysis and its interactions with algebraic geometry.  

The 1954 Fields medal committee was presided over by Hermann Weyl, with other members including E Bompiani, F Bureau, H Cartan, A Ostrowski, A Pleijel, G Szegö, and E C Titchmarsh.  Weyl begins his address to the conference with the statement: ``That at each International Mathematical Congress two gold medals be presented to two young mathematicians who have won distinction in recent years by outstanding work in our science.''  Weyl continued to describe the work of the 1954 winners Kodaira and Serre.  
\begin{quote}
	...for I realise how difficult it is for a man of my age to keep abreast of the rapid development in methods, problems and results which the young generation forces upon our old science; and without the help of friends inside and outside the Committee I could not have shouldered this burden at all. It rests more heavily on my than on my predecessors' shoulders; for while they reported on things within the circle of classical analysis, where every mathematician is at home, I must speak on achievements that have a less familiar conceptual basis. A report like this cannot help reflecting personal impressions.
\end{quote}
Having thus characterized the realm of ideas in which Kodaira and Serre are working, Weyl then mentions ``I find it convenient to explain as briefly as I can a number of universal concepts before entering upon some of our laureates' individual achievements. Be prepared then to have to listen now to a short lecture on cohomology, linear differential forms, faisceaux or sheaves Kahler manifolds and complex line bundles.''  Weyl continues to describe the work of Kodaira and Serre mentioning briefly the recently discovered links between the work of the pair; in this direction he remarks ``you may get the impression that our Committee did wrong in awarding the Fields Medals to two men whose research runs on such closely neighboring lines. This contact, however, has been established only during the last year and may well be a transient phenomenon.''  He closes with a strikingly positive description of Kodaira and Serre, writing about the latter:
\begin{quote}
	let me say this that never before have I witnessed such a brilliant ascension of a star in the mathematical sky as yours.
\end{quote}

Weyl's description above contributed to the shaping of Serre's reception by the broader mathematical community.  Weyl alludes above to ``help'' in preparing his presentation of the mathematics of Kodaira and Serre, and it perhaps comes as no surprise, given Weyl's proximity to Weil at the IAS that Weil was one of the people whom Weyl consulted.  

In a letter to Cartan dated Nov 8, 1953, Weil writes \cite[p. 338 ]{CartanWeil}:\footnote{``L'oncle Hermann (par lettre) et Bureau (de vive voix) m'ont consult\'e sur la m\'edaille Fields.  Il me parait ridicule et injuste de pousser Serre aux d\'epens de Borel et Thom, et, dans un autre genre, de N\'eron.  Serre est plus brillant, c'est entendu.  Mais je ne vois pas comment on pourrait, \`a l'heure actuelle, affirmer qu'il soit nettement et clairement sup\'erieur aux autres; et, si on n'en est pas s\^ur, il me parait injuste de se donner l'air de l'etre.  J'entends bien qu'il y a l\`a des questions de tactique, et qu'il y a plus de chances de gagner si on pr\'esente un candidate que si on en pr\'esente plusieurs.  Mais ceci n'est pas une partie de football, et apr\`es tout il importe peu que l'un ou l'autre ait la m\'edaille en chocolat.''}
\begin{quote}
Uncle Hermann (by letter) and Bureau (orally) consulted me on the Fields Medal. It seems to me ridiculous and unfair to push Serre at the expense of Borel and Thom, and, in another direction, N\'eron. Serre is more brilliant, that is understood. But I do not see how one could, at the present time, affirm that he is well and truly superior to the others; and, if one is not sure, it seems to me unfair to give the appearance of being so. I understand well that there are questions of tactics here, and that there is greater chance of winning if one presents a single candidate than if one presents several. But this is not a game of football, and after all it matters little whether one or the other has the chocolate medal.
\end{quote}

Cartan responds to this letter on the 16th of November \cite[p. 341]{CartanWeil}:\footnote{Le passage de ta lettre qui concerne les m\'edailles Fields m'a supris et contrari\'e.  Je ne m'\'etais jamais imagin\'e que tu mettais N\'eron et Thom sur le m\^eme niveau que Serre; \`a la rigueur, on pourrait discuter pour Borel, quoique ce ne soit pas mon avis.  En tout as, tu me connais assez pour ne pas mettre en doute que si je n'avais pas la certitude que Serre est \guillemetleft nettement at clairement sup\'erieur aux autres\guillemetright, je ne me donnerais pas l'air d'en \^etre s\^ur.  Il ne s'agie pas de tactique ici.  Il est curieux de constater (ce n'est pas la premi\`ere fois que je fais cette observation) que les personnalit\'es dont l'envergure est indiscutable sont parfois les plus discut\'ees.  }
\begin{quote}
The passage in your letter concerning the Fields medals surprised and upset me. I had never imagined that you would put N\'eron and Thom on the same level as Serre; maybe, we could discuss Borel, although that is not my opinion. In any case, you know me well enough not to doubt that if I were not certain that Serre is ``well and truly superior to the others'', I would not give the impression of being sure of it. This is not a question of tactics. It is curious to note (it is not the first time that I have made this observation) that the personalities whose stature is indisputable are sometimes the most discussed.
\end{quote}

Weil responds in a latter on November 20th \cite[p. 342]{CartanWeil}:\footnote{``J'ai sans doute en tort de penser que c'\'etait pour des raisons tactiques que tu concentrais tes efforts sur Serre en ce qui concerne la m\'edaille Fields; mais en tout cas je ne puis qu'admirer tes \guillemetleft certitudes \guillemetright, tout en trouant ridicule ta phrase qui dit \guillemetleft ce sont les personnalit\'es les plus indiscutables qui sont le plus disut\'ees\guillemetright.  Je ne vois nullement en quoi je \guillemetleft discute\guillemetright les m\'erites (incontestables, en effet) de Serre en leur comparant les m\'erites de Borel, Thom et N\'eron, qui, \`a moi, me paraissent non moins indiscutables.  Je suis particuli\`erement surpris que to te figures pouvoir avoir une opinion fond\'ee sur les m\'erites de N\'eron--je ne te savais pas si comp\'etent en g\'eometrie alg\'ebrique et en arithm\'etique; je puis t'affirmer, quant \`a moi, que ce qu'il a fait est toute premi\`ere importance, et tout \`a fait comparable en difficult\'e \`a ce qu'a fait Serre.  Pour Borel and Thom, je ne conteste pas tu comp\'etence, mais je te conseille un effort pour te d\'egager un peu de tous go\^uts personnels.''}
\begin{quote}
	I am probably wrong to think that it was for tactical reasons that you concentrated your efforts on Serre with regard to the Fields Medal; but in any case I can only admire your ``certainties'', while finding ridiculous your sentence which says ``it is the most indisputable personalities who are the most disputed''. I do not see at all in what way I ``dispute'' the merits (uncontestable, in fact) of Serre by comparing them to the merits of Borel, Thom and N\'eron, which, to me, seem no less indisputable. I am particularly surprised that you think you can have a well-founded opinion on the merits of N\'eron--I did not know you were so competent in algebraic geometry and arithmetic; I can tell you, for my part, that what he has done is of the utmost importance, and quite comparable in difficulty to what Serre has done. As for Borel and Thom, I do not contest your competence, but I advise you to make an effort to free yourself a little from all personal tastes.
\end{quote}
Weil continues with prescient:\footnote{La verit\'e est qu'il n'existe pas d'instrument de mesure qui, \`a l'heure qu'il est permette de dire avec certitude que l'un d'eux soit sup\`eriur \`a l'autre; tout ce qu'on peut en dire est du domaine du pari; dans dix ans ou dans vingt ans, il est possible qu'on puisse se faire un jugement clair sur leurs m\'erites respectifs (il est non moins possible qu'on ne le puisse pas non plus dans 10 ans ni dans 20 ans, si d'ici l\`a ils se maintiennent toujours \'egaux \`a eux-m\^emes, ce qui est \'evidemment \`a souhaiter.)} ``The truth is that there is no measuring instrument that at the moment allows us to say with certainty that one of them is superior to the other; all that can be said about it is a matter of guesswork; in ten years or twenty years, it is possible that we will be able to make a clear judgment on their respective merits (it is no less possible that we will not be able to do so in 10 years or in 20 years, if by then they remain equal to themselves, which is obviously to be hoped for.)''  

In the meantime, Weil and Cartan discuss a number of issues of academic politics around hiring.  Weil writes \cite[p. 347]{CartanWeil}:\footnote{``A ce propos, ta lettre \`a Dieudonn\'e indique que Delsarte et toi \^etes en train de vous adonner au jeu des combinaisons \`a multiple d\'etente; \`a mon avis vous \^etes \`a ce jeu-l\`a des amateurs, et des amateurs maladroits, et vous serez toujours battus par les professionnels; il est donc tout \`a fait inutile de perdre votre temps \`a ce genre de calculs, et illusoire de b\^atir l\`a-dessus des r\`egles d'action.''}
\begin{quote}
	In this regard, your letter to Dieudonné indicates that you and Delsarte are indulging in strategic gamesmanship; in my opinion you are amateurs in this game, and clumsy amateurs, and you will always be beaten by the professionals; it is therefore completely useless to waste your time on this kind of calculation, and illusory to build rules of action on it.
\end{quote}

The situation regarding the stature of the mathematicians discussed above as regards the 1958 Fields Medal is also fascinating.  For example, Thom and Borel, in addition to Hirzebruch were considered ``leading candidates'' based on correspondence possessed by Zariski \cite[p. 46]{BaranyFM195058}. Other candidates with algebraic/topological interests include, in no particular order, A. Grothendieck, J. Nash, J. Milnor and I. Shafarevich (though I should qualify here that, as the deliberations were taking place in 1957, this was likely more for Grothendieck's work in functional analysis and homological algebra, rather than later work in algebraic geometry).  These candidates and deliberations suggest the important role played by algebraic geometry/topology within the mathematics of the era.



\subsection{On the precipice of a self-sustaining domain}
\label{ss:conclusion}
While Serre's Faisceaux Algebriques Coherents \cite{SerreFAC} was being written in 1954, a host of parallel developments took place that we do not want to ignore.  Some of these notions were mentioned in passing in the preceding sections, in particular the work of Kodaira and Spencer applying techniques of harmonic integrals to the analysis of algebraic varieties \cite{KodairaSheaf,KodairaSpencer,KodairaSpencerII}.  In the preface to \cite{HirzebruchRR}, Hirzebruch reminisces on this period, discussing his interactions with Kodaira and Spencer and drawing a distinction between their approach to the theory of algebraic manifolds and the approach being worked out by Cartan--Serre, even though both groups could be said to be applying ``sheaf-theoretic'' techniques.  Weyl also references a distinction in our summary of his presentation of the 1954 Fields medals.

Hirzebruch used techniques learned from these groups together with work from Thom's cobordism theory to formulate and prove a version of the Riemann-Roch theorem \cite{HirzebruchRRPNAS}; this was announced in late 1953 and appeared in early 1954.  In conjunction with the work of Cartan and Serre on cohomology of coherent sheaves on suitably ``compact'' spaces, one can see an impetus to return to the study of {\em individual} vector bundles, rather than relations between vector bundles.  Indeed, Hirzebruch's formulation of the Riemann--Roch theorem provides a formula for the Euler characteristic of an analytic vector bundle in terms of topological data.

After the submission of \cite{SerreFAC}, the scope of homological techniques and sheaf cohomology continued to expand.  Eilenberg published a series of papers drawing members of the Japanese algebra school, e.g., Nakayama, into the orbit of homological techniques \cite{EilenbergCohFinDim,EilIN,EilN,EilNN}.  The links between homological algebra and commutative algebra mentioned in the introduction to Cartan--Eilenberg were deepened.  Auslander and Buchsbaum announced \cite{AuslanderBuchsbaumPNAS} their homological characterization of regular local rings as those rings with finite global dimension; this proof was later exposed by Serre \cite{SerreAusBuch}.  

The question asked by Serre in \cite{SerreFSBBKI} regarding comparison of algebraic and analytic fiber bundles, one of the two which we discussed in the previous section, was taken up again by Serre in \cite[\S 20]{SerreGAGA}.  He establishes \cite[Proposition 18 p. 31]{SerreGAGA} that every analytic vector bundle on a projective variety admits a unique analytic structure.  It must by this point have been well-known that the map in question was not a bijection even for line bundles over a non-singular affine curve.  Nevertheless, the question of comparison of algebraic and analytic vector bundles in non-projective situations seems to have lain dormant for many years until the work of Griffiths \cite{GriffithsI,GriffithsII} some 15 years later.  

The second question which we highlighted in the preceding section eventually became known as Serre's problem on projective modules, and it is fascinating to speculate about its role in the development of the theory;  a textbook treatment can be found in \cite{LamSerreProblem}.  One has to imagine that the tantalizing ease with which one can state the Serre problem, provided one accepts the notion of a projective module, combined with lack of an approach to the problem led to persistent fascination.  

Grothendieck asks \cite[p. 6]{GrothendieckSerre} in February 1955 (after Serre's FAC was submitted, but before it appeared in March 1955)
\begin{quote}
	is a finitely generated projective module over the ring in question (for instance a ring of polynomials or of holomorphic functions)
	free? Is this easy to see in interesting special cases?  If I understand correctly, in the case of polynomials, it is not even known whether this theorem is true, and one has to restrict oneself to graded rings to get a result.
\end{quote}
This question appears among a long list of other questions/criticisms that Serre answers carefully, but he seemingly skips over responding to this question.\todo{Check references!}

Progress on the Serre problem was slow relative to the pace of the development of sheaf theory, algebraic geometry, fiber bundles etc. during this period.  Serre's analogy between topological and algebraic vector bundles lay dormant for a few years, but all that changed in 1958.  

Seshadri published \cite{SeshadriSerreProb} establishing the first non-trivial case of the Serre problem: the case $r = 2$ in the notation used in the preceding section.  Seshadri's published proof was, according to the written text, a modified version of something he had earlier sent to Serre, incorporating suggestions of the latter.  The idea of the proof itself was once again inductive with respect to dimension: attempting to reduce the $2$-dimensional problem to the already solved $1$-dimensional problem.  This bit of progress, shortly after the initial formulation probably also provided hope that Serre's problem {\em did} admit a positive solution, and thereby perhaps led other people to its consideration.

Shortly thereafter, Serre revisited the theory of projective modules outside of the homological theory of commutative rings.  In \cite{SerreSplitting}, presented in the S\'eminaire Dubreil in May 1958, Serre returns to the identification of vector bundles on affine varieties with projective modules over the coordinate ring of a variety from \cite{SerreFAC}.  If the treatment of the earlier paper made the dictionary seem like an afterthought or supporting detail, the presentation here takes a rather different form.  Here, Serre's presentation seems to take more seriously the idea that projective modules are analogs of vector bundles in topology, and spells out the correspondence in considerably more detail.  Support for this interpretation is given by the main result of this paper, which is rather explicitly motivated by topological analogy: if $R$ is a commutative ring whose maximal ideal spectrum is connected and has dimension $d$, then every projective module of rank $r > d$ splits as the sum of a module of trivial module of rank $r-d$ and a projective module of rank $d$.  The proof of the corresponding fact in topology goes back to the idea of ``general position'' in topology.  One can even argue that Serre's proof is an adaptation of these ideas to algebraic geometry.  

From Serre's point of view, this result encodes much of the intuitive feel of the theory at this point in time.  For example, if $R$ is a local ring, then its maximal ideal spectrum has dimension $0$, so the splitting theorem implies that finitely generated projective $R$-modules are free, as was observed by Cartan--Eilenberg.  Another practical consequence had bearing on the Serre problem itself: one reduces this problem to the analysis of projective modules of rank smaller than the dimension.  This additional partial progress undoubtedly buoyed interest in the problem.

Yet another piece of support for the theory comes from the numerous attempts at classification problems that arise.  Grothendieck analyzed algebraic fiber bundles with structure a complex Lie group over the projective line in \cite{GrothendieckFibres} (submitted in October 1956).  Grothendieck discusses the case of vector bundles first, and it is unclear if he was aware of earlier work that could be interpreted as a solution to the problem.  Indeed, it was retroactively realized that the case of vector bundles is essentially contained in work of Dedekind--Weber as well as work of Birkhoff.  Shortly thereafter, Atiyah \cite{AtiyahVBEll} extends Grothendieck's results by classifying vector bundles over an elliptic curve.  Atiyah makes use of his previous work \cite{AtiyahKS} and phrases his results in terms of a Krull--Schmidt type theorem.

A further impetus to the theory came from Grothendieck's reworking of the Riemann--Roch theorem.  As exposed in Borel--Serre \cite[\S 4 p. 105]{BorelSerreRR}, for any algebraic variety $X$, one defines $K(X)$ as a quotient of the free abelian group on the set of isomorphism classes of coherent sheaves on $X$ modulo the relation that we identify the middle term in a short exact sequence as the sum of the two outer terms.  Two remarks are made immediately: first, a universal property for the group $K(X)$ is observed, in particular making it clear that $K(X)$ is the universal recipient for Euler characteristics of coherent sheaves.  

Second, the authors remark that one can make a similar definition for algebraic vector bundles on $X$ yielding a group $K_1(X)$, and that by construction there is a map $K_1(X) \to K(X)$.  Then, \cite[Th\'eor\`eme 2]{BorelSerreRR} states that if $X$ is a non-singular, irreducible, quasi-projective variety, then the map $K_1(X) \to K(X)$ is an isomorphism.  The proof given by Borel and Serre, effectively still the standard proof of this result, builds on many of the results we have discussed above in an essential way: local-to-global ideas are used to reduce the result to a local computation, where it is encompassed in the homological theory of regular local rings.  

\subsubsection*{Conclusion}
The analysis of vector bundles in general, and the classification problem we described just above lead to a host of new problems for the theory.  These problems, in conjunction with the Serre problem and questions arising from the nascent algebraic K-theory lead to a flurry of work.  All of this classification work is guided by the analogy with topology, and the problems seem structured so that understanding requires new mathematical architecture, rather than simply a recycling of ideas of the now classical bundle theory.  

That mathematicians were willing to engage with the new ideas around projective modules thus seems predicated on their arrival within a suitably receptive and interactive community.  For me this precariousness, by which I mean the mangle of factors we have highlighted that conspired to make the theory of projective modules relevant, including Bourbaki's push of ``axiomatic aesthetics'', fiber bundles as unifying ``important'' mathematical questions, analogy as a source for mathematical inspiration, either in Eilenberg and Mac Lane's mathematical conception of ``naturality'', or in Weil's more prosaic conception bridging mathematical subdisciplines, all serve to make the story more precious.  More strongly, one cannot separate the theory of projective modules from the community in which it emerged. 

I have, essentially arbitrarily, chosen to end the story here, the place which, in so many of my mathematical papers, the story {\em begins}.  As I have left them: projective modules in particular, and vector bundles on algebraic varieties in general are, around the end of the 1950s, at the vanguard of research in algebraic geometry and commutative algebra.  Bolstered by numerous high-status mathematicians routinely engaging with them and re-appearing in different sub-disciplines, the potential energy of the original definition, if one can speak of such a thing, has been converted to the kinetic energy of activity and knowledge production.

I return briefly to the discussion of Section~\ref{ss:whynow} as regards ``Why now?''.  A final, more personal, reason emerged only as I marked a stopping point, the arbitrariness of which reminds me of the self-conscious pronouncement of Mann's narrator in Doktor Faustus ``what in my own conscientious authorial opinion can really lay no claim to such segmentation'' \cite[p. 120]{mann1997doctor}, coincidentally written during the same period on which much of this text focuses.  While I embarked on this study feeling a certain anxiety about the value of the mathematics I do, the process of composition has led me to a rich web of knowledge into which mathematics feeds.  The discovery of that interconnectedness has (oddly?), I think, led me to value the subject more, underlining the fundamentally human nature of the pursuit.

\newpage
\addcontentsline{toc}{section}{Notes}
\theendnotes

\newpage
\begin{scriptsize}
\bibliographystyle{alpha}
\bibliography{ClayNotes}
\end{scriptsize}
\Addresses
\end{document}